\documentclass[]{article}  

\usepackage[margin=30truemm]{geometry}

\usepackage[whole]{bxcjkjatype}

\usepackage[cmex10]{amsmath} %
\usepackage{amssymb}  %
\usepackage{amsthm}

\usepackage{mathtools}%

\usepackage{cases} 

\usepackage{fixltx2e} 

\usepackage[subrefformat=parens]{subcaption} %
\usepackage{graphicx} %
\usepackage{epstopdf}

\usepackage{xcolor}

\usepackage{cite}

\usepackage{algorithm}
\usepackage{algorithmic}

\newtheorem{definition}{Definition}%
\newtheorem{theorem}{Theorem}%
\newtheorem{lemma}{Lemma}%
\newtheorem{proposition}{Proposition}%
\newtheorem{corollary}{Corollary}%
\newtheorem{assumption}{Assumption}%
\newtheorem{remark}{Remark}

\definecolor{myGreen}{rgb}{ 0, 0.7, 0.3 }
\definecolor{myBlue}{rgb}{ 0, 0.4, 1 }
\definecolor{myPurple}{rgb}{ 0.7, 0, 0.3 }
\definecolor{myGray}{gray}{ 0.7 }
\definecolor{myFigGray}{gray}{ 0.5 }

\newcommand{\cB}{\textcolor{myBlue}}

\newcommand{\MyHighlight}[1]{\textbf{#1}}

\usepackage[overlay]{textpos}%

\usepackage{pgfplots}
\pgfplotsset{compat=newest}
\usepackage{tikz}
\usetikzlibrary{positioning}
\usetikzlibrary{arrows}



\begin{document}
	
	\title{Second-Order Sampling-Based Stability Guarantee for Data-Driven Control Systems%
		\thanks{%
			This work has been submitted to the IEEE for possible publication. Copyright may be transferred without notice, after which this version may no longer be accessible.
			This work was partly supported by JSPS KAKENHI Grant Number JP18K04222.
			Parts of this paper have been published in conference proceedings \cite{ItoCDC19}. 
			We would like to thank Editage (www.editage.jp) for the English language editing.
	}} 
	
	\author{Yuji Ito\thanks{Yuji Ito is the corresponding author and with Toyota Central R\&D Labs., Inc., 41-1 Yokomichi, Nagakute-shi, Aichi 480-1192, Japan	(e-mail: ito-yuji@mosk.tytlabs.co.jp).}
		\and
		Kenji Fujimoto\thanks{Kenji Fujimoto is with the Department of Aeronautics and Astronautics, Graduate School of Engineering, Kyoto University, Kyotodaigakukatsura, Nishikyo-ku, Kyoto-shi, Kyoto 615-8540, Japan (e-mail: k.fujimoto@ieee.org).}
	}  
	
	\date{}
	
	\maketitle

\newcommand{\SymColor}[1]{#1}  

\newcommand{\El}[3][]{\SymColor{#1[}{#2}\SymColor{#1]}_{#3}}

\newcommand*{\MyTRANSPO}{\top}

\newcommand*{\wildcard}{\SymColor{\bullet}}

\newcommand{\NotationVec}{\SymColor{\boldsymbol{y}}}
\newcommand{\NotationbVec}{\SymColor{\boldsymbol{z}}}
\newcommand{\NotationMat}{\SymColor{\boldsymbol{Y}}}
\newcommand{\IDNotation}{\SymColor{a}}
\newcommand{\IDbNotation}{\SymColor{b}}
\newcommand{\DimANotation}{\SymColor{c}}
\newcommand{\DimBNotation}{\SymColor{d}}

\newcommand{\MyLabelAlgSR}{\SymColor{1}}

\newcommand{\ASSNLFA}{\SymColor{(H.1)}}
\newcommand{\ASSNLFB}{\SymColor{(H.2)}}
\newcommand{\ASSNLFC}{\SymColor{(H.3)}}
\newcommand{\ASSNLFD}{\SymColor{(H.4)}}

\newcommand{\ASSinNestdSet}{\SymColor{(A.1)}}
\newcommand{\ASSVCparamExist}{\SymColor{(A.2)}}
\newcommand{\ASSboundedMargin}{\SymColor{(A.3)}}
\newcommand{\ASSnonemptySetDist}{\SymColor{(A.4)}}

\newcommand{\LI}[1]{{#1}^{\SymColor{\mathrm{I}}}}

\newcommand{\setB}[1]{\SymColor{\partial}{#1}}
\newcommand{\dinf}[1]{ \SymColor{   \underline{\partial}        (} #1  \SymColor{)}}
\newcommand{\ddinf}[1]{\SymColor{   \underline{\partial^{2}}    (} #1  \SymColor{)}}
\newcommand{\dsup}[1]{ \SymColor{   \overline{\partial}         (} #1  \SymColor{)}}
\newcommand{\ddsup}[1]{\SymColor{   \overline{\partial^{2}}    (} #1  \SymColor{)}}

\newcommand{\dirddinf}[1]{\SymColor{   \underline{\mathrm{d}^{2}}   (} #1  \SymColor{)}}
\newcommand{\dirddsup}[1]{\SymColor{   \overline{\mathrm{d}^{2}}   (} #1  \SymColor{)}}

\newcommand{\zzzrawBLB}[2]{\SymColor{ \max }\{  {\zzzLM{#1}}  {\zzzUB{#2}} ,-{\zzzUM{#1}} {\zzzLB{#2}} \}  }
\newcommand{\zzzrawBUB}[2]{\SymColor{ \max }\{  {\zzzUM{#1}}  {\zzzUB{#2}} ,-{\zzzLM{#1}} {\zzzLB{#2}} \}  }

\newcommand{\ModelSetFreeParam}[1]{\SymColor{\alpha}_{#1}}

\newcommand{\hypNoise}{\SymColor{\beta_{\mathrm{n}}}}
\newcommand{\hypMag}{\SymColor{\beta_{\mathrm{k}}}}
\newcommand{\HJBweight}{\SymColor{\beta}}
\newcommand{\UBHJBweight}{\SymColor{\beta}_{\mathrm{U}}}

\newcommand{\LevSetVal}{\SymColor{\gamma}}
\newcommand{\LevSetRAVal}{\SymColor{\gamma_{\mathrm{A}}}}
\newcommand{\LevSetTRVal}{\SymColor{\gamma_{\mathrm{T}}}}

\newcommand{\GkCoefA}{\SymColor{\gamma_{1}}}
\newcommand{\GkCoefB}{\SymColor{\gamma_{2}}}
\newcommand{\GkCoefC}{\SymColor{\gamma_{3}}}

\newcommand{\hypCovMat}{\SymColor{\boldsymbol{\Gamma}}}
\newcommand{\sqrKcov}{\SymColor{\hypCovMat_{\mathrm{sqr}}}}

\newcommand{\KroDelta}[2]{\SymColor{\delta}_{#1,#2}}
\newcommand{\OffSetpenaltyFunc}{ \SymColor{\delta} }
\newcommand{\UBOffSetpenaltyFunc}{ \SymColor{\delta}_{\mathrm{U}} }

\newcommand{\LevSetEpsilon}{\SymColor{\epsilon}}
\newcommand{\VUBofMargin}{\SymColor{\epsilon}_{\MValF}}
\newcommand{\WUBofMargin}{\SymColor{\epsilon}_{\HdotV}}
\newcommand{\LBUBofMargin}{\SymColor{\widetilde{\epsilon}}}

\newcommand{\tempWPt}[1]{\SymColor{ \eta_{#1} }}

\newcommand{\tempUnitVec}{\SymColor{\boldsymbol{\nu}}}

\newcommand{\EigMax}[2][]{\SymColor{\lambda_{\mathrm{max}}}#1( #2 #1)}
\newcommand{\EigMin}[2][]{\SymColor{\lambda_{\mathrm{min}}}#1( #2 #1)}

\newcommand{\GPmean}{\SymColor{\boldsymbol{\mu}_{\mathrm{gp}}}}
\newcommand{\Mmean}{\SymColor{\boldsymbol{\mu}}}

\newcommand{\GPsd}{\SymColor{\boldsymbol{\sigma}_{\mathrm{gp}}}}
\newcommand{\Msd}{\SymColor{\boldsymbol{\sigma}}}
\newcommand{\LBGPsd}{\SymColor{\sigma}_{\mathrm{L}}}

\newcommand{\NewDist}{\SymColor{\tau}}
\newcommand{\UBDist}{\NewDist_{\mathrm{U}}}
\newcommand{\ODEFUBDist}{\SymColor{\overline{\NewDist}}}

\newcommand{\tempDist}{\SymColor{\chi}}
\newcommand{\SEkerExtPo}[1]{\tempDist_{#1}}

\newcommand{\NNnode}[1]{\SymColor{\boldsymbol{\xi}^{(#1)}}}
\newcommand{\PreliFunc}{\SymColor{\xi}}

\newcommand{\Basis}{\SymColor{\Phi}}
\newcommand{\PreliAFunc}{\SymColor{\phi}}
\newcommand{\PreliBFunc}{\SymColor{\psi}}
\newcommand{\subAbasis}[1]{\PreliAFunc_{#1}}
\newcommand{\subBbasis}[1]{\PreliBFunc_{#1}}
\newcommand{\tempPreliFunc}[1]{\SymColor{ \widetilde{\PreliFunc}_{#1} }}

\newcommand{\zzEachMBGset}[1]{\SymColor{\boldsymbol{\omega}(}  {#1}  \SymColor{)}}
\newcommand{\zzpEachMBGset}[1]{\SymColor{\boldsymbol{\omega}^{+}(}  {#1}  \SymColor{)}}

\newcommand{\penaltyFunc}{\SymColor{\zeta}}

\newcommand{\onlyDrift}{\SymColor{\boldsymbol{a}}}

\newcommand{\zzzBLB}[2]{\SymColor{\underline{a}}(#1,#2)}
\newcommand{\zzzBUB}[2]{\SymColor{\overline{a}}(#1,#2)}
\newcommand{\LCombiCoef}[1]{\SymColor{a_{#1}}}

\newcommand{\tempA}[1]{\SymColor{a}_{#1}}
\newcommand{\tempB}[1]{\SymColor{b}_{#1}}

\newcommand{\QFcoefB}{  \SymColor{\boldsymbol{b}_{\mathrm{q}}}   }
\newcommand{\QFcoefC}{  \SymColor{           {c}_{\mathrm{q}}}   }

\newcommand{\zzzzLB}{\SymColor{\underline{b}}}
\newcommand{\zzzzUB}{\SymColor{\overline{b}}}
\newcommand{\zzzLB}[1]{\zzzzLB(#1)}
\newcommand{\zzzUB}[1]{\zzzzUB(#1)}

\newcommand{\VCparam}{\SymColor{\boldsymbol{c}}}
\newcommand{\stabVCparam}{\SymColor{\boldsymbol{c}_{0}}}
\newcommand{\optVCparam}{\SymColor{\boldsymbol{c}_{\HJBweight}}}

\newcommand{\IDdat}{\SymColor{d}}
\newcommand{\IDbdat}{\SymColor{d^{\prime}}}

\newcommand{\Drift}{\SymColor{ \boldsymbol{f}      }}

\newcommand{\DDrift}[1]{\SymColor{\boldsymbol{f}_{\mathrm{d}}^{(#1)}}}

\newcommand{\zzzzLipConst}{\SymColor{g}}
\newcommand{\zzzLipConst}[1]{\zzzzLipConst(#1)}

\newcommand{\NLFunc}{\SymColor{h}}
\newcommand{\LipsNLFunc}{\NLFunc_{\mathrm{Lip}}}
\newcommand{\Globalddsup}{\partial^{*}(\NLFunc)}

\newcommand{\IDdom}{\SymColor{i}}
\newcommand{\IDPt}{\SymColor{j}}
\newcommand{\IDbPt}{\SymColor{j^{\prime}}}

\newcommand{\kernel}[2]{\SymColor{k}(#1,#2)}
\newcommand{\kernelVec}{\SymColor{\boldsymbol{k}_{\NumData}}}

\newcommand{\PenaltyperX}[5]{\SymColor{l}(#1,#2,#3,#4,#5)}
\newcommand{\SimplePenaltyperX}[2]{\SymColor{l}(#1,#2)}

\newcommand{\ODEFfunc}{\SymColor{{m}}}

\newcommand{\zzzzLM}{\SymColor{\underline{m}}}
\newcommand{\zzzzUM}{\SymColor{\overline{m}}}
\newcommand{\zzzLM}[1]{\zzzzLM(#1)}
\newcommand{\zzzUM}[1]{\zzzzUM(#1)}

\newcommand{\DimX}{\SymColor{n}}
\newcommand{\NumSimp}{\DimX+1}
\newcommand{\DimU}{\SymColor{n_{\mathrm{u}}}}
\newcommand{\DimVCparam}{\SymColor{n_{\mathrm{c}}}}
\newcommand{\NNdim}[1]{\DimX^{(#1)}}

\newcommand{\ODEFUBpow}{\SymColor{{p}}}

\newcommand{\NNidlayer}{\SymColor{q}}

\newcommand{\IDEl}{\SymColor{s}}
\newcommand{\IDbEl}{\SymColor{s^{\prime}}}

\newcommand{\MyT}{\SymColor{t}}

\newcommand{\Input}{\SymColor{\boldsymbol{u}}}

\newcommand{\NNweight}[2]{\SymColor{\boldsymbol{v}}_{#1,#2}}

\newcommand{\SimpVecCoef}{\SymColor{\boldsymbol{w}}}

\newcommand{\State}{\SymColor{\boldsymbol{x}}}
\newcommand{\dotState}{\SymColor{\dot{\State}}}
\newcommand{\SimpState}{\State_{\IDdom}(\SimpVecCoef)}
\newcommand{\SimpPt}[2]{\State_{#1,#2}}
\newcommand{\DState}[1]{\SymColor{\boldsymbol{x}^{(#1)}}}

\newcommand{\tempState}{\SymColor{ \widetilde{\State} }}
\newcommand{\drefState}{\SymColor{ \breve{\State} }}

\newcommand{\yBasis}{\SymColor{y}}
\newcommand{\tempyBasis}{\SymColor{\widetilde{\yBasis}}}
\newcommand{\MeanValueThmY}{\SymColor{  \tempyBasis_{\tempWPt{\ast}}  }}

\newcommand{\ArgpenaltyFunc}{\SymColor{z}}
\newcommand{\tempPPt}[2]{\SymColor{ \boldsymbol{z}_{#1,#2} }}


\newcommand{\QFcoefA}{  \SymColor{\boldsymbol{A}_{\mathrm{q}}}   }
\newcommand{\InMat}{\SymColor{\boldsymbol{B}}}

\newcommand{\ODEFUBCoef}{\SymColor{{C}}}

\newcommand{\LUMconst}[1]{\SymColor{C}_{#1}}

\newcommand{\NumData}{\SymColor{D}}

\newcommand{\HJB}{\SymColor{H}}

\newcommand{\Identity}[1]{\SymColor{\boldsymbol{I}}_{#1}}

\newcommand{\PerformanceTerm}{\SymColor{J}}

\newcommand{\kernelMat}{\SymColor{\boldsymbol{K}_{\NumData}}}

\newcommand{\ObjectiveF}{\SymColor{L}}

\newcommand{\MEdotV}{\SymColor{M}}

\newcommand{\NumNest}{\SymColor{N}}
\newcommand{\fNumNest}[1]{\SymColor{N}_{#1}}

\newcommand{\MyOrder}{\SymColor{O^{\prime}}}

\newcommand{\NNnumlayer}{\SymColor{Q}}

\newcommand{\SDdotV}{\SymColor{S}}

\newcommand{\NumBasis}{\SymColor{T}}

\newcommand{\MValFLaSalle}{\SymColor{\widetilde{V}}}
\newcommand{\MValF}{\SymColor{V}}
\newcommand{\quadMValF}{\SymColor{V}_{\mathrm{q}}}
\newcommand{\HdotV}{\SymColor{W}}

\newcommand{\VerticesestMDomStable}{\SymColor{X_{\mathrm{S}}}}
\newcommand{\NumVerticesestMDomStable}{\SymColor{\#(\VerticesestMDomStable)}}

\newcommand{\DomLIEachY}{\SymColor{{\LI{\mathbb{Y}}_{\IDdom}}}}
\newcommand{\DomLIY}{\SymColor{{\LI{\mathbb{Y}}}}}


\newcommand{\DomDrift}{\SymColor{\mathbb{F}}}

\newcommand{\SetALLIDdom}{\SymColor{\mathbb{I}}}
\newcommand{\SetStableIDdom}{\SymColor{\mathbb{I}_{\mathrm{S}}}}
\newcommand{\estSetStableIDdom}{\widehat{    \SetStableIDdom   }}

\newcommand{\NNReal}{\SymColor{\mathbb{R}_{\geq 0}}}

\newcommand{\SetVCparam}{\SymColor{\mathbb{S}}_{\VCparam}}

\newcommand{\SetDist}{ \SymColor{\mathbb{S}_{\NewDist}} }

\newcommand{\LevSetFuncLaSalle}{\SymColor{\widetilde{\mathbb{X}}_{\mathrm{Lv}}}}
\newcommand{\LevSetFunc}{\SymColor{\mathbb{X}_{\mathrm{Lv}}}}
\newcommand{\InvariantSet}{\SymColor{\mathbb{X}_{\mathrm{T}}}}
\newcommand{\RegionOfAttraction}{\SymColor{\mathbb{X}_{\mathrm{A}}}}
\newcommand{\MDomStable}{\SymColor{\mathbb{X}_{\mathrm{S}}}}
\newcommand{\estMDomStable}{\widehat{    \MDomStable   }}
\newcommand{\DomX}{\SymColor{\mathbb{X}}}
\newcommand{\tmpLevSetDom}{\DomX_{\mathrm{\LevSetEpsilon}}^{o}}
\newcommand{\tmpXBall}{\DomX_{\mathrm{b}}^{o}}
\newcommand{\DomBigX}{\SymColor{\widetilde{\DomX}}}
\newcommand{\DomSimpX}[1]{\DomX_{#1}}
\newcommand{\allDomSimpX}{\SymColor{ (\DomSimpX{\IDdom})_{\IDdom \in \SetALLIDdom}  } }

\newcommand{\DomY}{\SymColor{{\mathbb{Y}}}}
\newcommand{\DomBigY}{\SymColor{\widetilde{\mathbb{Y}}}}

\newcommand{\AlgNest}[3]{    \SymColor{ \mathcal{A}( } #1, #2 , #3 \SymColor{  ) }    }
\newcommand{\AlgpNest}[3]{    \SymColor{ \mathcal{A}^{+}( } #1, #2 , #3 \SymColor{  ) }    }

\newcommand{\FSETfunc}[2]{\SymColor{\mathcal{F}}(#2)} 
\newcommand{\pAFSETfunc}[2]{\SymColor{\mathcal{F}_{\mathrm{ndec}}^{+}(#2)}}
\newcommand{\pBFSETfunc}[2]{\SymColor{\mathcal{F}_{\mathrm{conv}}^{+}(#2)}}
\newcommand{\pFSETfunc}[2]{ \SymColor{\mathcal{F}^{+}(#2)}}
\newcommand{\NonNegativeFSETNLFunc}{\SymColor{\mathcal{F}_{\geq 0}}}

\newcommand{\FSETNLFunc}[1]{\SymColor{\mathcal{H}(#1)}}
\newcommand{\pAFSETNLFunc}[1]{\SymColor{\mathcal{H}_{\mathrm{ndec}}(#1)}}
\newcommand{\pBFSETNLFunc}[1]{\SymColor{\mathcal{H}_{\mathrm{conv}}(#1)}}


	\newcommand{\WFigConcept}{4.6in}
	\newcommand{\HFigConcept}{1.9in}
	
	\newcommand{\WFigSimplex}{3.2in}
	\newcommand{\HFigSimplex}{1.81in}
	
	\newcommand{\WFigLI}{3.8in}
	\newcommand{\HFigLI}{1.6in}
	
	\newcommand{\WFigLevelSets}{3.9in}
	\newcommand{\HFigLevelSets}{2.0in}


\newcommand{\DEFred}{red}
\newcommand{\DEFmyBlue}{myBlue}
\newcommand{\DEFmyGreen}{myGreen}
\newcommand{\DEFmyPurple}{myPurple}

\newcommand{\cfR}{\textcolor{\DEFred}}
\newcommand{\cfB}{\textcolor{\DEFmyBlue}}
\newcommand{\cfG}{\textcolor{\DEFmyGreen}}
\newcommand{\cfP}{\textcolor{\DEFmyPurple}}

\newcommand{\MyThin}{0.5pt/0.8}
\newcommand{\MyThick}{0.7pt/0.8}
\newcommand{\MyVeryThick}{0.9pt/0.8}

\newcommand{\MyFigScale}{1.0}

\definecolor{FigLSGray}{gray}{ 0.95 }
\definecolor{myFillBlue}{rgb}{ 0.88, 0.98, 1 }

\newcommand{\offsetFigRW}{+1.2}
\newcommand{\MyBendLine}[2]{
\draw[#1] 
(0,-7+#2) parabola bend(2,-0.5+#2) 
(3.2,-2.2+#2) parabola bend(5,-4.3+#2) 
(6.1,-3+#2) parabola bend(7,-2.0+#2) 
(8,-1.5+#2);
}
\newcommand{\JFigConcept}{%
\scalebox{0.85}{
	\begin{tikzpicture}[scale = \MyFigScale]
	\begin{axis}[axis y line=center, axis x line=middle, xlabel=$\State$,ylabel={}
	,xtick=\empty, ytick=\empty 
	,xmin=-1.7,xmax=9.0,ymin=-6,ymax=2	
	,width=\WFigConcept,height=\HFigConcept	,line width = \MyThin
	]

	\draw[line width = \MyThin,dashed] (1,10) -- (1,-10);
	\draw[line width = \MyThin,dashed] (2,10) -- (2,-10);
	\draw[line width = \MyThin,dashed] (3,10) -- (3,-10);
	\draw[line width = \MyThin,dashed] (4,10) -- (4,-10);
	\draw[line width = \MyThin,dashed] (5,10) -- (5,-10);
	\draw[line width = \MyThin,dashed] (6,10) -- (6,-10);
	\draw[line width = \MyThin,dashed] (7,10) -- (7,-10);
	\draw[line width = \MyThin,dashed] (8,10) -- (8,-10);

	{\MyBendLine{black, line width = \MyThick}{0}}
	{\MyBendLine{black, line width = \MyThick, dashed, \DEFmyBlue}{1.2}}
	
	\draw[\DEFmyGreen, line width = \MyVeryThick] 
	(0,-1.0)--(1,-1.0)--
	(1,0.73)--(2,0.73)--
	(3,0.73)--(3,-0.5)--
	(4,-0.5)--(4,-2.4)--
	(5,-2.4)--(5,-2.0)--
	(6,-2.0)--(6,-0.8)--
	(7,-0.8)--(7,-0.3)--
	(8,-0.3);

	\fill[\DEFred] (0,0) circle (2pt);
	\fill[\DEFred] (1,0) circle (2pt);
	\fill[\DEFred] (2,0) circle (2pt);
	\fill[\DEFred] (3,0) circle (2pt);
	\fill[\DEFred] (4,0) circle (2pt);
	\fill[\DEFred] (5,0) circle (2pt);
	\fill[\DEFred] (6,0) circle (2pt);
	\fill[\DEFred] (7,0) circle (2pt);
	\fill[\DEFred] (8,0) circle (2pt);

	\draw[line width = \MyThick, \DEFred] (5,0)--(7,-2.5);
	\node at (8.0,-3.0) %
	{\colorbox{white}{\cfR{Vertex $\SimpPt{\IDdom}{\IDPt}$}}};
		
	\draw[line width = \MyThick,  \DEFmyBlue] (5.7,-2.6)--(6.8,-3.7);
	\node at (7.8,-4.2) %
	{\colorbox{white}{\cfB{$\Basis(\State)+{\zzzUM{\Basis}}$}}};

	\draw[line width = \MyThick, \DEFmyGreen] (4.5,-2.4)--(6,-4.9);
	\node at (6.4,-5.4) %
	{\colorbox{white}{\cfG{Upper bound of $\Basis(\State)$}}};

	\draw[line width = \MyThick, black] (1.5,-1)--(2.5,-4.5);
	\node at (2.5,-4.9) %
	{\colorbox{white}{\shortstack{Stability condition\\$\Basis(\State)<0$}}};

	\draw[\DEFmyBlue,  line width = \MyThick,dashed] (-0.4,-0.5) -- (2,-0.5);		
	\draw[\DEFmyBlue,  line width = \MyThick,dashed] (-0.4,+0.73) -- (1,+0.73);		
	\draw[\DEFmyBlue, line width = \MyVeryThick, <->] (-0.29,+0.73) -- (-0.29,-0.5) node[left]%
	{\colorbox{white}{\cfB{\shortstack{\shortstack{Upper \\ margin} \\ ${\zzzUM{\Basis}}$}}}};

	\draw[\DEFred, line width = \MyVeryThick, <->] (4,+1.1) -- (5,+1.1);
	\node at (4.5 , +1.5) %
	{\cfR{$\NewDist$}};

	\node at (2 , +0.73+0.8) {\colorbox{white}{\cfG{Violation}}};
	\draw[\DEFmyGreen, line width = \MyVeryThick, <->] (1,+0.73+0.4) -- (3,+0.73+0.4);

	\end{axis}			
	\end{tikzpicture}
}
}

\definecolor{myFillRed}{rgb}{ 1, 0.9, 0.9 }
\newcommand{\MySimplexLines}[1]{\draw[line width = \MyThin, black] #1 ;}
\newcommand{\JFigSimplexes}{%
	\scalebox{0.85}{
		\begin{tikzpicture}[scale = \MyFigScale]
		\begin{axis}[axis y line=center, axis x line=middle, xlabel=${\El{\State}{1}}$,ylabel=${\El{\State}{2}}$
		,xtick=\empty, ytick=\empty 
		,xmin=-4.5,xmax=4.5,ymin=-2.1,ymax=2.8	
		,width=\WFigSimplex,height=\HFigSimplex	,line width = \MyThin
		]
		
		
		{\MySimplexLines{ (-1,2) -- (2,2)}}
		{\MySimplexLines{ (-2,1) -- (3,1)}}
		{\MySimplexLines{ (-3,-1) -- (2,-1)}}
		{\MySimplexLines{ (-2,-2) -- (1,-2)}}
		
		{\MySimplexLines{ (-3,-1) -- (-3,0)}}
		{\MySimplexLines{ (-2,-2) -- (-2,1)}}
		{\MySimplexLines{ (-1,-2) -- (-1,2)}}
		{\MySimplexLines{ (1,-2) -- (1,2)}}
		{\MySimplexLines{ (2,-1) -- (2,2)}}
		{\MySimplexLines{ (3,-0) -- (3,1)}}
		
		{\MySimplexLines{ (-3,0) -- (-1,2)}}
		{\MySimplexLines{ (-3,-1) -- (0,2)}}
		{\MySimplexLines{ (-2,-1) -- (1,2)}}
		{\MySimplexLines{ (-2,-2) -- (2,2)}}
		{\MySimplexLines{ (-1,-2) -- (2,1)}}
		{\MySimplexLines{ (0,-2) -- (3,1)}}
		{\MySimplexLines{ (1,-2) -- (3,0)}}
		
		\draw[line width = \MyThin, dashed] (-2,1) -- (-2.5,1.5);
		\draw[line width = \MyThin, dashed] (-1,2) -- (-1.5,2.5);
		\draw[line width = \MyThick, <->] (-2.25,1.25) -- (-1.25,2.25);		
		\node at (-2,2) {$\NewDist$};

		\draw[line width = \MyVeryThick, \DEFmyBlue] 
		(-3,0)--(-1,2)--(2,2)--(2,1)--(3,1)--(3,0)--(1,-2)--(-2,-2)--(-2,-1)--(-3,-1)--(-3,0);
		\draw[line width = \MyVeryThick, \DEFmyBlue,<-] (-2,-1.5)--(-3.5,-1.5) node[left]{$\cfB{\DomX}$};

		\draw[line width = \MyVeryThick, \DEFred,<-] (-2.25,+0.25)--(-3.5,1.5) node[left]{\cfR{$\DomSimpX{1}$}};
		\draw[line width = \MyVeryThick, \DEFred,<-] (-2.75,-0.25)--(-3.5,0.5) node[left]{\cfR{$\DomSimpX{2}$}};
		\draw[line width = \MyVeryThick, \DEFred,<-] (-2.25,-0.75)--(-3.5,-0.5) node[left]{\cfR{$\DomSimpX{3}$}};

		\node at (-0.2,-0.6) {\colorbox{white}{\cfR{$\SimpPt{\IDdom}{2}$}}};		
		\node at (1.2,-0.6) {\colorbox{white}{\cfR{$\SimpPt{\IDdom}{1}$}}};		
		\node at (1.2,1.5) {\colorbox{white}{\cfR{$\SimpPt{\IDdom}{3}$}}};		
		\fill[\DEFred] (-0,0) circle (3pt);
		\fill[\DEFred] (+1,0) circle (3pt);
		\fill[\DEFred] (+1,1) circle (3pt);
		\filldraw[line width = \MyVeryThick, draw=\DEFred, fill=myFillRed] (-0,0) -- (1,1) -- (1,0) -- (0,0);
		\draw[line width = \MyVeryThick, \DEFred,<-] (0.75,+0.25)--(2.25,1.75) node[right] {\cfR{$\DomSimpX{\IDdom}$}};

		
		\end{axis}			
		\end{tikzpicture}
	}
}

\newcommand{\JFigLI}{%
	\scalebox{0.85}{
		\begin{tikzpicture}[scale = \MyFigScale]
		\begin{axis}[axis y line=center, axis x line=middle, xlabel=$\State$,ylabel=$\PreliFunc(\State)$
		,xtick=\empty, ytick=\empty 
		,xmin=-1.4,xmax=7.5,ymin=-6.5,ymax=2.8	
		,width=\WFigLI,height=\HFigLI,line width = \MyThin
		]

		\draw[black, line width = \MyVeryThick] 
		(0,-3.5) parabola bend(1.5,-5.1) 
		(2,-4.0) parabola bend(2.5,-2.5) 
		(3,-3.0) parabola bend(3.5,-3.5)  
		(4,-2.0) parabola bend(4.7,-0.7)  
		(6,-1.5);
		\draw[\DEFmyBlue, line width = \MyVeryThick] (0,-3.5)  --               (2,-4)    --               (4,-2)  --                    (6,-1.5);
		\draw[line width = \MyThin, \DEFmyBlue, <-] (5.0,-1.75)--(6.2,-3.5) node[right]{\cB{$\LI{\PreliFunc}(\State)$}};	
		\draw[line width = \MyThin, black, <-] (3.85,-2.7)--(6.2,-5.5) node[right]{$\PreliFunc(\State)$};

		\draw[line width = \MyThin,dashed] (2,0) -- (2,-8);
		\draw[line width = \MyThin,dashed] (4,0) -- (4,-8);
		\draw[line width = \MyThin,dashed] (6,0) -- (6,-8);
		
		\node at (1 , -6) {${\DomSimpX{1}}$};
		\node at (3 , -6) {${\DomSimpX{2}}$};
		\node at (5 , -6) {${\DomSimpX{3}}$};	

		\fill[fill=\DEFred] (0,-3.5) circle (2.0pt);
		\fill[fill=\DEFred] (2,-4.0) circle (2.0pt);
		\fill[fill=\DEFred] (4,-2.0) circle (2.0pt);
		\fill[fill=\DEFred] (6,-1.5) circle (2.0pt);

		\draw[line width = \MyThin, \DEFred] (0,-3.5) -- (0,0) node[above left]{\cfR{$\SimpPt{1}{1}$}};	
		\draw[line width = \MyThin,\DEFred] (2,-4.0) -- (2,0) node[above]{\cfR{\shortstack{$\SimpPt{1}{2}$\\$\SimpPt{2}{1}$}}};
		\draw[line width = \MyThin,\DEFred] (4,-2.0) -- (4,0) node[above]{\cfR{\shortstack{$\SimpPt{2}{2}$\\$\SimpPt{3}{1}$}}};
		\draw[line width = \MyThin,\DEFred] (6,-1.5) -- (6,0) node[above]{\cfR{\shortstack{$\SimpPt{2}{2}$\\$\SimpPt{3}{1}$}}};
		
		\newcommand{\MyFigMargin}{1.5}	
		
		\draw[\DEFmyBlue, line width = \MyThick, <->] (-0.2,-3.5) -- (-0.2,-3.5-\MyFigMargin) node[left]{\cfB{$\zzzLM{\PreliFunc}$}};	
		\draw[\DEFmyBlue, line width = \MyThick, <->] (-0.2,-3.5) -- (-0.2,-3.5+\MyFigMargin) node[left]{\cfB{$\zzzUM{\PreliFunc}$}};

		\draw[\DEFmyBlue, line width = \MyThin,dashed] (0,-3.5-\MyFigMargin) -- (2,-4-\MyFigMargin);	
		\draw[\DEFmyBlue, line width = \MyThin,dashed] (0,-3.5+\MyFigMargin) -- (2,-4+\MyFigMargin);	
		
		\draw[\DEFmyBlue, line width = \MyThin,dashed] (2,-4-\MyFigMargin) -- (4,-2-\MyFigMargin);	
		\draw[\DEFmyBlue, line width = \MyThin,dashed] (2,-4+\MyFigMargin) -- (4,-2+\MyFigMargin);

		\draw[\DEFmyBlue, line width = \MyThin,dashed] (4,-2-\MyFigMargin) -- (6,-1.5-\MyFigMargin);
		\draw[\DEFmyBlue, line width = \MyThin,dashed] (4,-2+\MyFigMargin) -- (6,-1.5+\MyFigMargin);	
		
		
		\end{axis}			
		\end{tikzpicture}
	}
}

\newcommand{\FigLSscaleR}{0.61}
\newcommand{\MyLevelSetDraw}[2]{
	\draw[#1] (-0.7,+0.7) circle [x radius =#2* 0.6 * \FigLSscaleR in, y radius =#2* 0.3 * \FigLSscaleR in, rotate = +10];	
}

\newcommand{\JFigLevelSets}{%
	\scalebox{0.85}{
		\begin{tikzpicture}[scale = \MyFigScale]
		\begin{axis}[axis y line=none, axis x line=none, xlabel=\empty,ylabel=\empty
		,xtick=\empty, ytick=\empty 
		,xmin=-8.5,xmax=12.5,ymin=-8.5,ymax=10.3	
		,width=\WFigLevelSets,height=\HFigLevelSets	,line width = \MyThin
		]

		\filldraw[line width = \MyThick, draw=black, fill=FigLSGray] (-8,-8) rectangle (8,8);
		\node at (+10, 7) {${\DomX}$};	
		\draw[line width = \MyThick,black, <-] (7,6) -- (9.5,7);

		\filldraw[line width = \MyVeryThick, draw=\DEFmyBlue, fill=myFillBlue] 
		(-8,+4)--(-6,+6)--(-2,+6)--(-0,+8)--(+4,+8)--(+4,+6)--(+6,+6)--(+6,+4)--(+8,+4)--
		(+8,-2)--(+2,-8)--
		(-2,-8)--(-2,-6)--(-4,-6)--(-6,-6)--(-6,-4)--(-8,-4)--(-8,+4);
		
		\filldraw[line width = \MyVeryThick, draw=\DEFmyBlue, fill=FigLSGray] (-2,-0)--(+0,+2)--(+0,0)--(-2,-0);
		
		\node at (+10.2, 4) {$\cfB{\MDomStable}$};
		\draw[line width = \MyVeryThick,\DEFmyBlue, <-] (7,2) -- (9.5,4) ; 
		
		{\MyLevelSetDraw{\DEFmyGreen, line width = \MyThick, dashed}{1.5}}
		{\MyLevelSetDraw{\DEFmyGreen, line width = \MyThick, dashed}{2.7}}
						
		\node at (9.2, -3.5) {\colorbox{white}{
				{\cfG{Level set $\LevSetFunc(\LevSetVal)$ }}
		}};	
		\node at (9, -6.0) {\colorbox{white}{
				{\cfG{of $\MValF(\State)$}}
		}};				
		\draw[\DEFmyGreen, line width = \MyVeryThick, <-] (4.0,-2.8) -- (5.7,-3.5); 

		{\MyLevelSetDraw{\DEFred, line width = \MyVeryThick}{0.6}}
		{\MyLevelSetDraw{\DEFred, line width = \MyVeryThick}{2.3}}
		
		\node at (1.5 , 1) {$\cfR{\InvariantSet}$};
		\node at (5.6, 2) {\colorbox{myFillBlue}{$\cfR{\RegionOfAttraction}$}};	
		
		\end{axis}			
		\end{tikzpicture}
	}
}



\allowdisplaybreaks[4]

\begin{abstract}
	This study presents a sampling-based method to guarantee robust stability of general control systems with uncertainty.
	The method allows the system dynamics and controllers to be represented by various data-driven models, such as Gaussian processes and deep neural networks. 
	For nonlinear systems, stability conditions involve inequalities over an infinite number of states in a state space.
	Sampling-based approaches can simplify these hard conditions into inequalities discretized over a finite number of states.
	However, this simplification requires margins to compensate for discretization residuals.
	Large margins degrade the accuracy of stability evaluation, and obtaining appropriate margins for various systems is challenging.
	This study addresses this challenge by deriving second-order margins for various nonlinear systems containing data-driven models.
	Because the size of the derived margins decrease quadratically as the discretization interval decreases, the stability evaluation is more accurate than with first-order margins.
	Furthermore, this study designs feedback controllers by integrating the sampling-based approach with an optimization problem.
	As a result, the controllers can guarantee stability while simultaneously considering control performance.	
\end{abstract}
\section{Introduction} \label{sec_intro}

Data-driven system identification and controller design are timely topics in control engineering \cite{Nghiem2023ACCsurveyPIML}.
For example, Gaussian processes (GPs) \cite{Rasmussen06,YIto2024arXiv}, deep neural networks (DNNs) \cite{Aggarwal18}, kernel-based models \cite{Vovk2013,ItoAccess20}, and reservoir computing \cite{Yan2024Ncom,Zhai2023Ncom} are promising methods for representing complex dynamics with little prior knowledge.
Such models have been applied to various areas, including controller design \cite{Hewing2020review}, robotics \cite{Brunke2022review}, and manipulators \cite{Zhai2023Ncom}.
These successful applications motivate us to focus on the control of nonlinear systems using data-driven approaches.
A groundbreaking data-driven control approach has been proposed for linear systems \cite{Persis2020TAC}.
Meanwhile, control theory has also contributed to the machine learning community; for instance, Lyapunov stability has been employed for fast neural network (NN) inference \cite{Rodriguez2022ICML}.

Although a major challenge is to formally guarantee stability of general control systems containing data-driven models, existing methods have unfortunately focused on specific models and come with some limitations, which can be categorized into three topics. 
Firstly, verification methods have been applied to nonlinear and/or data-driven models 
\cite{Deka2022CDC,Alfarano2024NeurIPS,Jagtap2021TAC,Wajid2022L4DC,Jagtap2020CDC,Abate2021CSL,Chang19,Rose2023CDC}.
Solvers based on the sum of squares are applicable only for polynomial systems \cite{Deka2022CDC,Alfarano2024NeurIPS,Jagtap2021TAC}.	
Whereas satisfiability modulo theories (SMT) solvers are powerful tools for verifying general conditions \cite{Jagtap2021TAC,Wajid2022L4DC,Jagtap2020CDC,Chang19}, they often fail or time out, meaning their termination is not guaranteed \cite{Alfarano2024NeurIPS,Abate2021CSL}. 
Nonlinear optimization-based verification lacks formal guarantees owing to local optimality issues \cite{Deka2022CDC}. 
Numerical evaluations at only finite points result in the loss of exact stability guarantee \cite{ItoAccess20}.
Scenario-based approaches have provided stability guarantees only in a probabilistic sense \cite{Rose2023CDC}.
Secondly, while data-driven models have been analytically investigated for the evaluation of stability conditions, existing results remain conservative and are limited to specific models \cite{Scharnhorst2023TAC,Fazlyab2022TAC,Yin2022TAC,Jackson2020CDC}.
If target dynamics is represented by kernel functions, bounds of the function have been analyzed \cite{Scharnhorst2023TAC}.
Sector-based approximations of NNs' outputs have been obtained \cite{Fazlyab2022TAC}. 
Quadratic approximations restrict regions of attraction to ellipsoids for NN-based control \cite{Yin2022TAC}.
When a system is modeled using GPs combined with Markov decision processes, a probability bound on system safeness has been derived \cite{Jackson2020CDC}.
Thirdly,
learning dynamics while preserving proper properties does not clarify their boundaries, such as a region of attraction, and/or assumes that target dynamics already have the properties
\cite{Kolter2019NeurIPS,Kojima2022NeurIPS,Miyano2024CSL,Okamoto2024arXiv,Mohammadi2024ICLR,Mohammadi2024arXiv,Roth2025arXiv,Revay2024TAC,DAmico2024TAC}.
The preserved properties include asymptotic stability \cite{Kolter2019NeurIPS}, input-output stability \cite{Kojima2022NeurIPS}, passivity \cite{Miyano2024CSL}, dissipativity \cite{Okamoto2024arXiv}, contraction \cite{Mohammadi2024ICLR,Mohammadi2024arXiv}, and port-Hamiltonian properties \cite{Roth2025arXiv}.
Recurrent equilibrium networks (RENs) are generalized models that incorporate contracting properties and include DNNs \cite{Revay2024TAC,Wang2023CSL}.
Incremental input-state stability of recurrent NNs has been analyzed \cite{DAmico2024TAC}. 
Additionally, some methods rely on idealized assumptions, such as no regression loss \cite{Beppu21} and the ability to freely control all state variables \cite{Umlauft18CSL}.
Stochastic boundedness has been discussed only for strict-feedback systems \cite{Wang21}.

\tikzstyle{block} = [draw, rectangle, minimum height=2em, minimum width=3em]
\begin{figure}[!t]  
	\begin{center}
		\JFigConcept
		\vspace*{-0.1in}
	\end{center}
	\caption{Concept of sampling-based approaches for guaranteeing the stability condition $\Basis(\State)<0$ for all states $\State$.}
	\label{fig:RelatedWork} 
\end{figure}

Apart from the aforementioned methods, sampling-based approaches have the potential to guarantee stability of various dynamics with minimal conservativeness.
The concept is illustrated in Fig.\ref{fig:RelatedWork}.
Stability conditions are expressed as inequalities  $\Basis(\State)<0$ for an \textit{infinite} number of states $\State$.
These hard conditions can be simplified by discretizing the inequalities over a \textit{finite} number of vertices $\SimpPt{\IDdom}{\IDPt}$ and introducing upper bounds.
The bounds are constructed by adding \textit{margins} ${\zzzUM{\Basis}}$ to compensate for discretization residuals.
If small margins can be obtained, the stability conditions can be evaluated more precisely.

However, a critical challenge in sampling-based approaches is obtaining \textit{small} margins for \textit{various} functions.
Especially, \textit{second-order} margins are needed, as their sizes decrease quadratically with a decreasing sampling interval $\NewDist$ in the discretization.
These margins are superior to first-order margins, which are only proportional to $\NewDist$.
Moreover, even obtaining margins for various functions is challenging.
First-order margins based on Lipschitz constants have been used for specific models, 
such as GPs \cite{Berkenkamp16CDC,Berkenkamp17,Lederer19,Lederer2019CDC}
and 
NNs \cite{Richards18}, as well as for designing control barrier (CB) functions \cite{Awan2023CSL} and multiple CB certificates \cite{Nejati2023CSL}.
Taylor expansion-based margins are at best first-order \cite{Bobiti2018TAC}.
Continuous piecewise affine methods \cite{Marinosson02,Giesl15} restrict the class of applicable functions because they require bounds on second-order derivatives of the functions.
A delta-cover method requires the modulus of continuity, which is generally difficult to obtain for various functions  \cite{Marchi2022CSL,Marchi2021L4DC}.
In addition, other sampling-based methods \cite{Chen2021CDC,ItoIFAC17,Vinogradska17} have been applied only to specific models.

To overcome the aforementioned limitations, we propose a sampling-based method for guaranteeing stability for various control systems.
The proposed method can handle general system classes that include a wide range of data-driven models.
We derive second-order margins for functions in the general classes, which are less conservative than first-order margins such as those based on Lipschitz constants.
The main contributions of this study are summarized as follows:

\begin{enumerate}
	\item 
	Section \ref{sec_2ndOrderMargin}:
	We derive second-order margins of functions contained in general classes denoted by ${\FSETfunc{\NumBasis}{\NumNest}}$ and ${\pFSETfunc{\NumBasis}{\NumNest}}$ (Theorems \ref{thm:main_results_margins} and \ref{thm:main_results_margins_specific}).
	Examples of functions in these classes include GPs, DNNs, kernel-based models, polynomials, transcendental functions, and their multiple compositions (Table \ref{tab:results} and Theorems \ref{thm:ex_simple_functions}--\ref{thm:ex_GPsd}).
	Using the margins, we derive precise upper and lower bounds of the functions, which are less conservative than those obtained using first-order margins.
	These contribute to realize the generality of system classes in sampling-based analyses, having the potential for integration with other methods, such as \cite{Marchi2021L4DC,Awan2023CSL}.

	\item
	Section \ref{sec_solution1}:
	Using the proposed second-order margins, we analyze robust stability for nonlinear systems with uncertainty and Lyapunov functions belonging to the proposed classes ${\FSETfunc{\NumBasis}{\NumNest}}$ and ${\pFSETfunc{\NumBasis}{\NumNest}}$ (Theorems \ref{thm:find_stability_region} and \ref{thm:HdotV_margin}).
	The stability analysis reduces to finding stability regions where stability conditions hold.
	To achieve this, we propose a sampling-based method that incorporates the second-order margins.

	\item
	Section \ref{sec_solution2}:	
	We design stabilizing controllers using the sampling-based method (Theorem \ref{thm:general_control_design}).
	The stability guarantee is integrated with various performance indices, such as cost functions and optimality residuals.
	The control design is formulated as optimization problems regarding parameters of controllers and Lyapunov functions.
	Solving the problems yields controllers that stabilize systems under certain technical assumptions.

\end{enumerate}

This paper is an extended version of the authors' conference paper \cite{ItoCDC19}.
This paper focuses on stability and controller design for general data-driven systems whereas the conference paper has focused solely on GPs without controller design.

\textit{\bfseries Remainder of this paper:} 
The main problems are introduced in Section \ref{sec_problem}.
Section \ref{sec_2ndOrderMargin} describes a key method to solve the problems.	 
Solutions to the problems are presented in Section \ref{sec_solutions}.
Section \ref{sec_demonstration} demonstrates the effectiveness of the proposed method.
Finally, Section \ref{sec_conclusion} concludes this study.

\textit{\bfseries Tips for readers:}
This paper establishes a comprehensive theory.
For a brief understanding, it is easy to follow Sections \ref{sec_problem}, the beginning of \ref{sec_2ndOrderMargin}, \ref{sec_2ndOrderMargin_summary}, \ref{sec_solutions}, and \ref{sec_demonstration} before reading Sections \ref{sec_set_of_func}--\ref{sec_set_of_func_LB}.

%
\textit{\bfseries Notation:}
The following notation is used in this paper.
\begin{itemize}
	
	\item
	$\NNReal$: the set of nonnegative real numbers
	
	\item
	$\mathbb{N}:=\{1,2,\dots\}$: the set of natural numbers

	\item 
	$\Identity{\DimANotation}\in \mathbb{R}^{\DimANotation \times \DimANotation}$: the identity matrix
	
	\item
	$\El{\NotationVec}{\IDNotation}$: the $\IDNotation$th component of a vector $\NotationVec$
	\item
	$\El{\NotationMat}{\IDNotation,\IDbNotation}$: the component in the $\IDNotation$th row and $\IDbNotation$th column of a matrix $\NotationMat$

	\item
	${\EigMin{\NotationMat}}$ (resp. ${\EigMax{\NotationMat}}$): the minimum (resp. maximum) eigenvalue of a symmetric matrix $\NotationMat$

	\item
	${\partial \NotationVec(\NotationbVec) }/{\partial \NotationbVec^{\MyTRANSPO}} $: 
	the partial derivative of $\NotationVec:  \mathbb{R}^{\DimANotation} \to  \mathbb{R}^{\DimBNotation}$, 
	where
	$
	{\El{  {\partial \NotationVec(\NotationbVec) }/{\partial \NotationbVec^{\MyTRANSPO}}  }{\IDNotation,\IDbNotation}}
	=  {\partial {\El{\NotationVec(\NotationbVec)}{\IDNotation}}  }/
	{\partial {\El{\NotationbVec}{\IDbNotation}}   } 	
	$

\end{itemize}

\section{Problem setting} \label{sec_problem}

\subsection{Target systems with data-driven control problems} \label{sec_target_sys}

Consider the following nonlinear system:
\begin{align}
\dotState(\MyT)
 = \Drift(\State(\MyT),  \Input(\State(\MyT)) ) 
,\label{eq:def_sys}
\end{align}
where $\State(\MyT) \in  \mathbb{R}^{\DimX}$ is the state at the time $\MyT$.
Let $\Input: \mathbb{R}^{\DimX} \to  \mathbb{R}^{\DimU} $ be a state feedback controller to be designed in this study.
The nonlinear dynamics $\Drift: \mathbb{R}^{\DimX} \times \mathbb{R}^{\DimU} \to \mathbb{R}^{\DimX}$ 
is assumed to be $C^{1}$ continuous and satisfy $\Drift(0,0)=0$.

This study considers both cases that $\Drift$ is given and that unknown $\Drift$ is identified by a data-driven model using a data set.
Let $\DomDrift: \mathbb{R}^{\DimX}\times \mathbb{R}^{\DimU} \to 2^{\mathbb{R}^{\DimX}}$ be the model set  of $\Drift$ that is defined using nominal dynamics $\Mmean:\mathbb{R}^{\DimX}\times \mathbb{R}^{\DimU} \to \mathbb{R}^{\DimX}$ and uncertainty $\Msd:\mathbb{R}^{\DimX}\times \mathbb{R}^{\DimU} \to \NNReal^{\DimX}$ as follows:
\begin{align}	
&
\DomDrift(\State, \Input)
\nonumber\\&
:=
\Bigg\{  \Mmean(\State, \Input) 
+
\begin{bmatrix}
\ModelSetFreeParam{1} \quad\quad\;\;\;\; \\ \quad\;\; \ddots\quad\;\;\\ \quad\quad\;\;\;\; \ModelSetFreeParam{\DimX}
\end{bmatrix}
\Msd(\State, \Input) 
\Bigg| \; \forall \IDEl ,  - 1 \leq \ModelSetFreeParam{\IDEl} \leq 1  
\Bigg\}
.
\end{align}	
Let $\DomX \subset \mathbb{R}^{\DimX}$ be a given bounded closed set containing the origin $\{0\}$.
The following is assumed throughout this paper.

\begin{assumption}[{\MyHighlight{Uncertainty description}}]\label{ass:f_err}
Given a feedback controller $\Input$, the model set $\DomDrift$ satisfies
	\begin{align}	
	&
	\forall \State \in \DomX
	,\;
	\Drift(\State, \Input(\State))   \in \DomDrift(\State, \Input(\State))
	.
	\end{align}
\end{assumption}
\begin{remark}[{\MyHighlight{Model set}}]
If a deterministic $\Drift$ is given, $\DomDrift(\State, \Input(\State)) =\{ \Drift(\State, \Input(\State)) \}$ is obtained by $\Mmean=\Drift$ and $\Msd=0$.
If $\Drift$ is unknown, 
$\Mmean$ and  $\Msd$ correspond to the nominal dynamics and uncertainty of a data-driven model, respectively.  
\end{remark}

\subsection{Main problems} \label{sec_main_problem}

The main objective of this study is to guarantee the following stability for the system \eqref{eq:def_sys}.

\begin{definition}[\MyHighlight{Practical stability}] \label{def:stability}
	The system \eqref{eq:def_sys} is said to be practically stable 
	if there exist a \textit{region of attraction} $\RegionOfAttraction \subseteq \DomX$ and \textit{target region} $\InvariantSet \subseteq \RegionOfAttraction$ that satisfy  
	\begin{align}	
	\forall \State(0) \in \RegionOfAttraction , 
	\quad
	\lim_{\MyT \to \infty}
	\Big( \inf_{  \tempState \in  \InvariantSet      } \|\State(\MyT)  -  \tempState \| \Big)
	&= 0
	.\label{eq:def_asymptotic_stability}
	\end{align}
\end{definition}

This notion indicates that $\State(\MyT)$ reaches $\InvariantSet$ asymptotically for any initial $\State(0)$ on $\RegionOfAttraction$.
Guaranteeing the practical stability reduces to finding a stability region $\MDomStable$ defined below.
\begin{definition}[\MyHighlight{Stability region}] \label{def:stability_region}
For a $C^{1}$ continuous function $\MValF: \mathbb{R}^{\DimX} \to \mathbb{R} $, controller $\Input$, nominal dynamics $\Mmean$, and uncertainty $\Msd$,
let $\MDomStable \subseteq \DomX$ be a stability region satisfying 
\begin{align}
&
\forall \State \in \MDomStable
,\;  \MValF(\State) > 0 
,\; \HdotV(\State)  < 0  
,\label{eq:stab_conditions}
\end{align}
where
\begin{align}
\HdotV(\State)
&:=
\MEdotV(\State)
+
\SDdotV(\State)
,  \label{eq:def_HdotV}
\\
\MEdotV(\State)
&
:=
\frac{\partial \MValF (\State)}{\partial \State }^{\MyTRANSPO} \Mmean(\State, \Input(\State)) 	
, \label{eq:def_MEdotV} 
\\
\SDdotV(\State)
&
:=\sum_{\IDEl=1}^{\DimX} 
\Big| {\El[\Big]{\frac{\partial \MValF (\State)}{\partial \State }}{\IDEl}} 
 \Big|  
 {\El{\Msd(\State, \Input(\State))}{\IDEl}} 
. \label{eq:def_SDdotV} 
\end{align}	
Note that $\MDomStable$ is associated with $\MValF$, $\Input$, $\Mmean$, and $\Msd$.
\end{definition}

The condition \eqref{eq:stab_conditions} indicates a robust version of Lyapunov stability theory because 
$
\dot{\MValF}(\State) 
\leq \HdotV(\State) < 0
$ for every $\Drift \in \DomDrift$.
As illustrated in Fig. \ref{fig:LevelSets}, if a stability region $\MDomStable$ is found,
we can obtain a region of attraction $\RegionOfAttraction$ and target region $\InvariantSet$ as follows,
where similar results are seen in \cite[Theorem 2.5]{Bjornsson15}.

\begin{proposition}[\MyHighlight{Region of attraction and target region}] \label{thm:LaSalle}
	For any $C^{1}$ continuous $\MValF: \mathbb{R}^{\DimX} \to \mathbb{R} $ and any $\LevSetVal \in \mathbb{R}$,
	let  
	$\LevSetFunc(\LevSetVal) := \{ \State \in \DomX | \MValF(\State) \leq  \LevSetVal \} $.
		For any stability region $\MDomStable$ and parameters $(\LevSetRAVal, \LevSetTRVal) \in \mathbb{R}^{2}$, 
		the following $\RegionOfAttraction$ and $\InvariantSet$:
		\begin{align}
			\RegionOfAttraction	
			&= \LevSetFunc(\LevSetRAVal)
			, \label{eq:LaSalle_RegionOfAttraction}
			\\
			\InvariantSet	
			&=  \LevSetFunc(\LevSetTRVal)
			\supseteq \{0\}
			,
			\label{eq:LaSalle_InvariantSet}
		\end{align}
		are a region of attraction and a target region for the practical stability, respectively, if the following conditions are satisfied:
		\begin{align}
		\MValF(0) \leq \LevSetTRVal &< \LevSetRAVal  < \inf_{\State \in {\setB{\DomX}} } \MValF(\State)
		\label{eq:LaSalle_cond1}
		,\\
		\LevSetFunc(\LevSetRAVal) &\setminus  \LevSetFunc(\LevSetTRVal) \subseteq  \MDomStable
		\label{eq:LaSalle_cond2}
		,
		\end{align}
		where ${\setB{\DomX}}$ denotes the boundary of $\DomX$.	
	
\end{proposition}
\begin{proof}
The proof is given in Appendix \ref{pf:LaSalle}.
\end{proof}	

\tikzstyle{block} = [draw, rectangle, minimum height=2em, minimum width=3em]
\begin{figure}[!t]  
	\begin{center}		
		\JFigLevelSets
		\vspace*{-0.1in}
	\end{center}
	\caption{
		Illustration of  a stability region  $\MDomStable$, region of attraction $\RegionOfAttraction$, target region $\InvariantSet$, and level sets of $\MValF(\State)$.}
	\label{fig:LevelSets} 
\end{figure}

By virtue of Proposition \ref{thm:LaSalle}, guaranteeing the practical stability reduces to the following main problems for finding a stability region associated with designing $\MValF$ and $\Input$.

\textit{\bfseries Problem 1 (Stability analysis):}
Find 
a stability region $\MDomStable$ 
for a given Lyapunov function $\MValF$, given feedback controller $\Input$, given nominal dynamics $\Mmean$, and given uncertainty $\Msd$.

\textit{\bfseries Problem 2 (Controller design):}
Design a Lyapunov function $\MValF$ and feedback controller $\Input$ such that a stability region $\MDomStable$ contains a given candidate region $\estMDomStable \subseteq \DomX$ for a given nominal dynamics $\Mmean$ and given uncertainty $\Msd$.

\section{Proposed method: Sampling-based analysis using second-order margins}\label{sec_2ndOrderMargin}

\begin{table}
	\centering
	\renewcommand{\arraystretch}{1.1}	
	\caption{
		Examples of functions in ${\FSETfunc{\NumBasis}{\NumNest}}$ and ${\pFSETfunc{\NumBasis}{\NumNest}}$.}
	\label{tab:results}
	\begin{tabular}{| c|c |c |}
		\hline
		{Functions} & {Class} & {Details}  \\		\hline\hline

		\shortstack{ 
			Sum and Product $\sum_{\IDEl} {\subAbasis{\IDEl}} {\subBbasis{\IDEl}}$ 
			\\	
			$({\subAbasis{\IDEl}}, {\subBbasis{\IDEl}}\in {\FSETfunc{\NumBasis}{\NumNest-1}})$
		}
		& ${\FSETfunc{\NumBasis}{\NumNest}}$  & Def. \ref{def:nested_set} \\ \hline

		\shortstack{ 	Since, cosine, sigmoid, and \\ tanh  of ${\subAbasis{}} \in {\FSETfunc{\NumBasis}{\NumNest-1}}$}
		& ${\FSETfunc{\NumBasis}{\NumNest}}$ &
		\shortstack{ Def. \ref{def:nested_set},  \\  Prop. \ref{thm:ex_classH} }
		\\	\hline

		Quadratic functions & ${\FSETfunc{\NumBasis}{1}}$ & Thm. \ref{thm:ex_simple_functions} \\		\hline
		$\NumNest$th order polynomials & ${\FSETfunc{\NumBasis}{\NumNest-1}}$ & Thm. \ref{thm:ex_simple_functions} \\	\hline
		Squared exponential kernels & ${\FSETfunc{\NumBasis}{1}}$ & Thm. \ref{thm:ex_SEkernel}   \\		\hline
		Mean of GPs & ${\FSETfunc{\NumBasis}{2}}$ & Thm. \ref{thm:ex_GP}   \\		\hline
		Standard deviation of GPs & ${\pFSETfunc{\NumBasis}{4}}$ & Thm. \ref{thm:ex_GPsd}  \\		\hline
		$\NNnumlayer$-layered DNNs & ${\FSETfunc{\NumBasis}{\NNnumlayer}}$  & Thm. \ref{thm:ex_NN}  \\		\hline
		
		The function $\MEdotV$ in \eqref{eq:def_MEdotV}${}^{*}$  
		& ${\FSETfunc{\NumBasis}{ \max\{ {\fNumNest{\partial\MValF}} ,  \fNumNest{\Mmean} \} + 1}}$ 
		& Thm. \ref{thm:HdotV_margin}
		\\		\hline
		The function $\SDdotV$ in \eqref{eq:def_SDdotV}${}^{*}$ 
		& ${\pFSETfunc{\NumBasis}{ \max\{ {\fNumNest{\partial\MValF}} ,  \fNumNest{\Msd} \} + 2}}$  
		& Thm. \ref{thm:HdotV_margin}
		\\		\hline
		
	\end{tabular}
	\vspace*{+0.01in}	
	
	${}^{*}$Suppose ${\El{\partial\MValF/\partial\State}{\IDEl}}  \in {\FSETfunc{\NumBasis}{\fNumNest{\partial\MValF}}}$, ${\El{\Mmean(\wildcard,\Input(\wildcard))}{\IDEl}} \in {\FSETfunc{\NumBasis}{\fNumNest{\Mmean}}}$, and ${\El{\Msd(\wildcard,\Input(\wildcard))}{\IDEl}} \in {\FSETfunc{\NumBasis}{\fNumNest{\Msd}}}
	\cup {\pFSETfunc{\NumBasis}{\fNumNest{\Msd}+1}}$.
\end{table}

In this section, we propose an efficient method for solving Problems 1 and 2.
This method provides solutions to the following general problem, where let $\DomBigX \subset \mathbb{R}^{\DimX}$ be a bounded open set containing the set $\DomX$.

\textit{\bfseries Problem 3:}
For a continuous function $\Basis:\DomBigX \to \mathbb{R}$, 
find lower and/or upper bounds of $\Basis$ on the set $\DomX$, provided that $\DomX$ is equal to a union of simplexes $\allDomSimpX$ to be defined below.

\begin{remark}[{\MyHighlight{Application to stability analysis}}]
	Solutions to Problem 3 are promising for solving Problems 1 and 2.
	If $ \MValF$ and $\HdotV$ are adopted as $\Basis$, we can discriminate whether the conditions $ \MValF(\State)>0$ and $\HdotV(\State)<0$ in \eqref{eq:stab_conditions} are satisfied on each simplex $\DomSimpX{\IDdom}$, respectively, using lower and upper bounds of $ \MValF$ and $\HdotV$.
\end{remark}

Examples of functions $\Basis$ in Problem 3 are listed in Table. \ref{tab:results}, where the classes ${\FSETfunc{\NumBasis}{\NumNest}}$ and ${\pFSETfunc{\NumBasis}{\NumNest}}$ for $\NumNest\in \mathbb{N}$ are to be defined later.
In the following, Section \ref{sec_2ndOrderMargin_summary} presents main results to solve Problem 3.
Sections \ref{sec_set_of_func}--\ref{sec_set_of_func_LB} provide the details of the solutions.
We use several key definitions 
as follows and as illustrated in Figs. \ref{fig:simplexes} and \ref{fig:LI_margins}.

\begin{definition}[\MyHighlight{Simplexes (Fig. \ref{fig:simplexes})}]	\label{def:simplexes}
	Let $\SetALLIDdom \subset \mathbb{N}$ be a finite set.  
	For each $\IDdom \in \SetALLIDdom$, the simplex $\DomSimpX{\IDdom}$ and its member $\SimpState$ are defined as follows:
	\begin{align}
	\DomSimpX{\IDdom}&:=
	\Big\{
	\SimpState 
	\Big|
	\forall \IDPt 
	,\;  
	\El{\SimpVecCoef}{\IDPt}\geq 0 
	,\; 
	\sum_{\IDPt=1}^{\NumSimp}\El{\SimpVecCoef}{\IDPt}=1	
	\Big\}	
	,\label{eq:def_tilde_x_set}
	\\
	\SimpState
	&:=
	\sum_{\IDPt=1}^{\NumSimp}\El{\SimpVecCoef}{\IDPt}\SimpPt{\IDdom}{\IDPt}
	,\label{eq:def_tilde_x}	
	\end{align}	
	where the vertices $\SimpPt{\IDdom}{\IDPt} \in \DomX$ for $\IDPt= \{1,2,\dots,\NumSimp\}$ are chosen such that $\SimpPt{\IDdom}{\IDPt}-\SimpPt{\IDdom}{\NumSimp}$ for $\IDPt= \{1,\dots,\DimX\}$ are linearly independent and that the following relation holds:
	\begin{align}
	\DomX=\bigcup_{\IDdom \in \SetALLIDdom} \DomSimpX{\IDdom}
	.
	\end{align}
Let $\allDomSimpX$ denote the collection of simplexes that contains the information of all the vertices $\SimpPt{\IDdom}{\IDPt}$.
\end{definition}	
\begin{definition}[\MyHighlight{Maximum sampling interval (Fig. \ref{fig:simplexes})}]	\label{def:interval}	
	Given $\allDomSimpX$, the maximum sampling interval between vertices $\SimpPt{\IDdom}{\IDPt} \in \DomX$ is defined as follows:
	\begin{align}
	\NewDist:=
	\max_{\IDdom \in \SetALLIDdom,\IDPt,\IDbPt}
	{\|\SimpPt{\IDdom}{\IDPt}-\SimpPt{\IDdom}{\IDbPt}\|}
	. \label{eq:def_tau}
	\end{align}		
\end{definition}

\begin{definition}[\MyHighlight{Linear interpolations (Fig. \ref{fig:LI_margins})}] \label{def:LI}
	Given $\allDomSimpX$ and a continuous function $\PreliFunc: \DomBigX \to \mathbb{R}$, 
	let $\LI{\PreliFunc}: \DomX \to \mathbb{R}$ denote the linear interpolation of $\PreliFunc$ satisfying
	\begin{align}
	\forall \IDdom \in \SetALLIDdom
	,\;
	\forall \SimpState \in \DomSimpX{\IDdom}
	,\quad
	\LI{\PreliFunc}(\SimpState)
	:=
	\sum_{\IDPt=1}^{\NumSimp}\El{\SimpVecCoef}{\IDPt} \PreliFunc(\SimpPt{\IDdom}{\IDPt}) 
	.\label{eq:def_linear_interpolation}
	\end{align}	
\end{definition}

\begin{definition}[\MyHighlight{Lower and upper margins (Fig. \ref{fig:LI_margins})}] \label{def:margins} 
	Given $\allDomSimpX$, a continuous function $\PreliFunc:\DomBigX \to \mathbb{R}$, and its linear interpolation $\LI{\PreliFunc}$, let ${\zzzLM{\PreliFunc}} \geq 0$ and ${\zzzUM{\PreliFunc}} \geq 0$ be lower and upper margins of $\PreliFunc(\State)$, respectively, that satisfy
	\begin{align}	
	\forall \State \in \DomX
	,\quad
	-
	{\zzzLM{\PreliFunc}}
	\leq 
	\PreliFunc(\State) - \LI{\PreliFunc}(\State)
	\leq
	{\zzzUM{\PreliFunc}}
	.\label{eq:def_general_bounds}
	\end{align}	
	Note that the margins ${\zzzLM{\PreliFunc}}$ and ${\zzzUM{\PreliFunc}}$ are not unique.
\end{definition}
\begin{definition}[\MyHighlight{$\ODEFUBpow$th-order property $\MyOrder(\NewDist^{\ODEFUBpow})$}]\label{def:order}
Given $\ODEFUBpow \in \mathbb{N}$, a scalar $\ODEFfunc$ depending on $\allDomSimpX$ is said to be $\ODEFUBpow$th-order and denoted by $\ODEFfunc=\MyOrder(\NewDist^{\ODEFUBpow})$
	if there exists $\ODEFUBDist>0$ and $\ODEFUBCoef>0$ such that for every $\allDomSimpX$ with the corresponding $\NewDist$, we have $\NewDist \leq \ODEFUBDist
	\Rightarrow
	|\ODEFfunc| \leq \ODEFUBCoef \NewDist^{\ODEFUBpow}$.
For example,  $\ODEFfunc=\MyOrder(\NewDist^{2})$ indicate that $|\ODEFfunc|$ decreases quadratically to zero as $\NewDist\to 0$. 
\end{definition}

\tikzstyle{block} = [draw, rectangle, minimum height=2em, minimum width=3em]
\begin{figure}[!t]  
	\begin{center}
		\JFigSimplexes
		\vspace*{-0.1in}
	\end{center}	
	\caption{Two-dimensional illustration of simplexes $\DomSimpX{\IDdom}$ in  $\DomX$.}
	\label{fig:simplexes} 
\end{figure}

\tikzstyle{block} = [draw, rectangle, minimum height=2em, minimum width=3em]
\begin{figure}[!t]  
	\begin{center}		
		\JFigLI
		\vspace*{-0.1in}
	\end{center}
	\caption{One-dimensional illustration of $\PreliFunc(\State)$, its linear interpolation $\LI{\PreliFunc}(\State)$, lower margin ${\zzzLM{\PreliFunc}}$, and upper margin ${\zzzUM{\PreliFunc}}$.}
	\label{fig:LI_margins} 
\end{figure}

\subsection{Summary of solutions to Problem 3: Derinving bounds of $\Basis$}\label{sec_2ndOrderMargin_summary}

Solving Problem 3 starts with the following results.

\begin{proposition}[\MyHighlight{Sampling-based bounds}] \label{def:bounds} 
	For any $\allDomSimpX$, any continuous function $\Basis:\DomBigX \to \mathbb{R}$, any lower margin ${\zzzLM{\Basis}}$, and any upper margin ${\zzzUM{\Basis}}$,
	$\Basis$ is bounded on $\DomX$ as follows:
	\begin{align}	
	&
	\forall \IDdom \in \SetALLIDdom
	,\;
	\forall  \State \in \DomSimpX{\IDdom}
	,\;
	\nonumber\\&
	\min_{\IDPt } \Basis(\SimpPt{\IDdom}{\IDPt}) - {\zzzLM{\Basis}} 
	\leq
	\Basis(\State)
	\leq
	\max_{\IDPt } \Basis(\SimpPt{\IDdom}{\IDPt}) + {\zzzUM{\Basis}} 
	.\label{eq:meaning_LMUM}
	\end{align}		
\end{proposition}
\begin{proof}
This is satisfied because 
$\min_{\IDPt } \Basis(\SimpPt{\IDdom}{\IDPt}) 
\leq \LI{\Basis}(\State) \leq
\max_{\IDPt } \Basis(\SimpPt{\IDdom}{\IDPt})
$ holds for all $\State \in \DomSimpX{\IDdom}$ and all $\IDdom \in \SetALLIDdom$.
\end{proof}

Hence, if ${\zzzLM{\Basis}}$ and $ {\zzzUM{\Basis}} $ are obtained as significantly small values, $\Basis(\State)$ is precisely bounded using a finite number of evaluations of $\Basis(\SimpPt{\IDdom}{\IDPt}) $ at the vertices ${\SimpPt{\IDdom}{\IDPt}}$.
If a stability condition is described by $\Basis(\State)<0$, this condition is checked via a finite number of evaluations and a small upper margin is desirable as shown Fig. \ref{fig:RelatedWork}.
This motivates us to find small margins ${\zzzLM{\Basis}}$ and  ${\zzzUM{\Basis}}$ that depend on the maximum sampling interval $\NewDist$ deeply.

Our main results to derive such small margins are stated below.
For each $\NumNest \in \mathbb{N}$, let ${\FSETfunc{\NumBasis}{\NumNest}}$ be a set of various functions to be defined in Section \ref{sec_set_of_func}.
We propose Algorithm ${\AlgNest{\Basis}{\NumNest}{\allDomSimpX}}$ to be defined in Section \ref{sec_algorithm}.

\begin{theorem}[{\MyHighlight{Second-order margins for ${\FSETfunc{\NumBasis}{\NumNest}}$}}]\label{thm:main_results_margins}
	For any $\allDomSimpX$, any $\NumNest \in \mathbb{N}$, and any function $\Basis \in {\FSETfunc{\NumBasis}{\NumNest}}$,
	Algorithm ${\AlgNest{\Basis}{\NumNest}{\allDomSimpX}}$ provides a lower margin ${\zzzLM{\Basis}}$ and upper margin ${\zzzUM{\Basis}}$ that satisfy the second-order property in Definition \ref{def:order}:
	\begin{align}
		{\zzzLM{\Basis}}&=\MyOrder(\NewDist^{2})
		,\\
		{\zzzUM{\Basis}}&=\MyOrder(\NewDist^{2})
		.
	\end{align}
\end{theorem}
\begin{proof}
	The proof is given in Section \ref{pf:main_results_margins}.		
\end{proof}	

\begin{remark}[{\MyHighlight{Contribution of Theorem \ref{thm:main_results_margins}}}]
By virtue of Theorem \ref{thm:main_results_margins},
we obtain solutions to Problem 3 that are lower and upper bounds
$(\min_{\IDPt } \Basis(\SimpPt{\IDdom}{\IDPt}) - {\zzzLM{\Basis}})$
and
$(\max_{\IDPt } \Basis(\SimpPt{\IDdom}{\IDPt}) + {\zzzUM{\Basis}})$ from Proposition \ref{def:bounds}.
These bounds are precise in the second-order sense, that is,
\begin{align}
\forall \IDdom
,\;
\big|
\big(
\min_{\IDPt } \Basis(\SimpPt{\IDdom}{\IDPt})  - {\zzzLM{\Basis}}  
\big)
-\min_{\State \in \DomSimpX{\IDdom}} \Basis(\State)
\big|
&=  \MyOrder(\NewDist^{2})
,\\
\forall \IDdom 
,\;
\big|
\big(
\max_{\IDPt } \Basis(\SimpPt{\IDdom}{\IDPt}) + {\zzzUM{\Basis}}  
\big)
-\max_{\State \in \DomSimpX{\IDdom}} \Basis(\State) 
\big|
&= \MyOrder(\NewDist^{2})
.
\end{align}	
\end{remark}

In addition, Theorem \ref{thm:main_results_margins} is extended for another set $ {\pFSETfunc{\NumBasis}{\NumNest}}$ of functions and Algorithm ${\AlgpNest{\Basis}{\NumNest}{\allDomSimpX}}$, which will be defined in Section \ref{sec_set_of_func_LB}.

\begin{theorem}[{\MyHighlight{Second-order margins for ${\pFSETfunc{\NumBasis}{\NumNest}}$}}]\label{thm:main_results_margins_specific}
	For any $\allDomSimpX$, any $\NumNest \in \mathbb{N}$, and any $\Basis \in {\pFSETfunc{\NumBasis}{\NumNest}}$,
	Algorithm ${\AlgpNest{\Basis}{\NumNest}{\allDomSimpX}}$
	provides an upper margin ${\zzzUM{\Basis}}$ satisfying the second-order property ${\zzzUM{\Basis}}=\MyOrder(\NewDist^{2})$.
\end{theorem}
\begin{proof}
	The proof is similar to that of Theorem \ref{thm:main_results_margins}.
\end{proof}	

In the following, the details of the set ${\FSETfunc{\NumBasis}{\NumNest}}$, Algorithm ${\AlgNest{\Basis}{\NumNest}{\allDomSimpX}}$, and the proof of Theorem \ref{thm:main_results_margins} are described in Sections \ref{sec_set_of_func}, \ref{sec_algorithm}, and \ref{pf:main_results_margins}, respectively.
Variety of functions in ${\FSETfunc{\NumBasis}{\NumNest}}$  is demonstrated in Section \ref{sec_example_set_of_func}.
Section \ref{sec_set_of_func_LB} presents their extended versions: the set ${\pFSETfunc{\NumBasis}{\NumNest}}$ and Algorithm ${\AlgpNest{\Basis}{\NumNest}{\allDomSimpX}}$.

\textit{\bfseries Tips for readers:}
For briefly understanding this paper, it is easy to follow Sections \ref{sec_solutions} and \ref{sec_demonstration} before Sections \ref{sec_set_of_func}--\ref{sec_set_of_func_LB}.

\subsection{Details of  the set ${\FSETfunc{\NumBasis}{\NumNest}}$}\label{sec_set_of_func}

To describe the sets ${\FSETfunc{\NumBasis}{\NumNest}}$ of functions for $\NumNest \in \mathbb{N}$ focused on in Theorem \ref{thm:main_results_margins}, we use the following definitions. 
\begin{itemize}
	\item 
	Let $\NumBasis \in \mathbb{N}$ be a predefined natural number.

	\item 
	For any continuous function $\yBasis: \DomBigX \to \mathbb{R}$,
	let ${\FSETNLFunc{\yBasis}}$ be the set of functions $\NLFunc:\DomBigY \to \mathbb{R}$
	that satisfy {\ASSNLFA} and {\ASSNLFB}:
	\begin{enumerate}
		\item[\ASSNLFA] 
		The domain $\DomBigY \subseteq \mathbb{R}$ of $\NLFunc$ is an open set containing $\{ \yBasis(\State) \in \mathbb{R} | \State \in \DomBigX \}$.
		\item[\ASSNLFB]
		The function $\NLFunc:\DomBigY \to \mathbb{R}$ is $C^{2}$ continuous and
		for some $\Globalddsup>0$,
		for any $\allDomSimpX$,
		there exist \textit{known}\footnote{%
		This means that these constants are obtained from the form of the corresponding function (e.g., $\NLFunc$ or $\PreliFunc$) and given information.}
		bounds ${\dinf{\NLFunc}}$, ${\dsup{\NLFunc}}$, ${\ddinf{\NLFunc}}$, and ${\ddsup{\NLFunc}}$ that satisfy
		$|{\ddinf{\NLFunc}}|+|{\ddsup{\NLFunc}}| < \Globalddsup$ and
		\begin{align}
		\forall \tempyBasis \in \DomY,\;
		{\dinf{\NLFunc}}
		&\leq 
		\frac{\partial \NLFunc(\tempyBasis) }{\partial \yBasis} 
		\leq
		{\dsup{\NLFunc}}
		,\label{eq:ASSNLFB_eq1}
		\\
		\forall \tempyBasis \in \DomLIY,\;
		{\ddinf{\NLFunc}}
		&\leq 
		\frac{\partial^{2} \NLFunc(\tempyBasis) }{\partial \yBasis^{2}} 
		\leq
		{\ddsup{\NLFunc}}
		,\label{eq:ASSNLFB_eq2}
		\end{align}
		where $\DomY:=\{ {\yBasis}(\State) \in \mathbb{R} | \State \in \DomX \}$ and
		$\DomLIY:=\{ \LI{\yBasis}(\State) \in \mathbb{R} | \State \in \DomX \}$.
	\end{enumerate}

\end{itemize}

\begin{proposition}[\MyHighlight{Examples of  ${\FSETNLFunc{\yBasis}}$}]\label{thm:ex_classH}
	For any continuous function $\yBasis: \DomBigX \to \mathbb{R}$, 
	the functions
	$\sin(\yBasis)$,
	$\cos(\yBasis)$,
	$1/(1+\exp(-\yBasis))$ (sigmoid), 
	and
	$\tanh \yBasis = (\exp(2\yBasis)-1)/(\exp(2\yBasis)+1)$
	are contained in $ {\FSETNLFunc{\yBasis}}$. 	
\end{proposition}
\begin{proof}
	The proof is given in Appendix \ref{pf:ex_classH}.		
\end{proof}

\begin{definition}[\MyHighlight{Basis set ${\FSETfunc{\NumBasis}{1}}$}]\label{def:basis_set}
	Let ${\FSETfunc{\NumBasis}{1}}$ be {the} set of $C^{2}$ continuous functions $\PreliFunc: \DomBigX \to \mathbb{R}$ such that there exist \textit{known}$^{1}$ bounds ${\dirddinf{\PreliFunc}}$ and ${\dirddsup{\PreliFunc}}$ satisfying
	\begin{align}
		&
		\forall
		\State \in \DomX
		,\;
		\forall
		\tempUnitVec 
		\in \{  \tempUnitVec \in  \mathbb{R}^{\DimX} 
		|
		\|\tempUnitVec\|=1
		\}
		,\;
		\nonumber\\&
		{\dirddinf{\PreliFunc}}
		\leq
		\frac{ \partial^{2} \PreliFunc(  \State + \tempDist\tempUnitVec  ) }{\partial \tempDist^{2}}
		\Big|_{ \tempDist=0}
		\leq
		{\dirddsup{\PreliFunc}}
		.\label{eq:bounds_of_directional_derivative}
	\end{align}
\end{definition}

\begin{remark}\label{rem:bounds_of_directional_derivative}
Note that ${\dirddinf{\PreliFunc}}$ and ${\dirddsup{\PreliFunc}}$ represent bounds of second-order directional derivatives.
If the minimum and maximum eigenvalues of 
${\partial^{2}  \PreliFunc( \tempState  ) }/{\partial  \State \State^{\MyTRANSPO}} $
are bounded on $\DomX$, they can be set as ${\dirddinf{\PreliFunc}}$ and ${\dirddsup{\PreliFunc}}$, respectively.
Specifically, by setting $\tempState(\tempDist) = \State+\tempDist \tempUnitVec $ with $\partial \tempState(\tempDist)/\partial \tempDist = \tempUnitVec$, we obtain
\begin{align}
\frac{ \partial^{2} \PreliFunc(   \tempState(0)  ) }{\partial \tempDist^{2}}
&=
\tempUnitVec^{\MyTRANSPO}  
\frac{\partial^{2}  \PreliFunc( \tempState(0)  ) }{\partial  \State \State^{\MyTRANSPO}} 
\tempUnitVec
\leq
\sup_{ \tempState \in \DomX }
\EigMax[\Big]{ 
	\frac{\partial^{2}  \PreliFunc( \tempState  ) }{\partial  \State \State^{\MyTRANSPO}} 
}
.
\end{align}		
\end{remark}

\begin{definition}[\MyHighlight{Nested set ${\FSETfunc{\NumBasis}{\NumNest}}$}]\label{def:nested_set}
	For each $\NumNest \in \{2,3,\dots\}$, ${\FSETfunc{\NumBasis}{\NumNest}}$ denotes the set of functions $\Basis: \DomBigX \to \mathbb{R}$
	such that there exist \textit{known}\footnote{%
		This means that if $\Basis$ belongs to this set, we can identify the corresponding ${\subAbasis{\IDEl}}$, ${\subBbasis{\IDEl}}$, and $\NLFunc$ that satisfy the conditions.}
	functions 
	${\subAbasis{\IDEl}}\in {\FSETfunc{\NumBasis}{\NumNest-1}}$, 
	${\subBbasis{\IDEl}}\in {\FSETfunc{\NumBasis}{\NumNest-1}}$ for $\IDEl \in \{1,2,\dots,\NumBasis \}$,
	and 
	$\NLFunc \in {\FSETNLFunc{\yBasis}}$ 
	that satisfy the following relations:
	\begin{align}
		\Basis(\State)
		&=  
		\NLFunc 
		(\yBasis(\State))
		,  \label{eq:def_Basis}
		\\
		\yBasis(\State)
		&= 
		\sum_{\IDEl=1}^{\NumBasis} 
		{\subAbasis{\IDEl}}(\State) {\subBbasis{\IDEl}}(\State) 
		. \label{eq:def_yBasis}
	\end{align}
\end{definition}

\begin{proposition}[\MyHighlight{Key properties of ${\FSETfunc{\NumBasis}{\NumNest}}$}]\label{thm:FSETfunc_properties}
For any $\NumNest \in \mathbb{N}$, the following properties hold.
\begin{enumerate}
	\item\label{thm:C2continuous}
	Any $\Basis \in {\FSETfunc{\NumBasis}{\NumNest}}$ is $C^{2}$ continuous on $\DomBigX$.

	\item\label{thm:FSETfunc_nondec}
	We have ${\FSETfunc{\NumBasis}{\NumNest}} \subseteq {\FSETfunc{\NumBasis}{\NumNest+1}}$.
	
\end{enumerate}	
\end{proposition}
\begin{proof}
The proof is given in Appendix \ref{pf:FSETfunc_properties}.		
\end{proof}

\subsection{Examples of functions contained in the set ${\FSETfunc{\NumBasis}{\NumNest}}$}\label{sec_example_set_of_func}

We demonstrate the applicability of the proposed method.
It is firstly shown that some basic functions are contained in ${\FSETfunc{\NumBasis}{\NumNest}}$. 

\begin{theorem}[\MyHighlight{Basic functions in ${\FSETfunc{\NumBasis}{\NumNest}}$}] \label{thm:ex_simple_functions}
	$\;$
	\begin{enumerate}
		
		\item
		Any quadratic function $\PreliFunc(\State)=\State^{\MyTRANSPO} \QFcoefA \State + \State^{\MyTRANSPO} \QFcoefB + \QFcoefC $ with symmetric $\QFcoefA$ is contained in the basis set ${\FSETfunc{\NumBasis}{1}}$ and the settings of ${\dirddinf{\PreliFunc}}={\EigMin{ \QFcoefA }}$ and ${\dirddsup{\PreliFunc}}={\EigMax{ \QFcoefA }}$  satisfy the condition \eqref{eq:bounds_of_directional_derivative}.
		If $\PreliFunc(\State)$ is constant or linear with $\QFcoefA=0$, we have ${\dirddinf{\PreliFunc}}={\dirddsup{\PreliFunc}}=0$.

		\item
		Suppose $\NumBasis \geq \DimX +1$.
		For any $\NumNest \in \{3,4,\dots\}$,
		any $\NumNest$th order polynomial $\PreliFunc: \DomBigX \to \mathbb{R}$ is contained in ${\FSETfunc{\NumBasis}{\NumNest-1}}$.
		
		\item
		Given $\NumNest \in \{2,3,\dots\}$, ${\LCombiCoef{\IDEl}} \in \mathbb{R}$, and ${\subAbasis{\IDEl}} \in {\FSETfunc{\NumBasis}{\NumNest-1}}$,
		the linear combination $\PreliFunc(\State)=  
		\sum_{\IDEl=1}^{\NumBasis} 
		{\LCombiCoef{\IDEl}}
		{\subAbasis{\IDEl}}(\State)$  (whose structure is given) is contained in ${\FSETfunc{\NumBasis}{\NumNest}}$. 
	\end{enumerate}		
\end{theorem}
\begin{proof}
	The proof is given in Appendix \ref{pf:ex_simple_functions}.
\end{proof}

Next, let us define the following squared exponential (SE) kernel ${\kernel{  \wildcard   }{\DState{\IDdat}}}:\DomBigX\to\mathbb{R}$ as follows:
\begin{align}	
	\kernel{  \State   }{\DState{\IDdat}}
	:=	\hypMag \exp ( - ( \State -\DState{\IDdat})^{\MyTRANSPO} \hypCovMat^{-1}  (  \State -\DState{\IDdat})  /2 )	
	,\label{eq:def_SEkernel}
\end{align}
where $\hypMag > 0$ and a positive definite symmetric  matrix $\hypCovMat \in \mathbb{R}^{\DimX \times \DimX}$ are hyperparameters, and ${\DState{\IDdat}} \in \mathbb{R}^{\DimX}$ denotes a data point.
The SE kernel is popular and works well in data-driven approaches for system identification and control design.
We show that the SE kernel is contained in ${\FSETfunc{\NumBasis}{1}}$.
\begin{theorem}[\MyHighlight{Squared exponential kernels in ${\FSETfunc{\NumBasis}{1}}$}] \label{thm:ex_SEkernel}	
	The SE kernel ${\kernel{  \wildcard   }{\DState{\IDdat}}}$ in \eqref{eq:def_SEkernel} is contained in the basis set ${\FSETfunc{\NumBasis}{1}}$ and the following settings satisfy the condition \eqref{eq:bounds_of_directional_derivative}:	
	\begin{align}	
		{\dirddinf{\kernel{  \wildcard   }{\DState{\IDdat}}}}
		&=
		-\hypMag {\EigMax{\hypCovMat^{-1}}} 
		,\label{eq:SEkernel_dirddinf}
		\\
		{\dirddsup{\kernel{  \wildcard   }{\DState{\IDdat}}}}
		&=
		2\hypMag {\EigMax{\hypCovMat^{-1}}} \exp(-3/2) 
		.\label{eq:SEkernel_dirddsup}
	\end{align}		
\end{theorem}
\begin{proof}
	The proof is given in Appendix \ref{pf:ex_SEkernel}.	
\end{proof}

Next, consider the following GP \cite{Rasmussen06} with posterior mean $\GPmean:\DomBigX\to\mathbb{R}^{\DimX}$ and standard deviation $\GPsd:\DomBigX\to\mathbb{R}^{\DimX}$ that is developed using $\NumData$ data points:
\begin{align}	
	\GPmean(\State)
	&:=
	[\DDrift{1}, \DDrift{2}, \dots, \DDrift{\NumData}]
	\kernelMat^{-1}\kernelVec(\State)
	,\label{eq:def_GPmean}
	\\
	{\El{\GPsd(\State)}{\IDEl}}&:=\sqrt{ \kernel{\State}{\State} - \kernelVec(\State)^{\MyTRANSPO} \kernelMat^{-1} \kernelVec(\State)  } 
	,\label{eq:def_GPsd}
	\\
	\kernelVec(\State)
	&:= [\kernel{\State}{\DState{1}}, \dots, \kernel{\State}{\DState{\NumData}}]^{\MyTRANSPO} \in  \mathbb{R}^{\NumData} 
	,\label{eq:def_kernel_vec}
	\\	
	{\El{\kernelMat}{\IDdat,\IDbdat}}
	&:=
	\kernel{\DState{\IDdat}}{\DState{\IDbdat}}
	+ \KroDelta{\IDdat}{\IDbdat} \hypNoise
	,\label{eq:def_kernel_mat_1}
\end{align}
for $\IDEl \in \{1,\dots, \DimX\}$, $\IDdat \in \{1,\dots, \NumData\}$, and $\IDbdat \in \{1,\dots, \NumData\}$.
The symbols ${\DState{\IDdat}} \in \mathbb{R}^{\DimX}$ and ${\DDrift{\IDdat}} \in \mathbb{R}^{\DimX}$ denote the input and output of $\IDdat$th data point generated from a function with noise to be predicted, respectively. 
The symbols $\KroDelta{\IDdat}{\IDbdat}$ and $\hypNoise> 0$ are the Kronecker delta and the hyperparameter, respectively.

\begin{theorem}[\MyHighlight{Mean of Gaussian processes of ${\FSETfunc{\NumBasis}{2}}$}] \label{thm:ex_GP}	
	Suppose 
	that $\NumBasis \geq \NumData$ holds and 
	that $\kernel{ \wildcard }{\DState{\IDdat}}$ is contained in ${\FSETfunc{\NumBasis}{1}}$ or is the SE kernel in \eqref{eq:def_SEkernel}.
	Each component of the GP mean $\GPmean$ defined in \eqref{eq:def_GPmean} is contained in ${\FSETfunc{\NumBasis}{2}}$.
\end{theorem}
\begin{proof}
	The proof is given in Appendix \ref{pf:ex_GP}.
\end{proof}
\begin{remark}[\MyHighlight{Standard deviation of GPs}]
The standard deviation $\GPsd$ in \eqref{eq:def_GPsd} is analyzed in Theorem \ref{thm:ex_GPsd} in Section \ref{sec_set_of_func_LB}.	
\end{remark}

Finally, consider the following DNN \cite[Section 1.2]{Aggarwal18}:
\begin{align}
	{\El{\NNnode{\NNidlayer+1}(\State)}{\IDEl}}
	&=\NLFunc( {\NNweight{\NNidlayer}{\IDEl}^{\MyTRANSPO}} {\NNnode{\NNidlayer}}(\State)   ) 
	,\label{eq:def_NN}
\end{align}
for $ \NNidlayer \in \{1,\dots   \NNnumlayer-1  \}$ and $ \IDEl \in \{1,\dots   {\NNdim{\NNidlayer+1}}  \}$,
where ${\NNnode{1}}(\State):= [\State^{\MyTRANSPO},1]^{\MyTRANSPO}$, ${\NNnode{\NNidlayer}}(\State)$ for $\NNidlayer \in \{2,\dots   \NNnumlayer-1  \}$, and ${\NNnode{\NNnumlayer}}(\State)$ are 
the nodes in the input layer, hidden layer, and output layer, respectively, 
and ${\NNdim{\NNidlayer}}\in \mathbb{N}$ denotes the dimension of  $\NNidlayer$th  layer.
The symbols ${\NNweight{\NNidlayer}{\IDEl}} \in \mathbb{R}^{{\NNdim{\NNidlayer}}}$ and $\NLFunc$ indicate the weight vector and activation function, respectively.

\begin{theorem}[\MyHighlight{Deep neural networks in ${\FSETfunc{\NumBasis}{\NNnumlayer}}$}] \label{thm:ex_NN}	
	Given an activation function $\NLFunc$ in \eqref{eq:def_NN}, suppose that $\NLFunc \in {\FSETNLFunc{\yBasis}}$ holds for any continuous function $\yBasis: \DomBigX \to \mathbb{R}$.
	For any $\NNnumlayer$ and any ${\NNdim{\NNidlayer}}$,
	suppose that $\NumBasis \geq {\NNdim{\NNidlayer}}$ holds for any $\NNidlayer \in \{1,\dots   \NNnumlayer-1  \}$.  
	Each component of the output node ${\NNnode{\NNnumlayer}}$ of the DNN in \eqref{eq:def_NN} is contained in the set ${\FSETfunc{\NumBasis}{\NNnumlayer}}$.	
\end{theorem}
\begin{proof}
	The proof is given in Appendix \ref{pf:ex_NN}.	
\end{proof}
\begin{remark}
	Various functions such as the sigmoid activation can be adopted as activation functions $\NLFunc$ by Proposition \ref{thm:ex_classH}. 	
\end{remark}

\subsection{Details of Algorithm ${\AlgNest{\Basis}{\NumNest}{\allDomSimpX}}$}	\label{sec_algorithm}

This subsection presents Algorithm ${\AlgNest{\Basis}{\NumNest}{\allDomSimpX}}$ associated with some key lemmas and the following constants.
\begin{definition}[\MyHighlight{Key constants}] \label{def:collection} 
	For any $\allDomSimpX$, any continuous $\PreliFunc: \DomBigX \to \mathbb{R}$,
	any $ {\zzzLM{\PreliFunc}} $, and any $ {\zzzUM{\PreliFunc}}$, let us define the following constants if they exist.
	\begin{enumerate}
		\item 	
		Let $\zzzLB{\PreliFunc}$ and $\zzzUB{\PreliFunc}$ be the lower and upper bounds of $\PreliFunc$:
		\begin{align}	
		{\zzzLB{\PreliFunc}}
		&:= \min_{\IDdom \in \SetALLIDdom,\IDPt} \PreliFunc(\SimpPt{\IDdom}{\IDPt}) - {\zzzLM{\PreliFunc}} 
		,
		\\
		{\zzzUB{\PreliFunc}}
		&:= \max_{\IDdom \in \SetALLIDdom,\IDPt} \PreliFunc(\SimpPt{\IDdom}{\IDPt}) + {\zzzUM{\PreliFunc}} 
		.
		\end{align}

		\item
		Let $\zzzLipConst{\PreliFunc}$ be the sampling-based gradient of $\PreliFunc$:
		\begin{align}
		{\zzzLipConst{\PreliFunc}}
			&:=
			\max_{\IDdom \in \SetALLIDdom,\IDPt,\IDbPt}
			| \PreliFunc(\SimpPt{\IDdom}{\IDPt}) - \PreliFunc(\SimpPt{\IDdom}{\IDbPt}) |
			/\NewDist
			,
		\end{align}
		\item
		The collection of the constants are denoted by 	
		\begin{align}	
			{\zzEachMBGset{\PreliFunc}}
			&:=
			[
			\zzzLM{\PreliFunc}, \zzzUM{\PreliFunc}, \zzzLB{\PreliFunc}, \zzzUB{\PreliFunc}, \zzzLipConst{\PreliFunc}   ]^{\MyTRANSPO}
			.
		\end{align}	
	\end{enumerate}	 		
\end{definition}

\begin{remark}[\MyHighlight{Boundedness independent of $\NewDist$}]\label{rem:collection}
If $ {\zzzLM{\PreliFunc}} $ and $ {\zzzUM{\PreliFunc}}$ are $\ODEFUBpow$th-order  $\MyOrder(\NewDist^{\ODEFUBpow})$ for some $\ODEFUBDist$ and $\ODEFUBpow$, then $\zzzLB{\PreliFunc}$ and $\zzzUB{\PreliFunc}$ are clearly bounded on $\DomX$ regardless of $\NewDist \in (0, \ODEFUBDist]$. 
The sampling-based gradient ${\zzzLipConst{\PreliFunc}}$ is bounded regardless of $\allDomSimpX$ if $\PreliFunc$ is $C^{2}$ continuous on the bounded closed set $\DomX$, that is, Lipschitz continuous on $\DomX$.  
These bounded properties will be used in the proofs of Lemmas \ref{thm:margins_for_inner_products}, \ref{thm:margins_for_nonlinear_mappings}, and \ref{thm:margins_for_inner_products_specific}.
\end{remark}

Based on these definitions, we propose Algorithm ${\AlgNest{\Basis}{\NumNest}{\allDomSimpX}}$ that provides ${\zzEachMBGset{\Basis}}$ for any $\Basis \in {\FSETfunc{\NumBasis}{\NumNest}}$.
After decomposing $\Basis$ in Line 1,
we obtain ${\zzEachMBGset{\Basis}}$ of $\Basis \in {\FSETfunc{\NumBasis}{\NumNest}}$ using that of ${\subAbasis{\IDEl}}\in {\FSETfunc{\NumBasis}{\NumNest-1}}$ and ${\subBbasis{\IDEl}}\in {\FSETfunc{\NumBasis}{\NumNest-1}}$.
Lines 2--8 indicate a nested structure of this algorithm.
Specifically, ${\zzEachMBGset{\subAbasis{\IDEl}}}$ and ${\zzEachMBGset{\subBbasis{\IDEl}}}$ are calculated by using Algorithms ${\AlgNest{\subAbasis{\IDEl}}{\NumNest-1}{\allDomSimpX}}$ and ${\AlgNest{\subBbasis{\IDEl}}{\NumNest-1}{\allDomSimpX}}$, respectively.
Such a nested process continues until ${\subAbasis{\IDEl}}$ and ${\subBbasis{\IDEl}}$ are contained in ${\FSETfunc{\NumBasis}{1}}$.
For every ${\subAbasis{\IDEl}}$ and ${\subBbasis{\IDEl}}$ in ${\FSETfunc{\NumBasis}{1}}$, we obtain ${\zzEachMBGset{\subAbasis{\IDEl}}}$ and ${\zzEachMBGset{\subBbasis{\IDEl}}}$  in Line 4 using Lemma \ref{thm:margin_FSET1} that is presented below.
After the nested process, 
we finally obtain ${\zzEachMBGset{\Basis}}$ in Line 9--11, using Lemmas \ref{thm:margins_for_inner_products} and \ref{thm:margins_for_nonlinear_mappings} presented below.

\begin{algorithm}[!t]    
	\renewcommand{\algorithmicrequire}{\textbf{Input:}}
	\renewcommand{\algorithmicensure}{\textbf{Output:}}                       
	\caption{${\AlgNest{\Basis}{\NumNest}{\allDomSimpX}}$ }  

	\label{alg:nested_margins}
	\begin{algorithmic}[1]                  
		\REQUIRE $\Basis \in {\FSETfunc{\NumBasis}{\NumNest}}$, $\NumNest\geq 1$, and $\allDomSimpX$
		\ENSURE ${\zzEachMBGset{\Basis}}=[
		\zzzLM{\Basis}, \zzzUM{\Basis}, \zzzLB{\Basis}, \zzzUB{\Basis}, \zzzLipConst{\Basis}   ]^{\MyTRANSPO}$

		\STATE
		Decompose $\Basis\in {\FSETfunc{\NumBasis}{\NumNest}}$ into ${\subAbasis{\IDEl}}, {\subBbasis{\IDEl}}\in {\FSETfunc{\NumBasis}{\NumNest-1}}$, for $\IDEl \in \{1,2,\dots,\NumBasis \}$, and $\NLFunc$
		according to Definition \ref{def:nested_set}
		
		\FOR{all $\PreliFunc \in \{\subAbasis{1},\dots,\subAbasis{\NumBasis},\subBbasis{1},\dots,\subBbasis{\NumBasis}   \} $  }	
		
		\IF{ $\PreliFunc \in {\FSETfunc{\NumBasis}{1}} $ } 
		\STATE
		Obtain ${\zzEachMBGset{ \PreliFunc }}$ by Lemma \ref{thm:margin_FSET1}
		and Definition \ref{def:collection}
		\ELSE	
		\STATE
		Calculate ${\zzEachMBGset{ \PreliFunc }}$ by Algorithm	${\AlgNest{\PreliFunc}{\NumNest-1}{\allDomSimpX}}$
		\ENDIF

		\ENDFOR

		\STATE
		Obtain $({\zzzLM{\yBasis}}, {\zzzUM{\yBasis}})$ by Lemma \ref{thm:margins_for_inner_products} using $({\zzEachMBGset{\subAbasis{\IDEl}}}, {\zzEachMBGset{\subBbasis{\IDEl}}})$
		\STATE
		Calculate $({\zzzLM{\Basis}},{\zzzUM{\Basis}})$ by Lemma \ref{thm:margins_for_nonlinear_mappings} using $({\zzzLM{\yBasis}}, {\zzzUM{\yBasis}},{\zzzLipConst{\yBasis}})$, where ${\zzzLipConst{\yBasis}}$ is given by Definition \ref{def:collection}
		
		\STATE
		Obtain ${\zzzLB{\Basis}}$, ${\zzzUB{\Basis}}$, and  ${\zzzLipConst{\Basis}}$ 
		by Definition \ref{def:collection}
		
		\RETURN

	\end{algorithmic}
\end{algorithm}

\begin{lemma}[\MyHighlight{Second-order margins for ${\FSETfunc{\NumBasis}{1}}$}]\label{thm:margin_FSET1}
	For any $\PreliFunc \in {\FSETfunc{\NumBasis}{1}}$,
	second-order margins ${\zzzLM{\PreliFunc}}$ and ${\zzzUM{\PreliFunc}}$ are given by
	\begin{align}	
		{\zzzLM{\PreliFunc}}
		&=
		({ \DimX \NewDist^{2} }/{8})
		\max\{ 0,		{\dirddsup{\PreliFunc}}	\}
		=\MyOrder(\NewDist^{2})
		,\label{eq:margin_FSET1_LM}
		\\
		{\zzzUM{\PreliFunc}}
		&=
		({ \DimX \NewDist^{2} }/{8})
		\max\{ 0,	-{\dirddinf{\PreliFunc}}	\}
		=\MyOrder(\NewDist^{2})
		.\label{eq:margin_FSET1_UM}
	\end{align}	
\end{lemma}
\begin{proof}
	The proof is given in Appendix \ref{pf:margin_FSET1}.	
\end{proof}	

\begin{lemma}[\MyHighlight{Margins for inner products}] \label{thm:margins_for_inner_products}
	Given any $\NumNest \geq 1$,
	any ${\subAbasis{\IDEl}} \in {\FSETfunc{\NumBasis}{\NumNest}}$,
	and any ${\subBbasis{\IDEl}} \in {\FSETfunc{\NumBasis}{\NumNest}}$ for $\IDEl \in \{1,2,\dots,\NumBasis \}$,
	a lower and an upper margin of $\yBasis$ in \eqref{eq:def_yBasis} are given by	
	\begin{align}	
		{\zzzLM{\yBasis}}
		&=		
		\sum_{\IDEl=1}^{\NumBasis}
		\Big(
		{\zzzBLB{\subAbasis{\IDEl}}{\subBbasis{\IDEl}}}
		+
		{\zzzBLB{\subBbasis{\IDEl}}{\subAbasis{\IDEl}}}
		+
		\NewDist^{2} 	{\zzzLipConst{\subAbasis{\IDEl}}}{\zzzLipConst{\subBbasis{\IDEl}}}
		\Big)
		,\label{eq:chi_bound_L}
		\\
		{\zzzUM{\yBasis}}
		&=
		\sum_{\IDEl=1}^{\NumBasis}
		\Big(
		{\zzzBUB{\subAbasis{\IDEl}}{\subBbasis{\IDEl}}}
		+
		{\zzzBUB{\subBbasis{\IDEl}}{\subAbasis{\IDEl}}}
		+
		\NewDist^{2} 	{\zzzLipConst{\subAbasis{\IDEl}}}{\zzzLipConst{\subBbasis{\IDEl}}}
		\Big)
		,\label{eq:chi_bound_U}			
	\end{align}	
	where
	\begin{align}
	{\zzzBLB{\PreliAFunc}{\PreliBFunc}}
	&:=
	{\zzzrawBLB{\PreliAFunc}{\PreliBFunc}}
	,\label{eq:@multi_func_LB}
	\\
	{\zzzBUB{\PreliAFunc}{\PreliBFunc}}
	&:=
	{\zzzrawBUB{\PreliAFunc}{\PreliBFunc}}
	.\label{eq:@multi_func_UB}	
	\end{align}	
	Moreover, these ${\zzzLM{\yBasis}}$ and  ${\zzzUM{\yBasis}}$ are $\MyOrder(\NewDist^{2})$ if $\zzzLM{\subAbasis{\IDEl}}$, $\zzzUM{\subAbasis{\IDEl}}$, $\zzzLM{\subBbasis{\IDEl}}$, and $\zzzUM{\subBbasis{\IDEl}}$ for $\IDEl \in \{1,\dots,\NumBasis\}$ are $\MyOrder(\NewDist^{2})$.
\end{lemma}
\begin{proof}
	The proof is given in Appendix \ref{pf:margins_for_inner_products}.			
\end{proof}	

\begin{lemma}[\MyHighlight{Margins for nonlinear mappings}]  \label{thm:margins_for_nonlinear_mappings}
	Given any $C^{2}$ continuous function $\yBasis:\DomBigX \to \mathbb{R}$,
	any $\NLFunc \in {\FSETNLFunc{\yBasis}}$, 
	and 
	any $\Basis(\State)$ expressed by \eqref{eq:def_Basis},
	a lower and an upper margin of $\Basis$ are given by
	\begin{align}	
		{\zzzLM{\Basis}}
		&=		
		\max \big\{ 
		-
		\zzzUM{\yBasis}
		{\dinf{\NLFunc}} 
		,    
		\zzzLM{\yBasis}  
		{\dsup{\NLFunc}} 
		\big\}
		\nonumber\\&\quad
		+
		\NewDist^{2}
		\frac{ \DimX {\zzzLipConst{\yBasis}}^{2}}{8}
		\max\{0, {\ddsup{\NLFunc}} \}
		,\label{eq:LM_Basis_via_NLFunc}		
		\\
		{\zzzUM{\Basis}}
		&=
		\max \big\{ 
		-
		\zzzLM{\yBasis}
		{\dinf{\NLFunc}} 
		,    
		\zzzUM{\yBasis}  
		{\dsup{\NLFunc}} 
		\big\}
		\nonumber\\&\quad
		+
		\NewDist^{2}
		\frac{ \DimX  {\zzzLipConst{\yBasis}}^{2}}{8}
		\max\{0, -{\ddinf{\NLFunc}} \}
		,\label{eq:UM_Basis_via_NLFunc}			
	\end{align}		
	Moreover, these ${\zzzLM{\Basis}}$ and  ${\zzzUM{\Basis}}$ are second-order $\MyOrder(\NewDist^{2})$ if  ${\zzzLM{\yBasis}}$ and  ${\zzzUM{\yBasis}}$ are $\MyOrder(\NewDist^{2})$.
\end{lemma}
\begin{proof}
	The proof is given in Appendix \ref{pf:margins_for_nonlinear_mappings}.		
\end{proof}

\subsection{Proof of Theorem \ref{thm:main_results_margins}} \label{pf:main_results_margins}

We prove Theorem \ref{thm:main_results_margins} using mathematical induction.   
Because of Lemma \ref{thm:margin_FSET1}, any $\PreliFunc \in {\FSETfunc{\NumBasis}{1}}$ satisfies 
${\zzzLM{\PreliFunc}}=\MyOrder(\NewDist^{2})$
and 
${\zzzUM{\PreliFunc}}=\MyOrder(\NewDist^{2})$.
For an $\NumNest\geq 1$, suppose that 
for every $\PreliFunc \in {\FSETfunc{\NumBasis}{\NumNest}}$,
Algorithm ${\AlgNest{\PreliFunc}{\NumNest}{\allDomSimpX}}$
obtains 
${\zzzLM{\PreliFunc}}=\MyOrder(\NewDist^{2})$
and 
${\zzzUM{\PreliFunc}}=\MyOrder(\NewDist^{2})$.
Every $\Basis \in {\FSETfunc{\NumBasis}{\NumNest+1}}$ can be decomposed into ${\subAbasis{\IDEl}}, {\subBbasis{\IDEl}}\in {\FSETfunc{\NumBasis}{\NumNest}}$, and $\NLFunc$ by Line 1 of ${\AlgNest{\Basis}{\NumNest+1}{\allDomSimpX}}$.
Then, we obtain ${\zzEachMBGset{ \PreliFunc }}$ for every $\PreliFunc \in \{\subAbasis{1},\dots,\subAbasis{\NumBasis},\subBbasis{1},\dots,\subBbasis{\NumBasis}   \} $.
Next, 
Lemma \ref{thm:margins_for_inner_products} provides
${\zzzLM{\yBasis}}=\MyOrder(\NewDist^{2})$
and 
${\zzzUM{\yBasis}}=\MyOrder(\NewDist^{2})$.
Subsequently, Lemma \ref{thm:margins_for_nonlinear_mappings} provides
${\zzzLM{\Basis}}=\MyOrder(\NewDist^{2})$
and 
${\zzzUM{\Basis}}=\MyOrder(\NewDist^{2})$ for $\Basis \in {\FSETfunc{\NumBasis}{\NumNest+1}}$.
Thus,
for every $\PreliFunc \in {\FSETfunc{\NumBasis}{\NumNest+1}}$,
Algorithm ${\AlgNest{\PreliFunc}{\NumNest+1}{\allDomSimpX}}$
obtains 
${\zzzLM{\PreliFunc}}=\MyOrder(\NewDist^{2})$
and 
${\zzzUM{\PreliFunc}}=\MyOrder(\NewDist^{2})$.
Because of the mathematical induction, this statement holds for any $\NumNest$.
This completes the proof.

\subsection{Details of  ${\pFSETfunc{\NumBasis}{\NumNest}}$ and Algorithm ${\AlgpNest{\Basis}{\NumNest}{\allDomSimpX}}$}\label{sec_set_of_func_LB}

To describe the details of Theorem \ref{thm:main_results_margins_specific} presented in Section \ref{sec_2ndOrderMargin_summary},
we propose a modified version of the set ${\FSETfunc{\NumBasis}{\NumNest}}$ in this subsection.
Some definitions are introduced below.
\begin{itemize}
	\item
	Let $\NonNegativeFSETNLFunc$ be the set of nonnegative functions $\PreliFunc: \DomBigX \to \NNReal$.

	\item 	
	For any continuous function $\yBasis: \DomBigX \to \mathbb{R}$,
	let ${\pAFSETNLFunc{\yBasis}}$ (resp. ${\pBFSETNLFunc{\yBasis}}$) 	
	be the set of functions $\NLFunc:\DomBigY \to \mathbb{R}$
	satisfying {\ASSNLFA} and the following {\ASSNLFC} (resp. {\ASSNLFD}):	
\begin{enumerate}
	\item[\ASSNLFC]
	The function $\NLFunc:\DomBigY \to \mathbb{R}$ is monotonically nondecreasing $C^{2}$ continuous and
	for some $\Globalddsup>0$,
	for any $\allDomSimpX$,
	there exist \textit{known}$^{1}$ bounds 
	${\ddinf{\NLFunc}}$ and ${\dsup{\NLFunc}}$ that satisfy
	 \eqref{eq:ASSNLFB_eq2}, $|{\ddinf{\NLFunc}}|+|{\dsup{\NLFunc}}| < \Globalddsup$, and
	\begin{align}
	\forall \tempyBasis \in
	[\min_{\tempyBasis \in \DomLIY} \tempyBasis,
	 \max_{\tempyBasis \in \DomY}   \tempyBasis ]
	,\;
	\frac{\partial \NLFunc(\tempyBasis) }{\partial \yBasis} 
	\leq
	{\dsup{\NLFunc}}
	.\label{eq:ASSNLFB_eq1mod}
	\end{align}
	\item[\ASSNLFD]
	The function $\NLFunc:\DomBigY \to \mathbb{R}$ is convex and Lipschitz continuous with a \textit{known}$^{1}$ Lipschitz constant  $\LipsNLFunc\geq 0$, where the differentiability of $\NLFunc$ is not needed.
\end{enumerate}

	\item
	Let
	the part of ${\zzEachMBGset{\PreliFunc}}$
	be
	${\zzpEachMBGset{\PreliFunc}}
	:=
	[
	 \zzzUM{\PreliFunc}, \zzzUB{\PreliFunc}, \zzzLipConst{\PreliFunc}   ]^{\MyTRANSPO}
	$.

\end{itemize}

We define the set ${\pFSETfunc{\NumBasis}{\NumNest}}$ with its examples and properties.

\begin{definition}[{\MyHighlight{Nonnegative nested set ${\pFSETfunc{\NumBasis}{\NumNest}}$}}]\label{def:nonnegative_nested_set}	
	Let us define 
	${\pFSETfunc{\NumBasis}{1}}
	:=
	{\FSETfunc{\NumBasis}{1}} \cap \NonNegativeFSETNLFunc
	$.
	For $\NumNest \in \{2,3,\dots\}$, the following sets are defined.
\begin{enumerate}
	\item 
	Let ${\pAFSETfunc{\NumBasis}{\NumNest}}$ be the set of nonnegative functions $\Basis: \DomBigX \to \NNReal$ 
	such that there exist \textit{known}$^{2}$  functions 
	${\subAbasis{\IDEl}}\in {\pFSETfunc{\NumBasis}{\NumNest-1}}$, 
	${\subBbasis{\IDEl}}\in {\pFSETfunc{\NumBasis}{\NumNest-1}}$ for $\IDEl \in \{1,2,\dots,\NumBasis \}$,
	and 
	$\NLFunc \in {\pAFSETNLFunc{\yBasis}}$ 
	that satisfy \eqref{eq:def_Basis} and \eqref{eq:def_yBasis}.

	\item
	Let ${\pBFSETfunc{\NumBasis}{\NumNest}}$ be the set of nonnegative functions $\Basis: \DomBigX \to \NNReal$ 
	such that there exist \textit{known}$^{2}$  functions 
	${\subAbasis{\IDEl}}\in {\FSETfunc{\NumBasis}{\NumNest-1}}$, 
	${\subBbasis{\IDEl}}\in {\FSETfunc{\NumBasis}{\NumNest-1}}$ for $\IDEl \in \{1,2,\dots,\NumBasis \}$,
	and 
	$\NLFunc \in {\pBFSETNLFunc{\yBasis}}$ 
	that satisfy \eqref{eq:def_Basis} and \eqref{eq:def_yBasis}.
	\item
	Let ${\pFSETfunc{\NumBasis}{\NumNest}}$ be the following set: 
	\begin{align}
		{\pFSETfunc{\NumBasis}{\NumNest}}
		&:=
		{\pAFSETfunc{\NumBasis}{\NumNest}} \cup {\pBFSETfunc{\NumBasis}{\NumNest}}
		.
	\end{align}

\end{enumerate}
\end{definition}

\begin{theorem}[\MyHighlight{Standard deviation of GPs in ${\pFSETfunc{\NumBasis}{4}}$}] \label{thm:ex_GPsd}	
	Suppose 
	that $\NumBasis \geq \NumData+1$ holds,
	that $\kernel{ \wildcard }{\DState{\IDdat}}$ is contained in ${\FSETfunc{\NumBasis}{1}}$ or is the SE kernel in \eqref{eq:def_SEkernel}, 
	and
	that there exists $\LBGPsd>0$ satisfying ${\El{\GPsd(\State)}{\IDEl}} \geq \LBGPsd $ for all $\State \in \DomBigX$.
	Each component of the GP standard deviation $\GPsd$  in \eqref{eq:def_GPsd} is contained in ${\pFSETfunc{\NumBasis}{4}}$.	
\end{theorem}
\begin{proof}
	The proof is given in Appendix \ref{pf:ex_GPsd}.	
\end{proof}

\begin{proposition}[\MyHighlight{Key properties of ${\pFSETfunc{\NumBasis}{\NumNest}}$}]\label{thm:pFSETfunc_properties}
	For any $\NumNest \in \mathbb{N}$, the following properties hold.
	\begin{enumerate}
		\item\label{thm:Lipschitz_continuous}
		Any $\Basis \in {\pFSETfunc{\NumBasis}{\NumNest}}$ is continuous on $\DomBigX$ and Lipschitz continuous on $\DomX$.
		
		\item\label{thm:pFSETfunc_nondec}
		We have ${\pFSETfunc{\NumBasis}{\NumNest}} \subseteq {\pFSETfunc{\NumBasis}{\NumNest+1}}$ and ${\FSETfunc{\NumBasis}{\NumNest}} \cap \NonNegativeFSETNLFunc  	\subseteq {\pFSETfunc{\NumBasis}{\NumNest+1}}$.
		
	\end{enumerate}	
\end{proposition}
\begin{proof}
	The proof is given in Appendix \ref{pf:pFSETfunc_properties}.		
\end{proof}

\begin{remark}[\MyHighlight{Absolute map}]
	A useful example of $\NLFunc \in {\pBFSETNLFunc{\yBasis}}$ is the absolute mapping $\NLFunc(\yBasis)=|\yBasis|$.	
\end{remark}

We propose Algorithm ${\AlgpNest{\Basis}{\NumNest}{\allDomSimpX}}$ that uses Lemmas \ref{thm:margins_for_inner_products_specific} and \ref{thm:margins_for_nonlinear_mappings_specific} presented below.

\begin{algorithm}[!t]    
	\renewcommand{\algorithmicrequire}{\textbf{Input:}}
	\renewcommand{\algorithmicensure}{\textbf{Output:}}                       
	\caption{${\AlgpNest{\Basis}{\NumNest}{\allDomSimpX}}$ }  
	
	\label{alg:nonnegative_nested_margins}
	\begin{algorithmic}[1]                  
		\REQUIRE $\Basis \in {\pFSETfunc{\NumBasis}{\NumNest}}$, $\NumNest\geq 1$, and $\allDomSimpX$
		\ENSURE ${\zzpEachMBGset{ \Basis }}=[ \zzzUM{\Basis}, \zzzUB{\Basis}, \zzzLipConst{\Basis}   ]^{\MyTRANSPO}$
		
		\STATE
		Decompose $\Basis$ into ${\subAbasis{\IDEl}}$, ${\subBbasis{\IDEl}}$, for $\IDEl \in \{1,2,\dots,\NumBasis \}$, and $\NLFunc$
		according to Definition \ref{def:nonnegative_nested_set}
		
		\FOR{all $\PreliFunc \in \{\subAbasis{1},\dots,\subAbasis{\NumBasis},\subBbasis{1},\dots,\subBbasis{\NumBasis}   \} $  }

		\IF{ $\PreliFunc \in {\FSETfunc{\NumBasis}{1}} $ }
		\STATE
		Obtain ${\zzEachMBGset{ \PreliFunc }}$ by Lemma \ref{thm:margin_FSET1}
		and Definition \ref{def:collection}
		
		\ELSE
		\IF{ $\PreliFunc \in {\FSETfunc{\NumBasis}{\NumNest-1}} $ }
		\STATE
		Calculate ${\zzEachMBGset{ \PreliFunc }}$ via ${\AlgNest{\PreliFunc}{\NumNest-1}{\allDomSimpX}}$
		
		\ELSE			
		\STATE
		Calculate ${\zzpEachMBGset{ \PreliFunc }}$ via ${\AlgpNest{\PreliFunc}{\NumNest-1}{\allDomSimpX}}$
		
		\ENDIF	
		\ENDIF

		\ENDFOR

		\STATE
		Obtain ${\zzzUM{\yBasis}}$ or $({\zzzLM{\yBasis}}, {\zzzUM{\yBasis}})$ by Lemma \ref{thm:margins_for_inner_products} or \ref{thm:margins_for_inner_products_specific}
		\STATE
		Calculate ${\zzzUM{\Basis}}$ by Lemma \ref{thm:margins_for_nonlinear_mappings_specific} using $( {\zzzUM{\yBasis}},{\zzzLipConst{\yBasis}})$ or $( {\zzzUM{\yBasis}},{\zzzLM{\yBasis}})$, where ${\zzzLipConst{\yBasis}}$ is given by Definition \ref{def:collection}

		\STATE
		Obtain ${\zzzUB{\Basis}}$ and  ${\zzzLipConst{\Basis}}$ 
		by Definition \ref{def:collection}
		
		\RETURN

	\end{algorithmic}
\end{algorithm}

\begin{lemma}[\MyHighlight{Margins for specific inner products}] \label{thm:margins_for_inner_products_specific}
	For any $\NumNest\geq 1$,
	any ${\subAbasis{\IDEl}} \in {\pFSETfunc{\NumBasis}{\NumNest}} $,
	and
	any ${\subBbasis{\IDEl}} \in {\pFSETfunc{\NumBasis}{\NumNest}}$  for $\IDEl \in \{1,2,\dots,\NumBasis \}$,  
	an upper margin ${\zzzUM{\yBasis}}$ of $\yBasis(\State)$ in \eqref{eq:def_yBasis} is given by	
	\begin{align}	
		&
		{\zzzUM{\yBasis}}
		=
		\sum_{\IDEl=1}^{\NumBasis}
		\Big(
		{\zzzUM{\subAbasis{\IDEl}}} {\zzzUB{\subBbasis{\IDEl}}} 
		+
		{\zzzUM{\subBbasis{\IDEl}}} {\zzzUB{\subAbasis{\IDEl}}} 
	+
		\NewDist^{2} 	{\zzzLipConst{\subAbasis{\IDEl}}}{\zzzLipConst{\subBbasis{\IDEl}}}
		\Big)
		.\label{eq:chi_bound_U_spe}			
	\end{align}		
Moreover, this ${\zzzUM{\yBasis}}$ is second-order $\MyOrder(\NewDist^{2})$ if $\zzzUM{\subAbasis{\IDEl}}$ and $\zzzUM{\subBbasis{\IDEl}}$ for $\IDEl \in \{1,\dots,\NumBasis\}$ are $\MyOrder(\NewDist^{2})$.
\end{lemma}
\begin{proof}
	The proof is given in Appendix \ref{pf:margins_for_inner_products_specific}.		
\end{proof}	

\begin{lemma}[\MyHighlight{Margins for specific nonlinear mappings}]  \label{thm:margins_for_nonlinear_mappings_specific}
	For any continuous $\yBasis:\DomBigX \to \mathbb{R}$, 
	any $\NLFunc $,
	and any $\Basis$ described by \eqref{eq:def_Basis},
	upper margins of $\Basis$ are given as follows.
\begin{enumerate}
	\item 
	If $\yBasis$ is Lipschitz continuous on $\DomX$ and $\NLFunc \in {\pAFSETNLFunc{\yBasis}}$ holds, 
	then we have
	\begin{align}		
		{\zzzUM{\Basis}}
		&=	 
		\zzzUM{\yBasis}  
		{\dsup{\NLFunc}} 
		+
		\NewDist^{2}
		\frac{ \DimX {\zzzLipConst{\yBasis}}^{2}}{8}
		\max\{0, -{\ddinf{\NLFunc}} \}
		,\label{eq:UM_nondec_NLFunc}			
	\end{align}	
	where ${\zzzUM{\Basis}}$ is second-order $\MyOrder(\NewDist^{2})$ if  ${\zzzUM{\yBasis}}$  is $\MyOrder(\NewDist^{2})$.
	\item		
	If $\yBasis$ is $C^{2}$ continuous on $\DomBigX$ and $\NLFunc \in {\pBFSETNLFunc{\yBasis}}$ holds, 
	then we have
	\begin{align}		
		{\zzzUM{\Basis}}
		&=	 
		\LipsNLFunc
		\max \{ 
		\zzzLM{\yBasis}
		,    
		\zzzUM{\yBasis}  
		\}
		,\label{eq:UM_convex_NLFunc}					
	\end{align}	
	where ${\zzzUM{\Basis}}$ is $\MyOrder(\NewDist^{2})$ if ${\zzzLM{\yBasis}}$ and ${\zzzUM{\yBasis}}$  are $\MyOrder(\NewDist^{2})$.			
\end{enumerate}	
\end{lemma}
\begin{proof}
	The proof is given in Appendix \ref{pf:margins_for_nonlinear_mappings_specific}.
\end{proof}

The theoretical aspects of Theorem \ref{thm:main_results_margins_specific} has been established.
The proof of Theorem \ref{thm:main_results_margins_specific} is similar to that of Theorem \ref{thm:main_results_margins}.

\section{Solutions to Problems 1 and 2}\label{sec_solutions}

In this section, we propose a method to solve Problems 1 and 2, based on the results in Section \ref{sec_2ndOrderMargin}.
An overview of the proposed method is described in Section \ref{sec_overview_solving_problems123}.
Sections \ref{sec_solution1} and \ref{sec_solution2} provide solutions to Problems 1 and 2, respectively.

\subsection{Overview}\label{sec_overview_solving_problems123}

We find a stability region $\MDomStable$ via evaluating a finite number of inequalities as follows.

\begin{theorem}[\MyHighlight{Finding a stability region}] \label{thm:find_stability_region}
	Suppose that
	 a $C^{1}$ continuous Lyapunov function $\MValF$, 
	 function $\HdotV$ in Definition \ref{def:stability_region}, 
	 collection $\allDomSimpX$ of simplexes in Definition \ref{def:simplexes},
	 and 
	 margins ${\zzzLM{\MValF}}$ and ${\zzzUM{\HdotV}}$ in Definition \ref{def:margins}
	 are given.
	The following $\MDomStable$ is a stability region:
	\begin{align}		
	\MDomStable 
	&= 
	\bigcup_{\IDdom \in \SetStableIDdom} \DomSimpX{\IDdom}
	, \label{eq:def_MDomStable}
	\end{align}	
	where the set $\SetStableIDdom$ is defined as follows:
	\begin{align}		
	\SetStableIDdom
	&:=
	\Bigg\{
	\IDdom \in \SetALLIDdom 
	\Bigg|	 
	\begin{matrix}
	\min_{\IDPt \in \{1,\dots,\NumSimp\}  } \MValF(\SimpPt{\IDdom}{\IDPt}) - {\zzzLM{\MValF}}
	> 0
	,\\
	\max_{\IDPt \in \{1,\dots,\NumSimp\} }   \HdotV(\SimpPt{\IDdom}{\IDPt})	+ {\zzzUM{\HdotV}}
	< 0
	\end{matrix}	
	\Bigg\}
	.	\label{eq:cal_MDomStable}
	\end{align}	
\end{theorem}
\begin{proof}
	By Proposition \ref{def:bounds},
	for any $ \IDdom \in \SetStableIDdom$ and any $\State \in \DomSimpX{\IDdom}$,
	we have 
	$\MValF(\State) 
	\geq \min_{\IDPt
	  } \MValF(\SimpPt{\IDdom}{\IDPt}) - {\zzzLM{\MValF}}
	>0 $ and 
	$\HdotV(\State) 
	\leq \max_{\IDPt }   \HdotV(\SimpPt{\IDdom}{\IDPt})	+ {\zzzUM{\HdotV}}
	< 0$.
	Thus, $\DomSimpX{\IDdom}$ for $\IDdom \in \SetStableIDdom$ satisfies \eqref{eq:stab_conditions} and thus is contained in $\MDomStable$.
	This completes the proof.	
\end{proof}	
\begin{remark}[\MyHighlight{Contribution of Theorem \ref{thm:find_stability_region}}]
	Theorem \ref{thm:find_stability_region} finds a stability region $\MDomStable$ by evaluating only a finite number  of  $ \MValF(\SimpPt{\IDdom}{\IDPt}) $ and $ \HdotV(\SimpPt{\IDdom}{\IDPt}) $ for $\IDdom \in \SetALLIDdom$ and $\IDPt \in \{1,\dots,\NumSimp\}$.
\end{remark}

\begin{remark}[\MyHighlight{Challenge in solving Problems 1 and 2}]
	A crucial challenge in employing Theorem \ref{thm:find_stability_region} is to find margins ${\zzzLM{\MValF}}$ and ${\zzzUM{\HdotV}}$ as sufficiently small values.
	Such small margins are desirable for evaluating $\MDomStable$ accurately.  
	Our proposed method with Theorems \ref{thm:main_results_margins} and \ref{thm:main_results_margins_specific} in Section \ref{sec_2ndOrderMargin} is a promising tool to overcome this challenge.	
	Applying Theorems \ref{thm:main_results_margins} and \ref{thm:main_results_margins_specific} to control problems derives ${\zzzLM{\MValF}}$ and ${\zzzUM{\HdotV}}$ as second-order $\MyOrder(\NewDist^{2})$ in the subsequent subsections.
\end{remark}

\subsection{Solution to Problem 1: Finding a stability region for given Lyapunov functions and controllers}\label{sec_solution1}

We solve Problem 1 to find a stability region $\MDomStable$ if a Lyapunov function $\MValF(\State)$ and controller $\Input(\State)$ are given.
For various systems, controllers, and Lyapunov functions, which can be data-driven,
second-order margins ${\zzzLM{\MValF}}$ and ${\zzzUM{\HdotV}}$ are derived based on Theorems \ref{thm:main_results_margins} and \ref{thm:main_results_margins_specific}.
Combining the derived margins with Theorem \ref{thm:find_stability_region} provides a solution to Problem 1.

We introduce the following conditions to determine the classes of functions.

\begin{assumption}[\MyHighlight{Classes of functions}]\label{ass:FuncClass}
	There exist {known} positive integers 
	${\fNumNest{\Mmean}}$,	
	${\fNumNest{\Msd}}$,
	${\fNumNest{\MValF}}$, and
	${\fNumNest{\partial\MValF}}$ that satisfy the following properties  for every $\IDEl \in \{1,\dots,\DimX\}$.
	\begin{enumerate}
		\item
		For nominal dynamics, ${\El{\Mmean(\wildcard,\Input(\wildcard))}{\IDEl}} \in {\FSETfunc{\NumBasis}{\fNumNest{\Mmean}}}$ holds.	
		
		\item
		For uncertainty, ${\El{\Msd(\wildcard,\Input(\wildcard))}{\IDEl}} \in {\FSETfunc{\NumBasis}{\fNumNest{\Msd}}}
		\cup {\pFSETfunc{\NumBasis}{\fNumNest{\Msd}+1}}$ holds.

		\item
		For a Lyapunov function, $\MValF \in {\FSETfunc{\NumBasis}{\fNumNest{\MValF}}}$ and ${\El{\partial\MValF/\partial\State}{\IDEl}}  \in {\FSETfunc{\NumBasis}{\fNumNest{\partial\MValF}}}$ hold.

		\item
		Given $\NumBasis$ in the top of Section \ref{sec_set_of_func}, $\NumBasis \geq \DimX$ holds.
		
	\end{enumerate}
\end{assumption}

\begin{theorem}[\MyHighlight{Lyapunov inequalities}]\label{thm:HdotV_margin}
Given a Lyapunov function $\MValF: \mathbb{R}^{\DimX} \to \mathbb{R} $, controller $\Input$, nominal dynamics $\Mmean$, and uncertainty $\Msd$, suppose Assumption \ref{ass:FuncClass}.
The functions $\HdotV$, $\MEdotV$, and $\SDdotV$ in \eqref{eq:def_HdotV}--\eqref{eq:def_SDdotV} satisfy the following relations: 
	\begin{align}
	{\zzzUM{\HdotV}} & = {\zzzUM{\MEdotV}} + {\zzzUM{\SDdotV}}
	,\label{eq:sum_margins_HdotV}
	\\
	\MEdotV &\in {\FSETfunc{\NumBasis}{\fNumNest{\MEdotV}}}
	,\\
	\SDdotV &\in {\pFSETfunc{\NumBasis}{\fNumNest{\SDdotV}}}
	,
	\\
	{\fNumNest{\MEdotV}} 
	&:= \max\{ {\fNumNest{\partial\MValF}} ,  \fNumNest{\Mmean} \} + 1
	,\label{eq:def_fNumNest_MEdotV}
	\\
	{\fNumNest{\SDdotV}} 
	&:= \max\{ {\fNumNest{\partial\MValF}} ,  \fNumNest{\Msd} \} + 2
	.\label{eq:def_fNumNest_SDdotV}
	\end{align}	
\end{theorem}
\begin{proof}
The proof is given in Appendix \ref{pf:HdotV_margin}.		
\end{proof}	

\begin{remark}[\MyHighlight{Contribution of Theorem \ref{thm:HdotV_margin}}]
	By virtue of Theorem \ref{thm:HdotV_margin}, a second-order upper margin ${\zzzUM{\HdotV}}$ is obtained by using Algorithms ${\AlgNest{\MEdotV}{\fNumNest{\MEdotV}}{\allDomSimpX}}$ and ${\AlgpNest{\SDdotV}{\fNumNest{\SDdotV}}{\allDomSimpX}}$
	owing to Theorems \ref{thm:main_results_margins} and \ref{thm:main_results_margins_specific}.
	Algorithm ${\AlgNest{\MValF}{\fNumNest{\MValF}}{\allDomSimpX}}$ 
	derives a second-order lower margin ${\zzzLM{\MValF}}$ under Assumption \ref{ass:FuncClass}.
	After calculating these margins, we are ready to employ Theorem \ref{thm:find_stability_region}.
\end{remark}

We have solved Problem 1 via Algorithm {\MyLabelAlgSR} using Theorems \ref{thm:main_results_margins}, \ref{thm:main_results_margins_specific}, \ref{thm:find_stability_region}, and \ref{thm:HdotV_margin} as follows.

\begin{corollary}[\MyHighlight{Solution to Problem 1}]\label{thm:solution_to_Problem1}
	Given a Lyapunov function $\MValF: \mathbb{R}^{\DimX} \to \mathbb{R} $, controller $\Input$, nominal dynamics $\Mmean$, and uncertainty $\Msd$, suppose Assumption \ref{ass:FuncClass}.
	For a collection $\allDomSimpX$ of simplexes,
	Algorithm {\MyLabelAlgSR} presents a stability region $\MDomStable$ given by \eqref{eq:def_MDomStable}.
\end{corollary}

\begin{algorithm}[!t]    
	\renewcommand{\algorithmicrequire}{\textbf{Input:}}
	\renewcommand{\algorithmicensure}{\textbf{Output:}}                       
	\caption{\textbf{\MyLabelAlgSR}: Finding a stability region $\MDomStable$}  
	\begin{algorithmic}[1]                  
		\REQUIRE $\allDomSimpX$, $\MValF$, $\Input$,  $\Mmean$,  $\Msd$,
		${\fNumNest{\MValF}}$, ${\fNumNest{\partial\MValF}}$, $\fNumNest{\Mmean}$, $\fNumNest{\Msd}$
		
		\ENSURE $\MDomStable$ 
		
		\STATE
		Calculate ${{\MEdotV}}$, ${{\SDdotV}}$, ${\fNumNest{\MEdotV}}$, and ${\fNumNest{\SDdotV}}$ defined in \eqref{eq:def_MEdotV}, \eqref{eq:def_SDdotV}, \eqref{eq:def_fNumNest_MEdotV}, and \eqref{eq:def_fNumNest_SDdotV}, respectively

		\STATE
		Obtain ${\zzzLM{\MValF}}$ by Algorithm  ${\AlgNest{\MValF}{\fNumNest{\MValF}}{\allDomSimpX}}$
		\STATE
		Obtain ${\zzzUM{\MEdotV}}$ by Algorithm ${\AlgNest{\MEdotV}{\fNumNest{\MEdotV}}{\allDomSimpX}}$
		\STATE
		Obtain ${\zzzUM{\SDdotV}}$ by Algorithm ${\AlgpNest{\SDdotV}{\fNumNest{\SDdotV}}{\allDomSimpX}}$
		\STATE
		Obtain ${\zzzUM{\HdotV}}={\zzzUM{\MEdotV}}+{\zzzUM{\SDdotV}}$

		\STATE
		Obtain $\MDomStable$ in \eqref{eq:def_MDomStable} by calculating  $\SetStableIDdom$ in \eqref{eq:cal_MDomStable}
		
	\end{algorithmic}
\end{algorithm}

\subsection{Solution to Problem 2: Finding a stability region with the design of Lyapunov functions and controllers}\label{sec_solution2}

In this subsection, we solve Problem 2 to design a Lyapunov function $\MValF$ and feedback controller $\Input$ so that a stability region $\MDomStable$ contains a given candidate region $\estMDomStable$.
Our approach is to connect finding $\MDomStable$ with an optimization problem. 
We design $\MValF$ and $\Input$ appropriately via solving the optimization problem, thereby obtaining $\MDomStable$ via Algorithm {\MyLabelAlgSR}.

We propose the optimization problem in the following.
Let $\MValF(\State; \VCparam)$ and $\Input(\State; \VCparam)$ be parametric functions of a free parameter $\VCparam \in \SetVCparam$, where $\SetVCparam \subset \mathbb{R}^{\DimVCparam}$ is a bounded closed set and $(\wildcard;\VCparam)$ denotes the dependence of $\VCparam$.
Assume that $\MValF(\State; \VCparam)$ and $\HdotV(\State; \VCparam)$ are continuous in $\VCparam$ on $\SetVCparam$ for each $\State$.

Given $\estMDomStable$ and $\allDomSimpX$, let us consider the following minimization problem:
\begin{align}
\min_{\VCparam \in \SetVCparam} 
&\ObjectiveF(\VCparam;\HJBweight, \estMDomStable, \allDomSimpX)
,\label{eq:def_minimization}
\\
\ObjectiveF(\VCparam;\HJBweight, \estMDomStable, \allDomSimpX)
&:=
\HJBweight
\PerformanceTerm(\VCparam; \estMDomStable, \allDomSimpX)
\nonumber\\&\quad
+
 \sum_{ \State \in  \VerticesestMDomStable } 
{\PenaltyperX{\VCparam}{\State}{\VUBofMargin}{\WUBofMargin}{\OffSetpenaltyFunc}}
,
\\
{\PenaltyperX{\VCparam}{\State}{\VUBofMargin}{\WUBofMargin}{\OffSetpenaltyFunc}}
&:=
\penaltyFunc\big( \VUBofMargin(\State) -\MValF(\State;\VCparam);\OffSetpenaltyFunc \big)
\nonumber\\&\quad
+ \penaltyFunc\big( \WUBofMargin(\State) + \HdotV(\State;\VCparam);\OffSetpenaltyFunc \big)
,\label{eq:def_PenaltyperX}
\end{align}
where
$\PerformanceTerm(\wildcard; \estMDomStable,\allDomSimpX): \SetVCparam \to \NNReal$ is assumed to be continuous and is a control performance index to be decreased, and  $\HJBweight \in \NNReal$ is a free parameter.
Let $\sum_{ \State \in  \VerticesestMDomStable }$ denote the sum with respect to all the members of $\VerticesestMDomStable:=\{ {\SimpPt{\IDdom}{\IDPt}} |  {\SimpPt{\IDdom}{\IDPt}} \in \estMDomStable \} $, where ${\SimpPt{\IDdom}{\IDPt}} $ is determined from $\allDomSimpX$.
The functions $\VUBofMargin: \DomX \to (0, \infty)$ and $\WUBofMargin: \DomX \to (0, \infty)$ are continuous functions.
For a free parameter $\OffSetpenaltyFunc \in \NNReal$, 
let the penalty function $\penaltyFunc(\wildcard;\OffSetpenaltyFunc): \mathbb{R} \to  \NNReal$ be a nondecreasing continuous function to satisfy $\penaltyFunc(\ArgpenaltyFunc;\OffSetpenaltyFunc)\geq  \ArgpenaltyFunc + \OffSetpenaltyFunc$ for $\ArgpenaltyFunc >0$ and $0 \leq \penaltyFunc(\ArgpenaltyFunc;\OffSetpenaltyFunc) \leq \OffSetpenaltyFunc$ for  $\ArgpenaltyFunc \leq 0$.
	Typical examples of $\penaltyFunc$ contain the following $C^{1}$ continuous function with $\OffSetpenaltyFunc>0$:
	\begin{align}
	\penaltyFunc(\ArgpenaltyFunc; \OffSetpenaltyFunc)
	&=
	\begin{cases}
	\ArgpenaltyFunc + \OffSetpenaltyFunc & (\ArgpenaltyFunc> 0)
	\\
	( \ArgpenaltyFunc + 2\OffSetpenaltyFunc )^{2}/(4\OffSetpenaltyFunc)  & (- 2\OffSetpenaltyFunc  < \ArgpenaltyFunc \leq  0)
	\\
	0 & (\ArgpenaltyFunc \leq -2\OffSetpenaltyFunc )
	\end{cases}
	,
	\end{align}
	and the rectified linear function with $\OffSetpenaltyFunc=0$,
	that is, $\penaltyFunc(\ArgpenaltyFunc;0)=\ArgpenaltyFunc$ for $\ArgpenaltyFunc> 0$ and otherwise $\penaltyFunc(\ArgpenaltyFunc;0)=0$.
Before stating our results, we introduce the following assumptions.
\begin{assumption}[\MyHighlight{Feasibility for controller design}]\label{ass:controller_design}
Given $\estMDomStable$ and continuous functions $\VUBofMargin: \DomX \to (0, \infty)$ and $\WUBofMargin: \DomX \to (0, \infty)$, 
the following properties hold.

	\begin{enumerate}

\item[\ASSinNestdSet] 
For every $\VCparam \in \SetVCparam$, Assumption \ref{ass:FuncClass} holds.

		\item[\ASSVCparamExist]
		There exists $\VCparam \in \SetVCparam$ satisfying the condition:
		\begin{align}
		&
		\forall \State \in \estMDomStable
		,\;
		\MValF(\State; \VCparam) \geq \VUBofMargin(\State)   
		,\;
		\HdotV(\State; \VCparam) \leq -\WUBofMargin(\State) 
	.\label{eq:estSR_is_a_SR_offsetFuncs}
		\end{align}

		\item[\ASSboundedMargin]
		There exist $\ODEFUBDist>0$, ${\LUMconst{\MValF}}>0$, and ${\LUMconst{\HdotV}}>0$ such that 
		for every $\VCparam \in \SetVCparam$ and every $\allDomSimpX$ with the corresponding $\NewDist$, 
		second-order margins 
		$\zzzLM{\MValF;\VCparam}$ 
		and
		$\zzzUM{\HdotV;\VCparam}$ 
		given by Algorithm {\MyLabelAlgSR} (under \ASSinNestdSet)
		satisfy
		\begin{align}		
		\NewDist \leq \ODEFUBDist
		\Rightarrow
		\zzzLM{\MValF;\VCparam}
		\leq {\LUMconst{\MValF}} \NewDist^{2}
		,
		\zzzUM{\HdotV;\VCparam}
		\leq {\LUMconst{\HdotV}} \NewDist^{2}
		.
		\end{align}

		\item[\ASSnonemptySetDist]
		The set $\SetDist$ is nonempty, where $\SetDist \subset (0,\infty)$ is the set of possible values of $\NewDist>0$ such that
		there exist a collection $\allDomSimpX$ of simplexes and nonempty set $\estSetStableIDdom \subseteq \SetALLIDdom$ that satisfy the definition \eqref{eq:def_tau} 
		and $\estMDomStable 
		= 
		\bigcup_{\IDdom \in \estSetStableIDdom } 
		\DomSimpX{\IDdom}$.

	\end{enumerate}
\end{assumption}

We show that solving the minimization \eqref{eq:def_minimization} gives a Lyapunov function $\MValF$, controller $\Input$, and stability region $\MDomStable$ simultaneously.
Here, for $(\HJBweight,\estMDomStable,\allDomSimpX,\OffSetpenaltyFunc)$, let $\optVCparam$ be a minimizer to \eqref{eq:def_minimization}, denoting $\stabVCparam$ if $\HJBweight=0$.

\begin{theorem}[\MyHighlight{Controller design}] \label{thm:general_control_design}
	
There exists $\UBDist>0$ such that for every $\NewDist \in (0, \UBDist] \cap \SetDist $ and every $\allDomSimpX$ satisfying \eqref{eq:def_tau},  the following (i) holds under Assumption \ref{ass:controller_design}. 	
Moreover, 
the following (ii) holds for every $\allDomSimpX$ without Assumption \ref{ass:controller_design}.
	\begin{enumerate}

		\item			
		There exist $\UBHJBweight>0$ and $\UBOffSetpenaltyFunc>0$ 
		such that
		for any $\HJBweight \in [0, \UBHJBweight]$ and $\OffSetpenaltyFunc \in [0, \UBOffSetpenaltyFunc]$,		
		Algorithm {\MyLabelAlgSR} with $\MValF(\wildcard; \optVCparam) $, $\MEdotV(\wildcard; \optVCparam)$, and $\SDdotV(\wildcard; \optVCparam)$
		gives a stability region $\MDomStable$ satisfying $\estMDomStable \subseteq \MDomStable $.
		
		\item
		The following improvement with respect to the performance index $\PerformanceTerm$ holds:
		\begin{align}
		&
		\PerformanceTerm(\optVCparam; \estMDomStable, \allDomSimpX)
		\leq 
		\PerformanceTerm(\stabVCparam; \estMDomStable, \allDomSimpX)
		.\label{eq:Design_suboptimality}
		\end{align}

	\end{enumerate}

\end{theorem}
\begin{proof}
	The proof is given in Appendix \ref{pf:general_control_design}.	
\end{proof}

\begin{remark}[\MyHighlight{Solution to Problem 2 via Theorem \ref{thm:general_control_design}}]
By virtue of Theorem \ref{thm:general_control_design} (i), we solve Problem 2 as follows.
Solving the minimization \eqref{eq:def_minimization} can give an appropriate Lyapunov function $\MValF(\wildcard; \VCparam)$ and feedback controller $\Input(\wildcard; \VCparam)$ for small $\NewDist$.
Thereafter, Algorithm {\MyLabelAlgSR} gives a stability region $\MDomStable$ containing the candidate region $\estMDomStable$.
Section \ref{sec_demonstration} demonstrates that $\MValF(\wildcard; \VCparam)$, $\Input(\wildcard; \VCparam)$, and $\MDomStable$ can be successfully obtained even though the minimization \eqref{eq:def_minimization} may not be solved in a exact sense.  
Theorem \ref{thm:general_control_design} (ii) indicates the performance improvement for $\HJBweight>0$ in comparison with the case of $\HJBweight=0$ ignoring the performance. 
\end{remark}

\begin{remark}[\MyHighlight{Other contributions of Theorem \ref{thm:general_control_design}}]
	The proposed minimization \eqref{eq:def_minimization} overcomes several difficulties in optimization-based controller design as follows.
	If we have the minimum ${\PenaltyperX{\VCparam}{\State}{\VUBofMargin}{\WUBofMargin}{\OffSetpenaltyFunc}}|_{\OffSetpenaltyFunc=0}=0$ for every $\State  \in  \estMDomStable$, 
	that is, $\MValF(\State; \VCparam)  >  0$ and $\HdotV(\State; \VCparam)  < 0$,
	then $\estMDomStable$ is a stability region. 
	This fact motivates us to consider the minimization of
	$
	\int_{ \State  \in  \estMDomStable }
	{\PenaltyperX{\VCparam}{\State}{\VUBofMargin}{\WUBofMargin}{0}}
	\mathrm{d}\State
	$
	straightforwardly.
	However, this minimization suffers from three drawbacks.
	Firstly, it is difficult to evaluate and minimize the integral value
	$\int_{ \State  \in  \estMDomStable }
	(\dots)
	\mathrm{d}\State$. 
	Secondly, any control performance is not treated while a stability region may be found.
	Finally, because of the condition $\OffSetpenaltyFunc=0$, $\penaltyFunc$ is not differentiable or its gradient approaches to zero as ${\PenaltyperX{\VCparam}{\State}{\VUBofMargin}{\WUBofMargin}{\OffSetpenaltyFunc}} \to 0$.
	Theorem \ref{thm:general_control_design} shows that the minimization \eqref{eq:def_minimization} overcomes these drawbacks.
\end{remark}

\section{Demonstration with a numerical simulation}\label{sec_demonstration}

This section demonstrates designing a suboptimal controller with guaranteeing the stability of a data-driven system model.

\subsection{Plant system and setting \label{sec_plant} }

Consider the partially-unknown pendulum system:
\begin{align}	
	\dotState
	&= \Drift(\State,  \Input(\State) ) 
	=\onlyDrift(\State)
	+\InMat  \Input(\State)
	,\label{eq:sim_sys}
	\\
	\onlyDrift(\State)
	&:=
	\begin{bmatrix}
		{\El{\State}{2}} \\
		- 9.8 \sin {\El{\State}{1}}  -{\El{\State}{2}}
	\end{bmatrix}
	,
	\quad
	\InMat
	:=
	\begin{bmatrix}
		0 \\ 1
	\end{bmatrix}
	.
\end{align}
Supposing that $\InMat$ is known but $\onlyDrift$ is unknown, we develop the following nominal model $\Mmean$ of this pendulum by using $\GPmean(\State)$ in \eqref{eq:def_GPmean}, which is equivalent to kernel-ridge regression \cite{Vovk2013}:
\begin{align}	
	\Mmean(\State,\Input)
	&=
	[\DDrift{1}, \DDrift{2}, \dots, \DDrift{\NumData}]
	\kernelMat^{-1}\kernelVec(\State)+ \InMat \Input
	,
\end{align}
where $\NumData=121$.
The data points $\DState{\IDdat}$ are generated on $\DomX=[-8,8] \times [-8,8]$ at regular intervals.
Each output $\DDrift{\IDdat}$ is assumed to be obey the normal distribution with the mean $\onlyDrift(\DState{\IDdat})$  and covariance $ 0.05^{2} {\Identity{2}}$ that is independently distributed with respect to $\IDdat$.
We employ the SE kernel \eqref{eq:def_SEkernel} with the settings of $\hypMag = 1$, $\hypCovMat = 5 {\Identity{2}}$, and $\hypNoise=0.001$.

In the control design in Section \ref{sec_solution2}, a suboptimal controller is designed,
that is, the performance index $\PerformanceTerm$ indicates the residuals of the Hamilton-Jacobi-Bellman equation  $\HJB(\State;\VCparam)=0$ \cite{Lewis86} as follows:
\begin{align}
	\PerformanceTerm(\VCparam;  \estMDomStable, \allDomSimpX)
	&=
	\frac{1}{\NumVerticesestMDomStable}
	\sum_{ \State \in  \VerticesestMDomStable }
	\HJB(\State;\VCparam)^{2}
	,\\
	\HJB(\State;\VCparam)
	&:=
	\frac{\partial \MValF(\State; \VCparam)}{\partial\State^{\MyTRANSPO}} \GPmean(\State)
	-
	\frac{1}{2}
		\Big(\InMat^{\MyTRANSPO} \frac{\partial \MValF(\State; \VCparam)}{\partial\State} 
		\Big)^{2}
	\nonumber\\&\quad
	+
	5({\El{\State}{1}})^{2}
	+
	0.01({\El{\State}{2}})^{2}
	,\\
	\Input(\State; \VCparam)&=- \InMat^{\MyTRANSPO}\partial \MValF(\State; \VCparam)/\partial\State
	,\\
	\MValF(\State; \VCparam)&
	=\VCparam^{\MyTRANSPO} \kernelMat^{-1}  ( \kernelVec(\State)- \kernelVec(0) )
	,
\end{align}
where $\NumVerticesestMDomStable$ denotes the number of all the members of $\VerticesestMDomStable$.
As illustrated in Fig. \ref{fig:simplexes}, $\DomX$ is decomposed into 80,000 simplexes $\allDomSimpX$ with $\NewDist=0.08\sqrt{2}\approx 0.113$ at regular intervals, that is,
$\NewDist=\|\SimpPt{\IDdom}{2}-\SimpPt{\IDdom}{1}\|
=\|\SimpPt{\IDdom}{3}-\SimpPt{\IDdom}{1}\|
=\|\SimpPt{\IDdom}{3}-\SimpPt{\IDdom}{2}\|/\sqrt{2}$ holds.
The candidate set is defined by 
$
\estMDomStable =   \{  \State \in \DomX  |  \|\State\| >  0.1  \}$.

The minimization \eqref{eq:def_minimization} is solved by using a gradient method with a line search algorithm.
The objective function is determined from the settings of 
$
\VUBofMargin(\State)
:=
\WUBofMargin(\State)
:=
\|\State\|^{2} +0.1
$,
$\HJBweight=1$,
and $\OffSetpenaltyFunc=1.0\times 10^{-20}$.
The initial values of $\VCparam$ are determined by the following LQR settings.
We obtain a quadratic Lyapunov function $\quadMValF(\State)$ by solving $\HJB(\State)=0$ with replacing 
$\MValF(\State)$ and $\GPmean(\State)$ with $\quadMValF(\State)$ and $({\partial \GPmean(0)}/{\partial \State^{\MyTRANSPO} })  \State$, respectively, which yields the linear quadratic controller.
The initial value is set to ${\El{\VCparam}{\IDdat}}:=\quadMValF(\DState{\IDdat}) - \max_{\IDdat} \quadMValF(\DState{\IDdat})$.

In performing Algorithm {\MyLabelAlgSR}, we implement an efficient decomposition of $\DomX$ into simplexes.
Firstly, we perform Algorithm {\MyLabelAlgSR} using the 80,000 simplexes $\allDomSimpX$ as defined above.  
If some simplex $ \DomSimpX{\IDdom}$ is not a stability region, we perform Algorithm {\MyLabelAlgSR} with replacing $\DomX$ with $\DomSimpX{\IDdom}$.
Namely, $ \DomSimpX{\IDdom}$ is decomposed into 80,000 simplexes with $\NewDist\approx 5.7\times 10^{-4}$.
By iterating this decomposition and evaluations three times, we obtain fine results with small $\NewDist\approx 2.8 \times 10^{-6} $ with low computational burden.

\subsection{Results \label{sec_results_stability_region}}

\newcommand{\XFigODE}{-0.1}
\newcommand{\WFigODE}{3.6in}
\newcommand{\HFigODE}{3.6in}
\newcommand{\SpaceFigODE}{}

\newcommand{\WFigDist}{2.65in}
\newcommand{\HFigDist}{1.95in}
\newcommand{\SpaceFigDist}{}

We evaluate the stability region $\MDomStable$ obtained by the proposed method,
comparing the designed controller and Lyapunov function with their initial LQR versions.
In Fig. \ref{fig:stability_region}, \subref{fig:initial_sim} and \subref{fig:designed_sim} show the results for the initial and designed versions, respectively.
We see that employing the proposed method finds the stability regions $\MDomStable$ (white regions) successfully in both the results  in \subref{fig:initial_sim} and \subref{fig:designed_sim}.
Moreover, comparing both the results shows that the stability region $\MDomStable$ is enlarged by the proposed controller design 
because most area in \subref{fig:designed_sim} was $\MDomStable$ except for around the origin.
Indeed, the states were stabilized for all the initial states  in \subref{fig:designed_sim}.

\newcommand{\MyFigODE}[1]{
	\begin{tikzpicture}[scale=1.0]
	\begin{axis}[axis y line=none, axis x line=none
	,xtick=\empty, ytick=\empty 
	,xmin=-10,xmax=10,ymin=-10,ymax=10
	,width=\WFigODE,height=\HFigODE
	]
	
	
	\node at (\XFigODE,0.9) { 	\includegraphics[width=0.98\linewidth]{#1}  };
	\node[rotate=90] at (-9.0,1.5) {{\small \colorbox{white}{$\qquad {\El{\State}{2}} \qquad$}}};
	\fill[fill=white] (-2,-8.5) rectangle (3,-11);
	\node at (1.5,-9.3) {{\small $\qquad  {\El{\State}{1}} \qquad$}};				
	\end{axis}			
	\end{tikzpicture}	
	\SpaceFigODE
}

\begin{figure}[!t]		
	\begin{minipage}[b]{.5\linewidth}
		\centering		
		{\MyFigODE{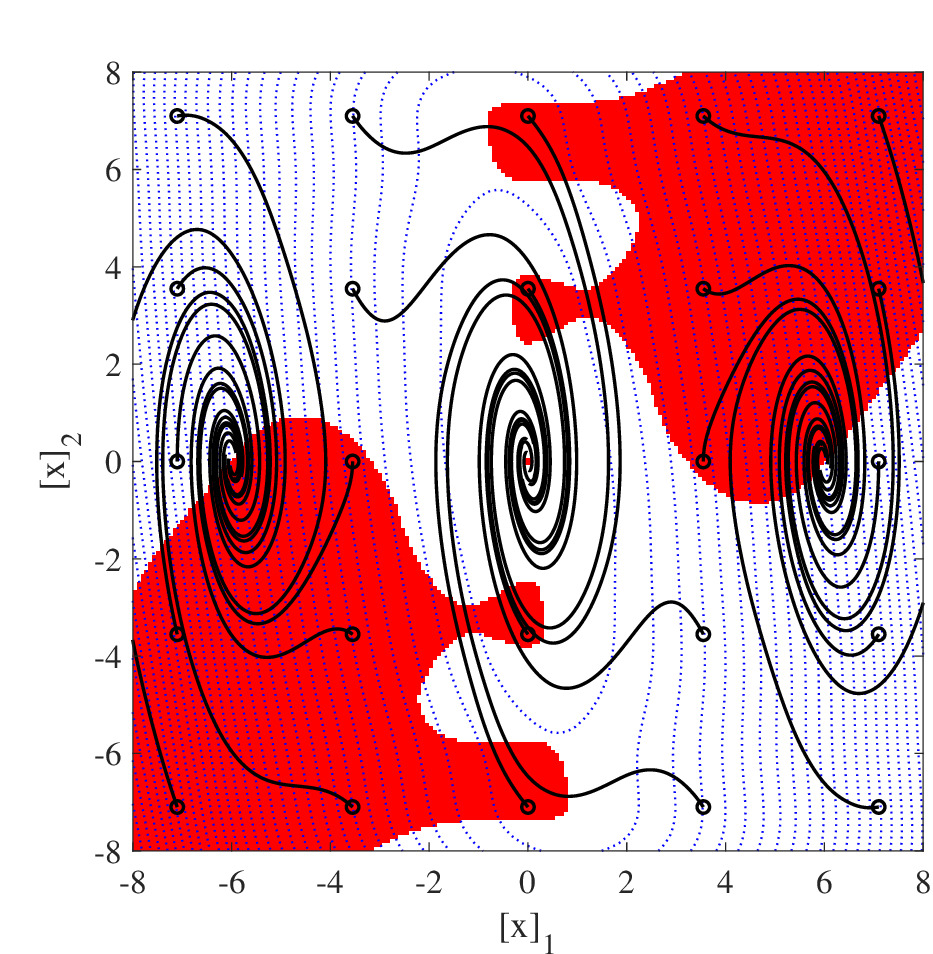}}
		\subcaption{With initial $\Input$ and $\MValF$}
		\label{fig:initial_sim}							
	\end{minipage}%
	\begin{minipage}[b]{.5\linewidth}
		\centering	
		{\MyFigODE{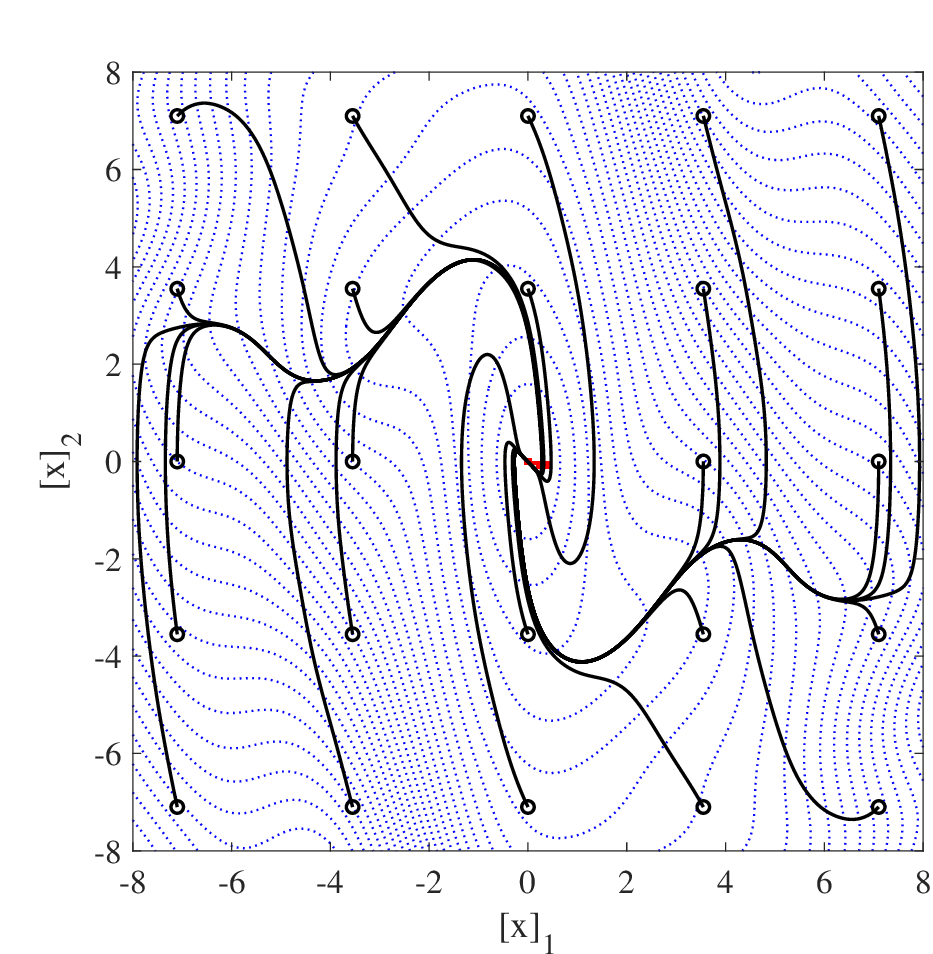}}	
		\subcaption{With designed $\Input$ and $\MValF$}
		\label{fig:designed_sim}							
	\end{minipage}%
	
	\caption{
		Stability evaluations for the initial and designed versions of controllers and Lyapunov functions.
		The white and red regions denote the obtained stability regions $\MDomStable$ and their complements, respectively, where the spatial resolution is $0.08$.
		The blued dashed and black lines indicate the contour lines of the Lyapunov function and the state trajectories, respectively.}
	
	\label{fig:stability_region}
\end{figure}

Next, we evaluate the proposed second-order margins.
Figure \ref{fig:@order_analysis} with the double logarithmic axes represents the margins for various values of the interval $\NewDist$.
The margins were successfully decreased along with the second-order red line and became significantly small compared with the first-order blue baseline.
This indicates that using the proposed margins gives a stability region $\MDomStable$ with high accuracy for small $\NewDist$.

\newcommand{\MyFigDist}[1]{
	\begin{tikzpicture}[scale=1.0]
	\begin{axis}[axis y line=none, axis x line=none
	,xtick=\empty, ytick=\empty 
	,xmin=-10,xmax=10,ymin=-10,ymax=10
	,width=\WFigDist,height=\HFigDist
	]
	
	
	\node at (0,-0.3) { 	\includegraphics[width=0.59\linewidth]{#1}  };
	\node[rotate=90] at (-8.5,2) {{\small \colorbox{white}{$\qquad \zzzLM{\MValF}  \qquad$}}};
	\fill[fill=white] (-2,-8) rectangle (3,-11);
	\node at (1.5,-9.2) {{\small $\qquad \NewDist \qquad$}};				
	\end{axis}			
	\end{tikzpicture}	
	\SpaceFigDist		
}

\begin{figure}[!t]		
	\begin{minipage}{.5\linewidth}
		\centering				
		{\MyFigDist{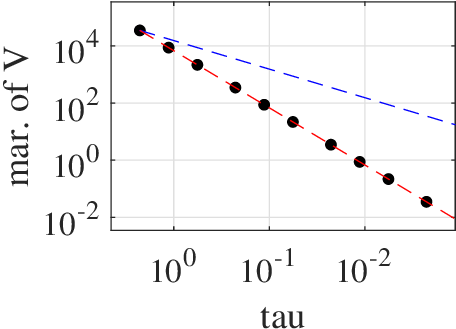}}
		\subcaption{$\zzzLM{\MValF}$}
	\end{minipage}%
	\begin{minipage}{.5\linewidth}
		\centering	
		{\MyFigDist{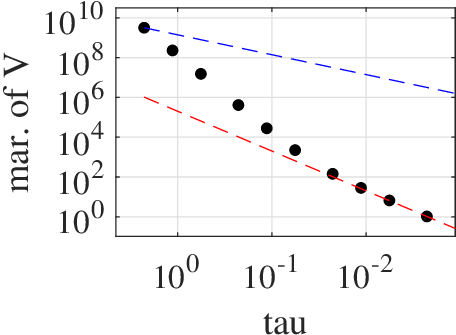}}			
		\subcaption{$\zzzUM{\HdotV}$}
	\end{minipage}%
	\caption{
		The obtained margins $\zzzLM{\MValF}$ and  $\zzzUM{\HdotV}$ for the designed $\MValF$ denoted by the black dots. 
		The upper blue and lower red dashed lines indicate (baseline) first-order and second-order decreasing properties, respectively.
		The margins were calculated on $[-0.8 \times 0.8]\times [-0.8 \times 0.8]$.
	}
	\label{fig:@order_analysis}
\end{figure}

\section{Conclusion}\label{sec_conclusion}

This study presented a control method to guarantee stability via a sampling-based manner. 
The method is applied to general (data-driven) control systems with uncertainty, such as GPs, DNNs, kernel-based models, polynomials, transcendental functions, and their compositions described in Table \ref{tab:results}.
Second-order margins for such systems are derived to realize the precise sampling-based stability analysis.  
Controllers are designed to guarantee the stability, considering control performance.

The proposed general theory can contribute to various types of systems and control methods.
The proposed system classes have the potential to include even more data-driven models beyond the aforementioned examples.
Because our controller design method is associated with general performance indices, it can be compatible with various other control methods.
Future work will involve further analysis and development to explore these broader impacts.

\appendix
\section{Proof of Proposition \ref{thm:LaSalle}} \label{pf:LaSalle}

Let $\LevSetEpsilon:= (\inf_{\State \in {\setB{\DomX}} } \MValF(\State) - \LevSetRAVal)/2 >0$.
Because
$ \LevSetFunc(\LevSetRAVal+\LevSetEpsilon/2)  \cap {\setB{\DomX}} \neq \emptyset$ indicates $\inf_{\State \in {\setB{\DomX}} } \MValF(\State)\leq \LevSetRAVal+\LevSetEpsilon/2 $ contradicting \eqref{eq:LaSalle_cond1}, 
using \eqref{eq:LaSalle_cond1} and  $\{0\} \subseteq \DomX$ provides $\{0\} \subseteq \LevSetFunc(\LevSetTRVal) \subseteq \LevSetFunc(\LevSetRAVal)\subseteq \LevSetFunc(\LevSetRAVal+\LevSetEpsilon/2) \subseteq \DomX \setminus {\setB{\DomX}} \subset \DomX$.
Let $\tmpLevSetDom$ be the interior of $\LevSetFunc(\LevSetRAVal+\LevSetEpsilon/2)$.
Because of $\{0\} \cap {\setB{\DomX}} = \emptyset$ and the continuity of $\MValF(\State)$,
$\tmpLevSetDom$ contains an open neighborhood of $\State=0$ and  is nonempty.
Let us define $\MValFLaSalle: \mathbb{R}^{\DimX} \to \mathbb{R}$ by
\begin{align}
&
\MValFLaSalle(\State) 
\nonumber\\&
:=
\begin{cases}
(    \min\{\max\{ \MValF(\State) - \LevSetTRVal   ,0 \} , \LevSetRAVal  + \LevSetEpsilon - \LevSetTRVal \}  )^{2}  &  (  \State  \in   \DomX   )
\\
(\LevSetRAVal  + \LevSetEpsilon - \LevSetTRVal  )^{2}  & (  \State  \notin  \DomX  )
\end{cases}
.
\end{align}

Firstly, we prove $\MValFLaSalle $ is $C^{1}$ continuous on nonempty open $\tmpLevSetDom$.
For any  $\State \in \tmpLevSetDom$ satisfying $\MValF(\State) > \LevSetTRVal$, 
there exists open neighborhood $\tmpXBall\subset \tmpLevSetDom$ of $\State$ such that 
every $\tempState \in \tmpXBall$ satisfies $\MValF(\tempState) > \LevSetTRVal$.
Thus, we have
$
\partial \MValFLaSalle(\State)/ \partial \State
=
2 ( \MValF(\State) - \LevSetTRVal ) \partial \MValF(\State)/ \partial \State
$.
For any  $\State \in \tmpLevSetDom$ satisfying $\MValF(\State) < \LevSetTRVal$,
we have
$
\partial \MValFLaSalle(\State)/ \partial \State
=
0
$ in the same manner.
If $\MValF(\State) = \LevSetTRVal$ holds,
using the mean value theorem yields
$(\MValFLaSalle(\tempState) - \MValFLaSalle(\State))/\|\tempState-\State\|
=
\max\{ \MValF(\tempState) - \LevSetTRVal  ,0 \}^{2}/\|\tempState-\State\|
\to 0 $ as $\tempState \to \State$. 
These are summarized as $
\partial \MValFLaSalle(\State)/ \partial \State
=
2 \max\{ \MValF(\State) - \LevSetTRVal  ,0 \} \partial \MValF(\State)/ \partial \State
$ and thus $\MValFLaSalle $ is $C^{1}$ continuous on $\tmpLevSetDom$.

Next, let us define the level set
$\LevSetFuncLaSalle(\LevSetVal) := \{ \State \in \mathbb{R}^{\DimX} | \MValFLaSalle(\State) \leq  \LevSetVal \} $.
Comparing $\MValFLaSalle$ with $\MValF$ yields $\LevSetFuncLaSalle( ( \LevSetVal - \LevSetTRVal )^{2} ) = \LevSetFunc(\LevSetVal) \subseteq \DomX$ for any $\LevSetVal \in [\LevSetTRVal , \LevSetRAVal+ \LevSetEpsilon/2]$ and thus this set is bounded.
Because  the continuity of $\MValF(\State)$ indicates $\LevSetFuncLaSalle( ( \LevSetRAVal - \LevSetTRVal )^{2} )= \LevSetFunc(\LevSetRAVal) \subseteq \tmpLevSetDom$, 
$\dot{\MValFLaSalle}(\State)$ is well-defined on the bounded closed set $\LevSetFuncLaSalle( ( \LevSetRAVal - \LevSetTRVal )^{2} )$.
Using the condition \eqref{eq:LaSalle_cond2} implies that
$
\dot{\MValFLaSalle}(\State) 
=
2( \MValF(\State) - \LevSetTRVal ) \dot{\MValF}(\State) 
\leq 
2( \MValF(\State) - \LevSetTRVal)
\HdotV(\State) 
< 0
$
for any $\State \in \LevSetFuncLaSalle( ( \LevSetRAVal - \LevSetTRVal )^{2} ) \setminus \LevSetFuncLaSalle(0)$
and $\dot{\MValFLaSalle}(\State) = 0$ for any $\State \in \LevSetFuncLaSalle(0)$.

Based on these results derived above, applying LaSalle invariance principle to $\MValFLaSalle(\State)$ shows that  
$\RegionOfAttraction=\LevSetFuncLaSalle( ( \LevSetRAVal - \LevSetTRVal )^{2} )$ 
in (\ref{eq:LaSalle_RegionOfAttraction})  and 
$\InvariantSet=\LevSetFuncLaSalle(0 )$ in (\ref{eq:LaSalle_InvariantSet}) 
are a region of attraction and target region, respectively.
This completes the proof.

\section{Proof of Proposition \ref{thm:ex_classH}} \label{pf:ex_classH}	
Setting $\DomBigY= \mathbb{R} $ satisfies {\ASSNLFA} for any continuous $\yBasis: \DomBigX \to \mathbb{R}$.
For  $\sin(\yBasis)$ and $\cos(\yBasis)$, the settings
${\dinf{\NLFunc}}={\ddinf{\NLFunc}}=-1$ and 
${\dsup{\NLFunc}}={\ddsup{\NLFunc}}=1$ satisfy {\ASSNLFB} because of
$|\partial \NLFunc(\yBasis)/\partial \yBasis| \leq 1 $ and $|\partial^{2} \NLFunc(\yBasis)/\partial \yBasis^{2}| \leq 1 $.
For the sigmoid function, 
because $\partial \NLFunc(\yBasis)/\partial \yBasis=  \NLFunc(\yBasis) (1-  \NLFunc(\yBasis))$
\cite[Sec. 1.2.1.6]{Aggarwal18}
and
$\partial^{2} \NLFunc(\yBasis)/\partial \yBasis^{2}
=
(1-2 \NLFunc(\yBasis))
(\partial \NLFunc(\yBasis)/\partial \yBasis)	
$,
{\ASSNLFB} holds by
${\dinf{\NLFunc}}=0$,
${\ddinf{\NLFunc}}=-1/4$,
and ${\dsup{\NLFunc}}={\ddsup{\NLFunc}}=1/4$.
For $\tanh \yBasis$, because $\partial \NLFunc(\yBasis)/\partial \yBasis=   1-  \NLFunc(\yBasis)^{2}$
\cite[Sec. 1.2.1.6]{Aggarwal18}
and
$\partial^{2} \NLFunc(\yBasis)/\partial \yBasis^{2}
=
- 2 \NLFunc(\yBasis)
(\partial \NLFunc(\yBasis)/\partial \yBasis)	
$,
{\ASSNLFB} holds by
${\dinf{\NLFunc}}=0$,
${\dsup{\NLFunc}}=1$,
${\ddinf{\NLFunc}}=-2$,
and ${\ddsup{\NLFunc}}=2$.
This completes the proof.	
%

\section{Proof of Proposition \ref{thm:FSETfunc_properties}} \label{pf:FSETfunc_properties}	

Supposing that  $\Basis \in {\FSETfunc{\NumBasis}{\NumNest}}$ is $C^{2}$ continuous,
any $\Basis \in {\FSETfunc{\NumBasis}{\NumNest+1}}$ is $C^{2}$ continuous because the corresponding ${\subAbasis{\IDEl}}$, ${\subBbasis{\IDEl}}$, and $\NLFunc$ are $C^{2}$ continuous.
Because of mathematical induction, the statement \ref{thm:C2continuous} is proved.	

Next, we prove the statement \ref{thm:FSETfunc_nondec}.
We consider $\NLFunc(\yBasis)=\yBasis$, $ {\subBbasis{1}}(\State) = 1$, and $ {\subBbasis{\IDEl}}(\State) = 0$  for $\IDEl\geq 2$, satisfying $ {\subBbasis{\IDEl}} \in {\FSETfunc{\NumBasis}{1}}$. 
Then, any ${\subAbasis{1}} \in {\FSETfunc{\NumBasis}{1}} $ satisfies ${\subAbasis{1}}=\Basis \in {\FSETfunc{\NumBasis}{2}}$. 
If ${\FSETfunc{\NumBasis}{1}} \subseteq \dots \subseteq {\FSETfunc{\NumBasis}{\NumNest}}$ holds for an $\NumNest$,
the chosen $ {\subBbasis{\IDEl}} \in {\FSETfunc{\NumBasis}{1}}$ are contained in ${\FSETfunc{\NumBasis}{\NumNest}}$ and 
any ${\subAbasis{1}} \in {\FSETfunc{\NumBasis}{\NumNest}} $ satisfies ${\subAbasis{1}}=\Basis \in {\FSETfunc{\NumBasis}{\NumNest+1}}$, This implies ${\FSETfunc{\NumBasis}{1}} \subseteq \dots \subseteq {\FSETfunc{\NumBasis}{\NumNest}} \subseteq {\FSETfunc{\NumBasis}{\NumNest+1}}$. 
Using mathematical induction, the proof is completed.

\section{Proof of Theorem \ref{thm:ex_simple_functions}} \label{pf:ex_simple_functions}

The statement (i) holds clearly because of Remark \ref{rem:bounds_of_directional_derivative}.
We prove the statement (ii) based on mathematical induction.
Suppose that for an $\NumNest$, any $\NumNest$th order polynomial function is contained in ${\FSETfunc{\NumBasis}{\NumNest-1}}$.
According to \eqref{eq:def_Basis} and \eqref{eq:def_yBasis},
any $(\NumNest+1)$th order polynomial function can be described by $\Basis \in {\FSETfunc{\NumBasis}{\NumNest}}$ that is decomposed into appropriate $\NumNest$th order polynomials ${\subAbasis{\IDEl}} \in {\FSETfunc{\NumBasis}{\NumNest-1}}$ for $\IDEl \in \{1, \dots, \DimX\}$,  each component ${\subBbasis{\IDEl}}={\El{\State}{\IDEl}}$ for $\IDEl \in \{1, \dots, \DimX\}$, and constants ${\subAbasis{\IDEl}}$ and ${\subBbasis{\IDEl}}$ for $\IDEl \in \{\DimX+1, \dots, \NumBasis\}$ with $\NLFunc(\yBasis)=\yBasis$.
Therefore, the statement (ii) holds via the mathematical induction with the statement (i).
Finally, the settings of ${\subBbasis{\IDEl}}={\LCombiCoef{\IDEl}}$ and $\NLFunc(\yBasis)=\yBasis$ gives the statement (iii).
This completes the proof.  

\section{Proof of Theorem \ref{thm:ex_SEkernel}} \label{pf:ex_SEkernel}

Let us define 
$\drefState:=\DState{\IDdat}- \State$,
$\GkCoefA:=\tempUnitVec^{\MyTRANSPO} \hypCovMat^{-1}\tempUnitVec>0$,
$\GkCoefB:=\tempUnitVec^{\MyTRANSPO} \hypCovMat^{-1} \drefState$,
and $\GkCoefC:=\drefState^{\MyTRANSPO} \hypCovMat^{-1} \drefState\geq 0$.
We obtain $\GkCoefA \GkCoefC \geq \GkCoefB^{2} $ because
$\GkCoefB \leq \| \sqrKcov \tempUnitVec \| \| \sqrKcov \drefState \| = \GkCoefA^{1/2} \GkCoefC^{1/2}  $ is satisfied by $\sqrKcov^{\MyTRANSPO}\sqrKcov=\hypCovMat^{-1}$.		
Substituting these definitions and relations into \eqref{eq:def_SEkernel} yields
$\kernel{     \tempDist\tempUnitVec + \State   }{\DState{\IDdat}}
=
\kernel{     \tempDist\tempUnitVec    }{\drefState}
$ and
\begin{align}	
\kernel{     \tempDist\tempUnitVec    }{\drefState}
&	
=
\hypMag \exp \big(-(  \tempDist\tempUnitVec - \drefState )^{\MyTRANSPO} \hypCovMat^{-1}  (  \tempDist\tempUnitVec - \drefState ) /2 \big)	
\nonumber\\&
=
\hypMag \exp \big( -\big(  
\GkCoefA \tempDist^{2} - 2 \GkCoefB \tempDist + \GkCoefC
\big)  /2 \big)
\nonumber\\&
=
\hypMag \exp \Big( \frac{- \GkCoefA }{2}\Big(  
\big( \tempDist -  \frac{\GkCoefB }{ \GkCoefA } \big)^{2} 
- \big(\frac{\GkCoefB }{ \GkCoefA } \big)^{2}  
+ \frac{\GkCoefC}{\GkCoefA}
\Big)  \Big)	
.
\end{align}	
The derivatives of ${\kernel{     \tempDist\tempUnitVec    }{\drefState}}$ are given as follows:
\begin{align}	
\frac{\partial}{\partial \tempDist} \kernel{     \tempDist\tempUnitVec    }{\drefState}
&=	
(- \GkCoefA \tempDist + \GkCoefB )
\kernel{     \tempDist\tempUnitVec    }{\drefState} 
.
\\
\frac{\partial^{2}}{\partial \tempDist^{2}} \kernel{     \tempDist\tempUnitVec    }{\drefState}
&=	
\Big(
-\GkCoefA
+
( - \GkCoefA \tempDist + \GkCoefB )^{2}
\Big)
\kernel{     \tempDist\tempUnitVec    }{\drefState}
.
\\
\frac{\partial^{3}}{\partial \tempDist^{3}} \kernel{     \tempDist\tempUnitVec    }{\drefState}
&=	
\Big(
- \GkCoefA
+
( - \GkCoefA \tempDist + \GkCoefB )^{2}
-2 \GkCoefA
\Big)
\nonumber\\&\quad \times
( - \GkCoefA \tempDist + \GkCoefB )
\kernel{     \tempDist\tempUnitVec    }{\drefState}
.
\end{align}
Thus,
we have ${\partial^{2}} \kernel{     \tempDist\tempUnitVec    }{\drefState}/{\partial \tempDist^{2}} \to 0 $ as $\tempDist\to \pm \infty$.
The minimum and maximum of  ${\partial^{2}} \kernel{     \tempDist\tempUnitVec    }{\drefState}/{\partial \tempDist^{2}} $ reduce to its extremums given at ${\SEkerExtPo{1}} := \GkCoefB/\GkCoefA$ and ${\SEkerExtPo{2}}:= (\GkCoefB \pm \sqrt{ 3 \GkCoefA } )/\GkCoefA$:
\begin{align}	
\frac{\partial^{2}}{\partial \tempDist^{2}} \kernel{     {\SEkerExtPo{1}} \tempUnitVec    }{\drefState}
&=	
-\GkCoefA
\hypMag \exp \Big( \frac{- 1 }{2 \GkCoefA }\Big(  
\GkCoefA \GkCoefC   -\GkCoefB^{2}   
\Big)  \Big)
\nonumber\\&
\geq
-
\hypMag
 {\EigMax{\hypCovMat^{-1}}}
.
\\
\frac{\partial^{2}}{\partial \tempDist^{2}} \kernel{     {\SEkerExtPo{2}} \tempUnitVec    }{\drefState}
&=	
2\GkCoefA
\hypMag \exp \Big( \frac{- \GkCoefA }{2}\Big(  
 \frac{3 }{ \GkCoefA }   - \big(\frac{\GkCoefB }{ \GkCoefA }\big)^{2}  + \frac{\GkCoefC}{\GkCoefA}
\Big)  \Big)	
\nonumber\\&
=
2\GkCoefA
\hypMag \exp \Big( \frac{- 1 }{2\GkCoefA}\Big(  
3   \GkCoefA
+
\GkCoefA \GkCoefC   -\GkCoefB^{2}   
\Big)  \Big)
\nonumber\\&
\leq
2
\hypMag 
 {\EigMax{\hypCovMat^{-1}}}
\exp ( {- 3 }/{2} )
,
\end{align}
where we used $\GkCoefA\leq {\EigMax{\hypCovMat^{-1}}}$ for $\|\tempUnitVec\|=1$. 
Note that 
$
\inf_{\drefState, \tempUnitVec, \tempDist }
{ \partial^{2}  \kernel{   \tempDist \tempUnitVec    }{\drefState}   }/{\partial \tempDist^{2}}
\leq
{ \partial^{2} {\kernel{  \tempDist\tempUnitVec + \State    }{\DState{\IDdat}}}     }/{\partial \tempDist^{2}}
|_{ \tempDist=0}
\leq
\sup_{\drefState, \tempUnitVec, \tempDist }
{ \partial^{2}  \kernel{   \tempDist\tempUnitVec    }{\drefState}   }/{\partial \tempDist^{2}}
$ for  $\|\tempUnitVec\|=1$.
Therefore, \eqref{eq:SEkernel_dirddinf} and \eqref{eq:SEkernel_dirddsup} satisfy \eqref{eq:bounds_of_directional_derivative}.
This completes the proof.

\section{Proof of Theorem \ref{thm:ex_GP}} \label{pf:ex_GP}

Theorem \ref{thm:ex_SEkernel} implies $\kernel{ \wildcard }{\DState{\IDdat}} \in {\FSETfunc{\NumBasis}{1}}$.
For each $\IDEl $,	
we choose  
$ {\subAbasis{\IDdat}}=\kernel{ \wildcard }{\DState{\IDdat}}$,
${\subBbasis{\IDdat}}={\El{
		[\DDrift{1},  \dots, \DDrift{\NumData}]
		\kernelMat^{-1}}{\IDEl,\IDdat}} \in {\FSETfunc{\NumBasis}{1}}$ for $\IDdat \leq \NumData$,
${\subAbasis{\IDdat}}={\subBbasis{\IDdat}}=0$ for $\IDdat>\NumData$, 
and $\NLFunc(\yBasis)=\yBasis$.	
Then, ${\El{\GPmean(\State)}{\IDEl}}
=\NLFunc(   
\sum_{\IDbEl=1}^{\NumBasis} 
{\subAbasis{\IDbEl}}(\State) {\subBbasis{\IDbEl}}(\State)    
)$, implying ${\El{\GPmean}{\IDEl}} \in  {\FSETfunc{\NumBasis}{2}}$.
This completes the proof.	

\section{Proof of Theorem \ref{thm:ex_NN}} \label{pf:ex_NN}

For each $\NNidlayer$ and $\IDEl$,	
we choose  
$[ {\subAbasis{1}}, \dots, {\subAbasis{\NNdim{\NNidlayer}}}]^{\MyTRANSPO}={\NNnode{\NNidlayer}}$,
$[ {\subBbasis{1}}, \dots, {\subBbasis{\NNdim{\NNidlayer}}}]^{\MyTRANSPO}={\NNweight{\NNidlayer}{\IDEl}}$,
and ${\subAbasis{\IDbEl}}={\subBbasis{\IDbEl}}=0$ for $\IDbEl>{\NNdim{\NNidlayer}}$.	
If ${\El{\NNnode{\NNidlayer}}{\IDbEl}} \in {\FSETfunc{\NumBasis}{\NNidlayer}}$ holds for $\IDbEl\in\{1,\dots,{\NNdim{\NNidlayer}}\}$,
we obtain ${\El{\NNnode{\NNidlayer+1}(\State)}{\IDEl}}
=\NLFunc(   
\sum_{\IDbEl=1}^{\NumBasis} 
{\subAbasis{\IDbEl}}(\State) {\subBbasis{\IDbEl}}(\State)    
)$ and ${\El{\NNnode{\NNidlayer+1}}{\IDEl}} \in  {\FSETfunc{\NumBasis}{\NNidlayer+1}}$.
In addition, ${\El{\NNnode{1}}{\IDEl}} \in  {\FSETfunc{\NumBasis}{1}}$ holds from Theorem \ref{thm:ex_simple_functions} (i).
Thus,  for $\IDEl\in\{1,\dots,{\NNdim{\NNnumlayer}}\}$,
${\El{\NNnode{\NNnumlayer}}{\IDEl}} \in  {\FSETfunc{\NumBasis}{\NNnumlayer}}$ holds based on mathematical induction.
This completes the proof.

\section{Proof of Lemma \ref{thm:margin_FSET1}} \label{pf:margin_FSET1}

	For each $\IDdom \in \SetALLIDdom$, let us introduce variables ${\tempPPt{\IDPt}{0}}\in \DomSimpX{\IDdom}$, ${\tempPPt{\IDPt}{1}} \in \DomSimpX{\IDdom}$, and ${\tempWPt{\IDPt}} \in [0,1]$ for $\IDPt \in \{0,1,\dots, \DimX\}$.
	For the given $C^{2}$ continuous $\PreliFunc \in {\FSETfunc{\NumBasis}{1}}$, let us define ${\tempPreliFunc{\IDPt}}: [0,1] \to \mathbb{R}$ as follows:
	\begin{align}	
		&
		{\tempPreliFunc{\IDPt}}(  {\tempWPt{}}  )
		:=
		\PreliFunc(  {\tempWPt{}} {\tempPPt{\IDPt}{0}} 	+ (1-{\tempWPt{}}) {\tempPPt{\IDPt}{1}}  )
		.
	\end{align}	
	The relation for the linear interpolation \cite[Theorem 2.1.3]{Villiers2012} is applied to $\PreliFunc$ as follows;
	For some ${\tempWPt{\ast}} \in (0,1)$, we have
	\begin{align}	
		&
		\PreliFunc( 	{\tempWPt{\IDPt}} {\tempPPt{\IDPt}{0}} 	+ (1-{\tempWPt{\IDPt}}) {\tempPPt{\IDPt}{1}}  )
		- (
		{\tempWPt{\IDPt}} \PreliFunc( {\tempPPt{\IDPt}{0}}  ) 
		+ (1-{\tempWPt{\IDPt}}) \PreliFunc( {\tempPPt{\IDPt}{1}}  ) 
		)
		\nonumber\\
		&
		=
		{\tempPreliFunc{\IDPt}}(  {\tempWPt{\IDPt}}  )
		-
		(	 {\tempWPt{\IDPt}} {\tempPreliFunc{\IDPt}}(1) 	+ (1- {\tempWPt{\IDPt}}) {\tempPreliFunc{\IDPt}}(0) 	)
		\nonumber\\
		&
		=
		\frac{-1}{2}
		\frac{\partial^{2} {\tempPreliFunc{\IDPt}}( {\tempWPt{\ast}} )  }{\partial {\tempWPt{}}^{2}} 
		({\tempWPt{\IDPt}}-1)
		(0-{\tempWPt{\IDPt}})
		.\label{eq:LinearInterpolation_0to1}
	\end{align}
	For the given $ {\tempPPt{\IDPt}{0}}$ and $ {\tempPPt{\IDPt}{1}}$, let us define the linear function $\tempDist( {\tempWPt{}}):=   {\tempWPt{}} \| {\tempPPt{\IDPt}{0}} - {\tempPPt{\IDPt}{1}} \|$.
	Because ${\tempPPt{\IDPt}{0}}$ and ${\tempPPt{\IDPt}{1}}$ are contained in convex $\DomSimpX{\IDdom}$, 
	$0 \leq \tempDist({\tempWPt{\ast}})  \leq \| {\tempPPt{\IDPt}{0}} - {\tempPPt{\IDPt}{1}} \| \leq \NewDist$ holds for ${\tempWPt{\ast}} \in (0,1)$.
	Let us define $\tempUnitVec := ( {\tempPPt{\IDPt}{0}} - {\tempPPt{\IDPt}{1}} )/ \| {\tempPPt{\IDPt}{0}} - {\tempPPt{\IDPt}{1}} \| $
	if ${\tempPPt{\IDPt}{0}} \neq {\tempPPt{\IDPt}{1}}$
	and $\tempUnitVec:= [1,0,\dots,0] \in \mathbb{R}^{\DimX}$ if ${\tempPPt{\IDPt}{0}} = {\tempPPt{\IDPt}{1}}$
	so that $\|\tempUnitVec\|=1$ and $\tempDist({\tempWPt{}}) \tempUnitVec = {\tempWPt{}} ({\tempPPt{\IDPt}{0}}-{\tempPPt{\IDPt}{1}})$ hold.
	Because
	$ {\tempPreliFunc{\IDPt}}(  {\tempWPt{\ast}}  )
	=
	\PreliFunc(  {\tempWPt{\ast}} {\tempPPt{\IDPt}{0}} 	+ (1-{\tempWPt{\ast}}) {\tempPPt{\IDPt}{1}}  )
	=
	\PreliFunc( \tempDist({\tempWPt{\ast}}) \tempUnitVec + {\tempPPt{\IDPt}{1}}  ) 
	$
	holds, we obtain
	\begin{align}	
		\frac{\partial^{2} {\tempPreliFunc{\IDPt}}( {\tempWPt{\ast}} )  }{\partial {\tempWPt{}}^{2}} 
		&
		=
		\frac{\partial^{2} \PreliFunc( \tempDist({\tempWPt{\ast}}) \tempUnitVec + {\tempPPt{\IDPt}{1}}  )   }{\partial {\tempWPt{}}^{2}} 
		\nonumber\\
		&
		=
		\frac{ \partial  }{\partial {\tempWPt{}} }
		\Big(
		\frac{ \partial    \PreliFunc( \tempDist({\tempWPt{\ast}}) \tempUnitVec + {\tempPPt{\IDPt}{1}}  )   }{\partial \tempDist }
		\frac{ \partial  \tempDist({\tempWPt{\ast}})  }{\partial  {\tempWPt{}} }
		\Big)
		\nonumber\\&
		=
		\frac{ \partial^{2} \PreliFunc( \tempDist({\tempWPt{\ast}}) \tempUnitVec + {\tempPPt{\IDPt}{1}}  )   }{\partial \tempDist^{2}}
		\| {\tempPPt{\IDPt}{0}} - {\tempPPt{\IDPt}{1}} \|^{2}
		.
	\end{align}
	Note that  $(\tempDist({\tempWPt{\ast}}) \tempUnitVec + {\tempPPt{\IDPt}{1}} ) \in \DomSimpX{\IDdom}$ holds because it is a point in the line segment between ${\tempPPt{\IDPt}{0}}  $ and ${\tempPPt{\IDPt}{1}}  $.
	Therefore, 
	if	$\sup_{  {\tempPPt{\IDPt}{0}}  , {\tempPPt{\IDPt}{1}} , {\tempWPt{\ast}} }
	{ \partial^{2} \PreliFunc( \tempDist({\tempWPt{\ast}}) \tempUnitVec + {\tempPPt{\IDPt}{1}}  )   }/{\partial \tempDist^{2}} \geq 0$ holds, using the inequality $\| {\tempPPt{\IDPt}{0}} - {\tempPPt{\IDPt}{1}} \| \leq \NewDist$ yields
	\begin{align}	
		\frac{\partial^{2} {\tempPreliFunc{\IDPt}}( {\tempWPt{\ast}} )  }{\partial {\tempWPt{}}^{2}} 
		&
		\leq
		\sup_{  {\tempPPt{\IDPt}{0}}  , {\tempPPt{\IDPt}{1}}  \in \DomSimpX{\IDdom}, {\tempWPt{\ast}} \in (0,1) }
		\frac{ \partial^{2} \PreliFunc( \tempDist({\tempWPt{\ast}}) \tempUnitVec + {\tempPPt{\IDPt}{1}}  )   }{\partial \tempDist^{2}}
		\NewDist^{2}
		\nonumber\\
		&
		\leq
		\NewDist^{2}
		\sup_{ 
			\State \in \DomSimpX{\IDdom}
			,\tempUnitVec 
			\in \{  \tempUnitVec \in  \mathbb{R}^{\DimX} 
			|
			\|\tempUnitVec\|=1
			\}
		}
		\frac{ \partial^{2} \PreliFunc(  \tempDist\tempUnitVec + \State  ) }{\partial \tempDist^{2}}
		\Big|_{ \tempDist=0}
		\nonumber\\&
		\leq	
		\NewDist^{2}
		{\dirddsup{\PreliFunc}}
		.\label{eq:Bound_partial2_tempPreliFunc_pre}
	\end{align}	
	In a similar manner, $\NewDist^{2} {\dirddinf{\PreliFunc}}  \leq {\partial^{2} {\tempPreliFunc{\IDPt}}( {\tempWPt{\ast}} )  }/{\partial {\tempWPt{}}^{2}} $ holds
	if $\inf_{  {\tempPPt{\IDPt}{0}}  , {\tempPPt{\IDPt}{1}} , {\tempWPt{\ast}} }
	{ \partial^{2} \PreliFunc( \tempDist({\tempWPt{\ast}}) \tempUnitVec + {\tempPPt{\IDPt}{1}}  )   }/{\partial \tempDist^{2}} \leq 0$.
	In addition, $0\leq ({\tempWPt{\IDPt}}-1)(0-{\tempWPt{\IDPt}}) \leq 1/4$ holds for ${\tempWPt{\IDPt}} \in [0,1]$.
	Therefore, we obtain bounds of the term in  \eqref{eq:LinearInterpolation_0to1}:
	\begin{align}	
		\frac{-\NewDist^{2}}{8}
		\max\{0, {\dirddsup{\PreliFunc}} \}
		&
		\leq
		\frac{-1}{2}
		\frac{\partial^{2} {\tempPreliFunc{\IDPt}}( {\tempWPt{\ast}} )  }{\partial {\tempWPt{}}^{2}} 
		({\tempWPt{\IDPt}}-1)
		(0-{\tempWPt{\IDPt}})
		\nonumber\\
		&
		\leq
		\frac{\NewDist^{2}}{8}
		\max\{0, -{\dirddinf{\PreliFunc}} \}
		.\label{eq:Bound_partial2_tempPreliFunc}
	\end{align}	
	The combination of \eqref{eq:LinearInterpolation_0to1} with \eqref{eq:Bound_partial2_tempPreliFunc} is represented as follows:
	\begin{align}	
		&
		\PreliFunc( 	{\tempWPt{\IDPt}} {\tempPPt{\IDPt}{0}} 	+ (1-{\tempWPt{\IDPt}}) {\tempPPt{\IDPt}{1}}  )
		\nonumber\\
		&\leq
		(
		{\tempWPt{\IDPt}} \PreliFunc( {\tempPPt{\IDPt}{0}}  ) 
		+ (1-{\tempWPt{\IDPt}}) \PreliFunc( {\tempPPt{\IDPt}{1}}  ) 
		)
		+
		\frac{\NewDist^{2}}{8}
		\max\{0, -{\dirddinf{\PreliFunc}} \}
		.\label{eq:LinearInterpolation_0to1_B}
	\end{align}
	If  $\DimX = 1$ holds, replacing ${\tempPPt{\IDPt}{0}}$, ${\tempPPt{\IDPt}{1}}$, and ${\tempWPt{\IDPt}}$ in \eqref{eq:LinearInterpolation_0to1_B} with ${\SimpPt{\IDdom}{1}}$, ${\SimpPt{\IDdom}{2}}$, and ${\El{\SimpVecCoef}{1}}$, respectively,  completes the proof. 
	Otherwise, the relation \eqref{eq:LinearInterpolation_0to1_B} is extended for the case of $\DimX \geq 2$ as follows.
	Let us define 
	\begin{align}	
		{\tempPPt{\IDPt-1}{1}} :=	{\tempWPt{\IDPt}} {\tempPPt{\IDPt}{0}} 	+ (1-{\tempWPt{\IDPt}}) {\tempPPt{\IDPt}{1}} 
		,\;
		\forall\IDPt \in \{1,2,\dots, \DimX\}
		.\label{eq:def_tempPPt}
	\end{align}
	This yields the following relation:
	\begin{align}	
		\PreliFunc(  {\tempPPt{0}{1}}    )
		&\leq
		{\tempWPt{1}} \PreliFunc( {\tempPPt{1}{0}}  ) 
		+ (1-{\tempWPt{1}}) \PreliFunc( {\tempPPt{1}{1}}  ) 
		+
		\frac{\NewDist^{2}}{8}
		\max\{0, -{\dirddinf{\PreliFunc}} \}
		\nonumber\\
		&\leq
		{\tempWPt{1}} \PreliFunc( {\tempPPt{1}{0}}  ) 
		+  (1-{\tempWPt{1}}){\tempWPt{2}} \PreliFunc( {\tempPPt{2}{0}}  ) 
		\nonumber\\&\quad
		+  (1-{\tempWPt{1}})(1-{\tempWPt{2}})  \PreliFunc( {\tempPPt{2}{1}}  ) 
		\nonumber\\&\quad
		+
		(1+ (1-{\tempWPt{1}})  )
		\frac{\NewDist^{2}}{8}
		\max\{0, -{\dirddinf{\PreliFunc}} \}
		\nonumber\\
		&\qquad \vdots
		\nonumber\\
		&\leq
		{\tempWPt{1}} 
		\PreliFunc( {\tempPPt{1}{0}}  ) 
		+
		\Bigg( 
		\sum_{\IDPt=2}^{\DimX}
		\Big(  
		{\tempWPt{\IDPt}} \prod_{\IDbPt=1}^{\IDPt-1} (1-{\tempWPt{\IDbPt}})  
		\Big) 
		\PreliFunc( {\tempPPt{\IDPt}{0}}  ) 
		\Bigg)
		\nonumber\\&\quad
		+
		\Big( 
		\prod_{\IDbPt=1}^{\DimX}(1-{\tempWPt{\IDbPt}})
		\Big) 
		\PreliFunc( {\tempPPt{\DimX}{1}}) 		
		+
		\frac{\DimX\NewDist^{2}}{8}
		\max\{0, -{\dirddinf{\PreliFunc}} \}
		. \label{eq:LI_multi_dim}
	\end{align}
	Because the definition \eqref{eq:def_tempPPt} does not constrain ${\tempPPt{\DimX}{1}}$, ${\tempPPt{\IDPt}{0}}$, and ${\tempWPt{\IDPt}}$ for $\IDPt \in \{1,\dots, \DimX\}$, we choose 
	${\tempPPt{\DimX}{1}} ={\SimpPt{\IDdom}{\DimX+1}}$ and 	${\tempPPt{\IDPt}{0}} ={\SimpPt{\IDdom}{\IDPt}}$ and set ${\tempWPt{\IDPt}} \in [0,1]$ such that 
	$\El{\SimpVecCoef}{\IDPt} =
	{\tempWPt{\IDPt}}\sum_{\IDbPt=\IDPt}^{\NumSimp}\El{\SimpVecCoef}{\IDbPt}$
	holds.
	For these settings, we prove the following relation for all $\IDPt \in \{2,\dots, \DimX\}$:
	\begin{align}
		1-  \sum_{\IDbPt=1}^{\IDPt-1}\El{\SimpVecCoef}{\IDbPt}
		=\prod_{\IDbPt=1}^{\IDPt-1} (1-{\tempWPt{\IDbPt}})
		.\label{eq:tempWPt_property}
	\end{align}
	For $\IDPt=2$, \eqref{eq:tempWPt_property} holds because 
	$\El{\SimpVecCoef}{1} =
	{\tempWPt{1}}\sum_{\IDbPt=1}^{\NumSimp}\El{\SimpVecCoef}{\IDbPt}={\tempWPt{1}}$
	indicates $ (1- \El{\SimpVecCoef}{1} )= (1- {\tempWPt{1}}  )$ in \eqref{eq:tempWPt_property}.
	Supposing that \eqref{eq:tempWPt_property} holds for a given $\IDPt \in \{2,\dots, \DimX-1\}$,
	using 
	$\sum_{\IDbPt=\IDPt}^{\NumSimp}\El{\SimpVecCoef}{\IDbPt}=1-  \sum_{\IDbPt=1}^{\IDPt-1}\El{\SimpVecCoef}{\IDbPt}$
	yields
	\begin{align}
		\El{\SimpVecCoef}{\IDPt} 
		&
		=
		{\tempWPt{\IDPt}}
		\Big(
		1-  \sum_{\IDbPt=1}^{\IDPt-1}\El{\SimpVecCoef}{\IDbPt}
		\Big)
		= {\tempWPt{\IDPt}} \prod_{\IDbPt=1}^{\IDPt-1} (1-{\tempWPt{\IDbPt}})
		,\label{eq:tempWPt_property_B}
		\\
		1-  \sum_{\IDbPt=1}^{\IDPt}\El{\SimpVecCoef}{\IDbPt}
		&=   \prod_{\IDbPt=1}^{\IDPt-1} (1-{\tempWPt{\IDbPt}})   -   \El{\SimpVecCoef}{\IDPt}  
		= \prod_{\IDbPt=1}^{\IDPt} (1-{\tempWPt{\IDbPt}})
		.\label{eq:tempWPt_property_C}
	\end{align}
	This implies that \eqref{eq:tempWPt_property} holds for $\IDPt+1$.
	Therefore,  \eqref{eq:tempWPt_property} holds for all $\IDPt \in \{2,\dots, \DimX\}$ based on mathematical induction.
	In addition, \eqref{eq:tempWPt_property_B} holds for all $\IDPt \in \{2,\dots, \DimX\}$
	and
	${\El{\SimpVecCoef}{\NumSimp}}= 1-  \sum_{\IDbPt=1}^{\DimX}\El{\SimpVecCoef}{\IDbPt} =  \prod_{\IDbPt=1}^{\DimX} (1-{\tempWPt{\IDbPt}})$ holds because of \eqref{eq:tempWPt_property_C}.
	Using these relations, we obtain 
	\begin{align}	
		{\tempPPt{0}{1}}
		&=
		{\tempWPt{1}} 
		{\tempPPt{1}{0}} 
		\nonumber\\&\quad
		+
		\sum_{\IDPt=2}^{\DimX}
		{\tempWPt{\IDPt}}
		\Big( \prod_{\IDbPt=1}^{\IDPt-1} (1-{\tempWPt{\IDbPt}})  \Big)
		{\tempPPt{\IDPt}{0}} 
		+
		\Big(   \prod_{\IDbPt=1}^{\DimX}(1-{\tempWPt{\IDbPt}})   \Big)
		{\tempPPt{\DimX}{1}} 		
		\nonumber\\&
		=
		\sum_{\IDPt=1}^{\NumSimp}\El{\SimpVecCoef}{\IDPt}\SimpPt{\IDdom}{\IDPt}
		=
		\SimpState
		.\label{eq:tempWPt_property_D}
	\end{align}
	Therefore, because of Definition \ref{def:LI}, substituting \eqref{eq:tempWPt_property_D} into \eqref{eq:LI_multi_dim} gives
		$\PreliFunc(\State) - \LI{\PreliFunc}(\State)
	\leq
		({ \DimX \NewDist^{2}}/{8})
		\max\{0, -{\dirddinf{\PreliFunc}} \}
		$.
	The lower bound is derived in the same manner.
	This completes the proof.

\section{Proof of Lemma \ref{thm:margins_for_inner_products}} \label{pf:margins_for_inner_products}

	Multiplying \eqref{eq:def_general_bounds} by $\subAbasis{\IDEl}(\State)$ and replacing $\PreliFunc$ with $\subBbasis{\IDEl}$ derive 
	$
	-\subAbasis{\IDEl}(\State) {\zzzLM{\subBbasis{\IDEl}}}
	\leq
	\subAbasis{\IDEl}(\State)
	(	\subBbasis{\IDEl}(\State)	-\LI{\subBbasis{\IDEl}}(\State)	)
	\leq
	\subAbasis{\IDEl}(\State)	{\zzzUM{\subBbasis{\IDEl}}}
	$ if $\subAbasis{\IDEl}(\State)\geq 0$,
	where the signs of these inequalities are reversed if  $\subAbasis{\IDEl}(\State)< 0$.
	These inequalities are integrated as follows, using $\subAbasis{\IDEl}(\State) \leq {\zzzUB{\subAbasis{\IDEl}}}$ and $- \subAbasis{\IDEl}(\State) \leq - {\zzzLB{\subAbasis{\IDEl}}}$:
	\begin{align}
	-
	{\zzzBLB{\subBbasis{\IDEl}}{\subAbasis{\IDEl}}}
	&
	\leq
	-
	\max \{   {\zzzLM{\subBbasis{\IDEl}}} \subAbasis{\IDEl}(\State)
	,       -  {\zzzUM{\subBbasis{\IDEl}}} \subAbasis{\IDEl}(\State) \}
	\nonumber\\&
	=
	\min \{ -  {\zzzLM{\subBbasis{\IDEl}}}  \subAbasis{\IDEl}(\State)
	,         {\zzzUM{\subBbasis{\IDEl}}} \subAbasis{\IDEl}(\State)    \}
	\nonumber\\&
	\leq
	\subAbasis{\IDEl}(\State) 
	(	\subBbasis{\IDEl}(\State)	-\LI{\subBbasis{\IDEl}}(\State)	)
	\nonumber\\&
	\leq
	\max \{  {\zzzUM{\subBbasis{\IDEl}}} \subAbasis{\IDEl}(\State) 
	,       - {\zzzLM{\subBbasis{\IDEl}}}  \subAbasis{\IDEl}(\State) \}
	\nonumber\\&
	\leq
	{\zzzBUB{\subBbasis{\IDEl}}{\subAbasis{\IDEl}}}
	.	\label{eq:bounds_inequality_with_coef_temp1}
	\end{align}	
	In the same manner,
	because of $\LI{\subBbasis{\IDEl}}(\State) \leq {\zzzUB{\subBbasis{\IDEl}}}$ and $-\LI{\subBbasis{\IDEl}}(\State) \leq -{\zzzLB{\subBbasis{\IDEl}}}$,
	multiplying \eqref{eq:def_general_bounds} by $\LI{\subBbasis{\IDEl}}(\State)$ and replacing $\PreliFunc$ in (\ref{eq:def_general_bounds}) with $\subAbasis{\IDEl}(\State)$ yield the following relation
	\begin{align}
	-
	{\zzzBLB{\subAbasis{\IDEl}}{\subBbasis{\IDEl}}}
	&
	\leq
	\LI{\subBbasis{\IDEl}}(\State)
	(	\subAbasis{\IDEl}(\State) - \LI{\subAbasis{\IDEl}}(\State) )
	\leq
	{\zzzBUB{\subAbasis{\IDEl}}{\subBbasis{\IDEl}}}	
	.	\label{eq:bounds_inequality_with_coef_temp2}
	\end{align}
	The sum of (\ref{eq:bounds_inequality_with_coef_temp1}) and (\ref{eq:bounds_inequality_with_coef_temp2}) is summarized as 
	\begin{align}
	-
	(
	{\zzzBLB{\subAbasis{\IDEl}}{\subBbasis{\IDEl}}}
	+
	{\zzzBLB{\subBbasis{\IDEl}}{\subAbasis{\IDEl}}}
	)
	&
	\leq
	\subAbasis{\IDEl}(\State)
	\subBbasis{\IDEl}(\State)	
	-
	\LI{\subAbasis{\IDEl}}(\State)
	\LI{\subBbasis{\IDEl}}(\State)
	\nonumber\\&
	\leq
	(
	{\zzzBUB{\subAbasis{\IDEl}}{\subBbasis{\IDEl}}}
	+
	{\zzzBUB{\subBbasis{\IDEl}}{\subAbasis{\IDEl}}}
	)
	.	\label{eq:bounds_inequality_with_coef_result1}
	\end{align}

	Next, for any real numbers ${\tempA{\IDPt}}$, ${\tempB{\IDPt}}$, and ${\tempWPt{\IDPt}}$ for $\IDPt \in \{1,\dots,\NumSimp \}$ satisfying $\sum_{\IDPt=1 }^{\NumSimp}{\tempWPt{\IDPt}}=1$, we obtain the following relation:
	\begin{align}	
	&
	\sum_{\IDPt=1 }^{\NumSimp}{\tempWPt{\IDPt}}  {\tempA{\IDPt}}
	\sum_{\IDbPt=1}^{\NumSimp}{\tempWPt{\IDbPt}} {\tempB{\IDbPt}}
	-
	\sum_{\IDPt=1}^{\NumSimp}{\tempWPt{\IDPt}} {\tempA{\IDPt}}{\tempB{\IDPt}}
	\nonumber\\&
	=
	\sum_{\IDPt=1}^{\NumSimp} {\tempWPt{\IDPt}}( {\tempWPt{\IDPt}} -1 ) {\tempA{\IDPt}}{\tempB{\IDPt}}
	+
	\sum_{\IDPt=1}^{\NumSimp}
	\sum_{\IDbPt=1,\IDbPt\neq \IDPt}^{\NumSimp}
	{\tempWPt{\IDPt}}{\tempWPt{\IDbPt}} {\tempA{\IDPt}}{\tempB{\IDbPt}}
	\nonumber\\&
	=
	\sum_{\IDPt=1}^{\NumSimp} {\tempWPt{\IDPt}} 
	\Big( - \sum_{\IDbPt=1,\IDbPt \neq \IDPt}^{\NumSimp} {\tempWPt{\IDbPt}} \Big) 
	{\tempA{\IDPt}}{\tempB{\IDPt}}
	+
	\sum_{\IDPt=1}^{\NumSimp}
	\sum_{\IDbPt=1,\IDbPt\neq \IDPt}^{\NumSimp}
	{\tempWPt{\IDPt}}{\tempWPt{\IDbPt}} {\tempA{\IDPt}}{\tempB{\IDbPt}}
	\nonumber\\&
	=
	\sum_{\IDPt=1}^{\NumSimp}\sum_{\IDbPt=1,\IDbPt\neq \IDPt}^{\NumSimp}
	{\tempWPt{\IDPt}} {\tempWPt{\IDbPt}} {\tempA{\IDPt}}( {\tempB{\IDbPt}} - {\tempB{\IDPt}} )
	\nonumber\\&
	=
	\sum_{\IDPt=1}^{\NumSimp}\sum_{\IDbPt=\IDPt+1}^{\NumSimp}
	{\tempWPt{\IDPt}} {\tempWPt{\IDbPt}} {\tempA{\IDPt}}( {\tempB{\IDbPt}} - {\tempB{\IDPt}} )
	+
	\sum_{\IDPt=1}^{\NumSimp}\sum_{\IDbPt=1}^{\IDPt-1}
	{\tempWPt{\IDPt}} {\tempWPt{\IDbPt}} {\tempA{\IDPt}}( {\tempB{\IDbPt}} - {\tempB{\IDPt}} )
	\nonumber\\&
	=
	\sum_{\IDPt=1}^{\NumSimp}\sum_{\IDbPt=\IDPt+1}^{\NumSimp}
	{\tempWPt{\IDPt}} {\tempWPt{\IDbPt}} {\tempA{\IDPt}}( {\tempB{\IDbPt}} - {\tempB{\IDPt}} )
	+
	\sum_{\IDbPt=1}^{\NumSimp}\sum_{\IDPt=1}^{\IDbPt-1}
	{\tempWPt{\IDbPt}} {\tempWPt{\IDPt}} {\tempA{\IDbPt}}( {\tempB{\IDPt}} - {\tempB{\IDbPt}} )
	\nonumber\\&
	=
	\sum_{\IDPt=1}^{\NumSimp}\sum_{\IDbPt=\IDPt+1}^{\NumSimp}
	{\tempWPt{\IDPt}}{\tempWPt{\IDbPt}} ( {\tempA{\IDPt}} - {\tempA{\IDbPt}} )( {\tempB{\IDbPt}} - {\tempB{\IDPt}} )
	, \label{eq:linear_interplation_transform_ex}
	\end{align} 
	where the last equality is satisfied because 
	$\sum_{\IDPt=1}^{\NumSimp}\sum_{\IDbPt=\IDPt+1}^{\NumSimp} $ and 
	$\sum_{\IDbPt=1}^{\NumSimp}\sum_{\IDPt=1}^{\IDbPt-1}$ 
	are equivalent to $\sum_{\IDbPt=1}^{\NumSimp}\sum_{\IDPt=1, \IDPt<\IDbPt}^{\NumSimp}$.
	By replacing ${\tempA{\IDPt}}$, ${\tempB{\IDPt}}$, and ${\tempWPt{\IDPt}}$ in \eqref{eq:linear_interplation_transform_ex} with ${\subAbasis{\IDEl}(\SimpPt{\IDdom}{\IDPt})}$, ${\subBbasis{\IDEl}(\SimpPt{\IDdom}{\IDPt})}$, and ${\El{\SimpVecCoef}{\IDPt}}$, respectively, and by using the sampling-based gradients in Definition \ref{def:collection},
	we obtain 
	\begin{align}	
	&
	\Big|
	\LI{\subAbasis{\IDEl}}(\SimpState)
	\LI{\subBbasis{\IDEl}}(\SimpState)
	-
	\sum_{\IDPt=1}^{\NumSimp}\El{\SimpVecCoef}{\IDPt} \subAbasis{\IDEl}(\SimpPt{\IDdom}{\IDPt})\subBbasis{\IDEl}(\SimpPt{\IDdom}{\IDPt})
	\Big|
	\nonumber\\&
	=
	\Big|
	\sum_{\IDPt=1}^{\NumSimp}\sum_{\IDbPt=\IDPt+1}^{\NumSimp}
	\El{\SimpVecCoef}{\IDPt}\El{\SimpVecCoef}{\IDbPt} 
	\nonumber\\&\qquad\times
	( \subAbasis{\IDEl}(\SimpPt{\IDdom}{\IDPt}) - \subAbasis{\IDEl}(\SimpPt{\IDdom}{\IDbPt}) )
	( \subBbasis{\IDEl}(\SimpPt{\IDdom}{\IDbPt}) - \subBbasis{\IDEl}(\SimpPt{\IDdom}{\IDPt}) )
	\Big|
	\nonumber\\&
	\leq
	\NewDist^{2} {\zzzLipConst{\subAbasis{\IDEl}}}{\zzzLipConst{\subBbasis{\IDEl}}}
	, \label{eq:linear_interplation_transform}
	\end{align}
	because $0 \leq \sum_{\IDPt=1}^{\NumSimp}\sum_{\IDbPt=\IDPt+1}^{\NumSimp} \El{\SimpVecCoef}{\IDPt}\El{\SimpVecCoef}{\IDbPt} \leq \sum_{\IDPt=1}^{\NumSimp}\sum_{\IDbPt=1}^{\NumSimp} \El{\SimpVecCoef}{\IDPt}\El{\SimpVecCoef}{\IDbPt} = 1$ holds.
	By substituting (\ref{eq:linear_interplation_transform}) into (\ref{eq:bounds_inequality_with_coef_result1}), we obtain
	\begin{align}
	&
	-
	(
	{\zzzBLB{\subAbasis{\IDEl}}{\subBbasis{\IDEl}}}
	+
	{\zzzBLB{\subBbasis{\IDEl}}{\subAbasis{\IDEl}}}
	+
	\NewDist^{2} {\zzzLipConst{\subAbasis{\IDEl}}}{\zzzLipConst{\subBbasis{\IDEl}}}
	)
	\nonumber\\&
	\leq
	\subAbasis{\IDEl}(\SimpState)
	\subBbasis{\IDEl}(\SimpState)	
	-
	\sum_{\IDPt=1}^{\NumSimp}\El{\SimpVecCoef}{\IDPt} 
	\subAbasis{\IDEl}(\SimpPt{\IDdom}{\IDPt})\subBbasis{\IDEl}(\SimpPt{\IDdom}{\IDPt})
	\nonumber\\&
	\leq
	(
	{\zzzBUB{\subAbasis{\IDEl}}{\subBbasis{\IDEl}}}
	+
	{\zzzBUB{\subBbasis{\IDEl}}{\subAbasis{\IDEl}}}
	+
	\NewDist^{2} {\zzzLipConst{\subAbasis{\IDEl}}}{\zzzLipConst{\subBbasis{\IDEl}}}
	)
	.	\label{eq:bounds_inequality_with_coef_result2}
	\end{align}
	Therefore, ${\zzzLM{\yBasis}}$ and ${\zzzUM{\yBasis}}$ given in (\ref{eq:chi_bound_L}) and (\ref{eq:chi_bound_U}) are indeed the margins because substituting them into the sum of (\ref{eq:bounds_inequality_with_coef_result2}) with respect to $\IDEl \in \{1,\dots,\NumBasis\}$ yields
	\begin{align}
	&
	\sum_{\IDEl=1}^{\NumBasis}
	\subAbasis{\IDEl}(\SimpState)\subBbasis{\IDEl}(\SimpState)
	-
	\sum_{\IDPt=1}^{\NumSimp}\El{\SimpVecCoef}{\IDPt} 
	\sum_{\IDEl=1}^{\NumBasis}
	\subAbasis{\IDEl}(\SimpPt{\IDdom}{\IDPt})\subBbasis{\IDEl}(\SimpPt{\IDdom}{\IDPt})
	\nonumber\\
	&
	=	
	\yBasis(\SimpState)	-	\LI{\yBasis}(\SimpState)
	\leq
	{\zzzUM{\yBasis}}
	.	  \label{eq:bounds_inequality_with_coef_result3}
	\end{align}	
	In addition, if $\zzzLM{\subAbasis{\IDEl}}$, $\zzzUM{\subAbasis{\IDEl}}$, $\zzzLM{\subBbasis{\IDEl}}$, and $\zzzUM{\subBbasis{\IDEl}}$ for $\IDEl \in \{1,\dots,\NumBasis\}$ are $\MyOrder(\NewDist^{2})$,
	${\zzzBLB{\subAbasis{\IDEl}}{\subBbasis{\IDEl}}}$,
	${\zzzBLB{\subBbasis{\IDEl}}{\subAbasis{\IDEl}}}$,
	${\zzzBUB{\subAbasis{\IDEl}}{\subBbasis{\IDEl}}}$,
	${\zzzBUB{\subBbasis{\IDEl}}{\subAbasis{\IDEl}}}$, and
	$\NewDist^{2} {\zzzLipConst{\subAbasis{\IDEl}}}{\zzzLipConst{\subBbasis{\IDEl}}}$ are $\MyOrder(\NewDist^{2})$ 
	because for any $\PreliFunc  \in \{ {\subAbasis{\IDEl}}, {\subBbasis{\IDEl}} \}$, 
	${\zzzLB{\PreliFunc}}$, ${\zzzUB{\PreliFunc}}$, and ${\zzzLipConst{\PreliFunc}}$ are bounded  regardless of $\NewDist$ as described in Remark \ref{rem:collection}.
	This implies that ${\zzzLM{\yBasis}}$ and  ${\zzzUM{\yBasis}}$ are $\MyOrder(\NewDist^{2})$. 
	This completes the proof.

\section{Proof of Lemma \ref{thm:margins_for_nonlinear_mappings}} \label{pf:margins_for_nonlinear_mappings}

	We use the proof of Lemma \ref{thm:margin_FSET1} in Appendix \ref{pf:margin_FSET1} by replacing $\PreliFunc$ and $\DomSimpX{\IDdom}$ with $\NLFunc$ and 
	$\DomLIEachY:=\{ \LI{\yBasis}(\State) \in \mathbb{R} | \State \in \DomSimpX{\IDdom} \}
	=[\min_{\IDPt}  \yBasis(\SimpPt{\IDdom}{\IDPt})    ,
	\max_{\IDPt}  \yBasis(\SimpPt{\IDdom}{\IDPt}) ]$, 
	respectively, that is, ${\tempPPt{\IDPt}{0}}$ and ${\tempPPt{\IDPt}{1}}$ are one-dimensional and included in $\DomLIEachY$.
	Note that $\| {\tempPPt{\IDPt}{0}} - {\tempPPt{\IDPt}{1}} \| \leq \NewDist {\zzzLipConst{\yBasis}} $ holds because of the definitions of ${\zzzLipConst{\yBasis}}$ and $\DomLIEachY$.
	Recall that $\NLFunc  \in {\FSETNLFunc{\yBasis}}$ is $C^{2}$ continuous on the open set $\DomBigY$ containing 
	$\DomLIY:=\{ \LI{\yBasis}(\State) \in \mathbb{R} | \State \in \DomX \}
	=[\min_{\IDdom \in \SetALLIDdom,\IDPt}  \yBasis(\SimpPt{\IDdom}{\IDPt})    ,	
		\max_{\IDdom \in \SetALLIDdom,\IDPt}  \yBasis(\SimpPt{\IDdom}{\IDPt}) ]
      \supseteq \DomLIEachY$.
	Considering one-dimensional version of \eqref{eq:bounds_of_directional_derivative} gives
	\begin{align}
		\forall
		\tempyBasis \in \DomLIY
		,\;
		\forall
		\tempUnitVec 
		\in \{ -1, 1	\}
		,
		&
		\nonumber\\
		\frac{ \partial^{2} \NLFunc(  \tempDist\tempUnitVec + \tempyBasis  ) }{\partial \tempDist^{2}}
		\Big|_{ \tempDist=0}
		&=
		\frac{ \partial }{\partial \tempDist }
		\Big(
		\frac{ \partial \NLFunc(  \tempDist\tempUnitVec + \tempyBasis  ) }{\partial \yBasis}
		\frac{ \partial (  \tempDist\tempUnitVec + \tempyBasis  ) }{\partial \tempDist}
		\Big)
		\Big|_{ \tempDist=0}
		\nonumber\\&=
		\frac{ \partial^{2} \NLFunc(  \tempDist\tempUnitVec + \tempyBasis  ) }{\partial \yBasis^{2}}
		\tempUnitVec^{2}
		\Big|_{ \tempDist=0}
		\nonumber\\&
		\leq
		\sup_{ \tempyBasis \in \DomLIY}
		\frac{\partial^{2}  \NLFunc( \tempyBasis  ) }{\partial  \yBasis^{2}} 
		\leq
		{\ddsup{\NLFunc}}
		.\label{eq:bounds_of_directional_derivative_one_dim}
	\end{align}	
If	$\sup_{  {\tempPPt{\IDPt}{0}}  , {\tempPPt{\IDPt}{1}}, {\tempWPt{\ast}} \in (0,1)  }
{ \partial^{2} \NLFunc( \tempDist({\tempWPt{\ast}}) \tempUnitVec + {\tempPPt{\IDPt}{1}}  )   }/{\partial \tempDist^{2}} \geq 0$ holds,
 \eqref{eq:Bound_partial2_tempPreliFunc_pre} is then modified as 
  $
 {\partial^{2} {\tempPreliFunc{\IDPt}}( {\tempWPt{\ast}} )  }/{\partial {\tempWPt{}}^{2}} 
 \leq	
 (\NewDist {\zzzLipConst{\yBasis}})^{2}
 {\ddsup{\NLFunc}}
 $,
	where the lower bound $(\NewDist {\zzzLipConst{\yBasis}})^{2}{\ddinf{\NLFunc}}$ is obtained in the same manner.
	The settings of 
	${\tempPPt{\DimX}{1}} =\yBasis( {\SimpPt{\IDdom}{\DimX+1}} )$ and ${\tempPPt{\IDPt}{0}} =\yBasis( {\SimpPt{\IDdom}{\IDPt}} )$
	yield
	${\tempPPt{0}{1}}=\sum_{\IDPt=1}^{\NumSimp}\El{\SimpVecCoef}{\IDPt} \yBasis( {\SimpPt{\IDdom}{\IDPt}} ) = \LI{\yBasis}(\SimpState)$
	based on \eqref{eq:tempWPt_property_D}.
	Therefore, we obtain the following relation:
	\begin{align}	
	&- 
	\frac{ \DimX (\NewDist{\zzzLipConst{\yBasis}})^{2}}{8}
	\max\{0, {\ddsup{\NLFunc}} \}
	\nonumber\\&
	\leq
	\NLFunc \circ  \LI{\yBasis}(\SimpState) 
	-
	\LI{( \NLFunc \circ  {\yBasis} )} (\SimpState)  
	\nonumber\\&
	\leq              
	\frac{ \DimX (\NewDist{\zzzLipConst{\yBasis}})^{2}}{8}
	\max\{0, -{\ddinf{\NLFunc}} \}
	.\label{eq:ineq_margins_NLFunc_partA}
	\end{align}

	Next, for some ${\tempWPt{\ast}} \in (   0, 1   )$,
	letting $\MeanValueThmY:={\tempWPt{\ast}} {\yBasis}(\SimpState)  + (1-{\tempWPt{\ast}})\LI{\yBasis}(\SimpState) \in \DomY$,
	the mean value theorem gives
	\begin{align}
	&
	\NLFunc \circ  {\yBasis}(\SimpState) 
	-
	\NLFunc \circ  \LI{\yBasis}  (\SimpState) 
	\nonumber\\
	&
	=
	\frac{\partial \NLFunc( 
		\MeanValueThmY
		)}{\partial \yBasis}
	( {\yBasis}(\SimpState)-\LI{\yBasis}(\SimpState) )
	\nonumber\\&
	\leq
	\max \Big\{ 
	-
	\zzzLM{\yBasis}
	\min_{\tempyBasis \in \DomY}
	\frac{\partial \NLFunc( \tempyBasis )}{\partial \yBasis} 
	,    
	\zzzUM{\yBasis}  
	\max_{\tempyBasis \in \DomY}
	\frac{\partial \NLFunc( \tempyBasis )}{\partial \yBasis} 
	\Big\}
	\nonumber\\&
	\leq
	\max \Big\{ 
	-
	\zzzLM{\yBasis}
	{\dinf{\NLFunc}} 
	,    
	\zzzUM{\yBasis}  
	{\dsup{\NLFunc}} 
	\Big\}
	.\label{eq:ineq_margins_NLFunc_partB}
	\end{align}
Using \eqref{eq:ineq_margins_NLFunc_partA} and \eqref{eq:ineq_margins_NLFunc_partB} yields ${\zzzUM{\Basis}}$ in \eqref{eq:UM_Basis_via_NLFunc}.
In the same manner, ${\zzzLM{\Basis}}$ in \eqref{eq:LM_Basis_via_NLFunc} is derived.

	In addition, $\zzzLM{\yBasis}$ and $\zzzUM{\yBasis}$ are $\MyOrder(\NewDist^{2})$, ${\zzzLM{\Basis}}$ and  ${\zzzUM{\Basis}}$ are also $\MyOrder(\NewDist^{2})$ because ${\zzzLipConst{\yBasis}}$ is bounded as described in Remark \ref{rem:collection}
	and because $|{\ddinf{\NLFunc}}|+|{\ddsup{\NLFunc}}| < \Globalddsup$ holds regardless of $\allDomSimpX$.
	This completes the proof.

\section{Proof of Theorem \ref{thm:ex_GPsd}} \label{pf:ex_GPsd}

Because of $\kernel{ \wildcard }{\DState{\IDdat}} \in {\FSETfunc{\NumBasis}{1}}$ in Theorem \ref{thm:ex_SEkernel},
we choose  
$[ {\subAbasis{1}}, \dots, {\subAbasis{\NumData}}]^{\MyTRANSPO}=\kernelVec$,
${\subBbasis{\IDdat}}={\El{-\kernelMat^{-1}}{\IDEl,\IDdat}}$
for $\IDdat \leq \NumData$,
${\subAbasis{\IDdat}}={\subBbasis{\IDdat}}=0$ for $\IDdat>\NumData$, 
and $\NLFunc(\yBasis)=\yBasis$.	
Then, ${\El{-\kernelMat^{-1}\kernelVec(\State)}{\IDEl}}
=\NLFunc(   
\sum_{\IDdat=1}^{\NumBasis} 
{\subAbasis{\IDdat}}(\State) {\subBbasis{\IDdat}}(\State)    
)$ and thus ${\El{-\kernelMat^{-1}\kernelVec}{\IDEl}} \in  {\FSETfunc{\NumBasis}{2}}$ holds.
We next choose
$[ {\subAbasis{1}}, \dots, {\subAbasis{\NumData}}]^{\MyTRANSPO}=\kernelVec$,
$ {\subAbasis{\NumData+1}} = \kernel{\State}{\State} $,
$[ {\subBbasis{1}}, \dots, {\subBbasis{\NumData}}]^{\MyTRANSPO}={-\kernelMat^{-1}\kernelVec} $,
$ {\subBbasis{\NumData+1}} = 1 $, 
${\subAbasis{\IDdat}}={\subBbasis{\IDdat}}=0$ for $\IDdat>\NumData+1$,
and $\NLFunc(\yBasis)=\yBasis$.	
Then, $\yBasis(\State)
=\sum_{\IDdat=1}^{\NumBasis} 
{\subAbasis{\IDdat}}(\State) {\subBbasis{\IDdat}}(\State)    
=\kernel{\State}{\State} - \kernelVec(\State)^{\MyTRANSPO} \kernelMat^{-1} \kernelVec(\State)
=
{\El{\GPsd(\State)}{\IDEl}}^{2}$ holds
and this $\yBasis$ is contained in ${\pFSETfunc{\NumBasis}{3}}$ because of $\NLFunc \in {\pBFSETNLFunc{\yBasis}}$ and $\yBasis(\State)  \geq \LBGPsd^{2}\geq 0$.
Finally, we choose 
$ {\subAbasis{1}} ={\El{\GPsd}{\IDEl}}^{2}$,
$ {\subBbasis{1}} = 1 $, 
${\subAbasis{\IDdat}}={\subBbasis{\IDdat}}=0$ for $\IDdat>1$,
and $\NLFunc(\yBasis)=\yBasis^{1/2}$.
Bounds of $\NLFunc$ can be given by
${\dsup{\NLFunc}}
=(1/2)(\min_{\tempyBasis \in \DomLIY} \tempyBasis)^{-1/2}
$
and
${\ddinf{\NLFunc}}
=(-1/4)(\min_{\tempyBasis \in \DomLIY} \tempyBasis)^{-3/2}
$ that satisfy
$|{\dsup{\NLFunc}}|+|{\ddinf{\NLFunc}}|< \Globalddsup:=1+\max\{  \LBGPsd^{-1/2},  \LBGPsd^{-3/2} \}$.
Therefore, 
${\El{\GPsd(\State)}{\IDEl}}
=\NLFunc(\yBasis(\State))$ and ${\El{\GPsd}{\IDEl}} \in {\pFSETfunc{\NumBasis}{4}}$ holds because of $\NLFunc \in {\pAFSETNLFunc{\yBasis}}$.
This completes the proof.		

\section{Proof of Proposition \ref{thm:pFSETfunc_properties}} \label{pf:pFSETfunc_properties}	

Note that any $C^{1}$ continuous function on $\DomBigX$ is Lipschitz continuous on the bounded closed set $\DomX$.
If ${\subAbasis{\IDEl}}$ and ${\subBbasis{\IDEl}}$ are Lipschitz continuous on  $\DomX$, 
$| {\subAbasis{\IDEl}}(\State){\subBbasis{\IDEl}}(\State)
- {\subAbasis{\IDEl}}(\tempState){\subBbasis{\IDEl}}(\tempState) |
\leq
| {\subAbasis{\IDEl}}(\State) |
|{\subBbasis{\IDEl}}(\State)- {\subBbasis{\IDEl}}(\tempState) |
+
| {\subBbasis{\IDEl}}(\tempState) |
| {\subAbasis{\IDEl}}(\State)- {\subAbasis{\IDEl}}(\tempState)|
$ holds that implies the Lipschitz continuity.
If ${\yBasis}$ and $\NLFunc$ are Lipschitz continuous,
the inequality
$| \NLFunc(\yBasis(\State)) - \NLFunc(\yBasis(\tempState)) |
\leq
\LipsNLFunc
| \yBasis(\State) - \yBasis(\tempState) |
$ implies the Lipschitz continuity.

Using these properties, we prove the statement \ref{thm:Lipschitz_continuous} based on mathematical induction.
Supposing that every $\Basis \in {\pFSETfunc{\NumBasis}{\NumNest}}$ is Lipschitz continuous with an $\NumNest$,
every $\Basis \in {\pAFSETfunc{\NumBasis}{\NumNest+1}}$ is Lipschitz continuous because the corresponding ${\subAbasis{\IDEl}}$, ${\subBbasis{\IDEl}}$, and $\NLFunc$ are Lipschitz continuous.
Every $\Basis \in {\FSETfunc{\NumBasis}{\NumNest}}$ is $C^{2}$ continuous and thus Lipschitz continuous owing to Proposition \ref{thm:FSETfunc_properties} \ref{thm:C2continuous}.
Thus, any $\Basis \in {\pBFSETfunc{\NumBasis}{\NumNest+1}}$ is Lipschitz continuous because $\NLFunc\in {\pBFSETNLFunc{\yBasis}}$	is Lipschitz continuous.
Therefore, any $\Basis \in {\pFSETfunc{\NumBasis}{\NumNest+1}}$  is Lipschitz continuous.
This is satisfied for any $\NumNest$ based on mathematical induction.
In the same manner, the continuity on $\DomBigX$ is proved.

Next, we prove the statement \ref{thm:pFSETfunc_nondec}.
The properties ${\pFSETfunc{\NumBasis}{\NumNest}} \subseteq {\pFSETfunc{\NumBasis}{\NumNest+1}}$ and ${\FSETfunc{\NumBasis}{\NumNest}}  \cap \NonNegativeFSETNLFunc  \subseteq {\pBFSETfunc{\NumBasis}{\NumNest+1}}$ are proved in a similar manner to Proposition \ref{thm:FSETfunc_properties} \ref{thm:FSETfunc_nondec}.
Thus, we obtain
${\FSETfunc{\NumBasis}{\NumNest}} \cap \NonNegativeFSETNLFunc
\subseteq
{\pBFSETfunc{\NumBasis}{\NumNest+1}}
\subseteq
{\pFSETfunc{\NumBasis}{\NumNest+1}}$.
This completes the proof.

\section{Proof of Lemma \ref{thm:margins_for_inner_products_specific}} \label{pf:margins_for_inner_products_specific}

Because ${\subAbasis{\IDEl}}$ and ${\subBbasis{\IDEl}}$ are nonnegative,
the upper bound $(
{\zzzBUB{\subAbasis{\IDEl}}{\subBbasis{\IDEl}}}
+
{\zzzBUB{\subBbasis{\IDEl}}{\subAbasis{\IDEl}}}
)$ in \eqref{eq:bounds_inequality_with_coef_result1} is replaced with 
$(
{\zzzUM{\subAbasis{\IDEl}}} {\zzzUB{\subBbasis{\IDEl}}} 
+
{\zzzUM{\subBbasis{\IDEl}}} {\zzzUB{\subAbasis{\IDEl}}}
)$.
In addition, as described in Remark \ref{rem:collection}, ${\zzzLipConst{\subAbasis{\IDEl}}}$ and ${\zzzLipConst{\subBbasis{\IDEl}}}$ are bounded regardless of $\NewDist$ because Proposition \ref{thm:pFSETfunc_properties} \ref{thm:Lipschitz_continuous} leads to the Lipschitz continuity of ${\subAbasis{\IDEl}}$ and ${\subBbasis{\IDEl}}$.
Thus, the statements are proved in a manner similar to the proof of Lemma \ref{thm:margins_for_inner_products} in Appendix \ref{pf:margins_for_inner_products}.		

\section{Proof of Lemma \ref{thm:margins_for_nonlinear_mappings_specific}} \label{pf:margins_for_nonlinear_mappings_specific}

We prove the statements (i) and (ii) based on the proof of Lemma \ref{thm:margins_for_nonlinear_mappings}  in Appendix \ref{pf:margins_for_nonlinear_mappings}.
Firstly, we prove the statement (i).
Using the nondecreasing property of $\NLFunc \in {\pAFSETNLFunc{\yBasis}}$, we obtain ${\dinf{\NLFunc}}\geq 0$.
By this nonnegativity, we can consider only the case of $( {\yBasis}(\SimpState)-\LI{\yBasis}(\SimpState) )\geq 0$ in \eqref{eq:ineq_margins_NLFunc_partB}.
This implies $\MeanValueThmY\geq \min_{\tempyBasis \in \DomLIY} \tempyBasis $ and thus ${\dsup{\NLFunc}}$ in the condition {\ASSNLFC} is applicable for \eqref{eq:ineq_margins_NLFunc_partB}.
Note that the $C^{2}$ continuity of $\yBasis$ on $\DomBigX$ in  the proof of Lemma \ref{thm:margins_for_nonlinear_mappings} is used to bounding  ${\zzzLipConst{\yBasis}}$.
This $C^{2}$ continuity is replaced with the Lipschitz continuity on $\DomX$ in this proof.
Consequently, substituting these results into \eqref{eq:UM_Basis_via_NLFunc} yields \eqref{eq:UM_nondec_NLFunc}.

Next, we prove the statement (ii).
Because of the convexity of $\NLFunc  \in {\pBFSETNLFunc{\yBasis}}$,
we employ Jensen's inequality \cite[Section 3.1.8]{Boyd04} instead of \eqref{eq:ineq_margins_NLFunc_partA}, yielding
$\NLFunc( \LI{\yBasis}(\SimpState) )
=
\NLFunc (
\sum_{\IDPt=1}^{\NumSimp}\El{\SimpVecCoef}{\IDPt} \yBasis(\SimpPt{\IDdom}{\IDPt}) 
)
\leq 
\sum_{\IDPt=1}^{\NumSimp}\El{\SimpVecCoef}{\IDPt}
\NLFunc (
\yBasis(\SimpPt{\IDdom}{\IDPt}) 
)
=
\LI{(\NLFunc \circ \yBasis)} (\SimpState)$.
Using the Lipschitz continuity replaces the mean value theorem in \eqref{eq:ineq_margins_NLFunc_partB} with the inequality: 
$
\NLFunc \circ  {\yBasis}(\SimpState) 
-
\NLFunc \circ  \LI{\yBasis}  (\SimpState) 
=
\LipsNLFunc
| {\yBasis}(\SimpState)-\LI{\yBasis}(\SimpState) |
\leq
\LipsNLFunc
\max \{ 
\zzzLM{\yBasis}
,    
\zzzUM{\yBasis}  
\}
$.
Combining these results gives \eqref{eq:UM_convex_NLFunc}.

Finally, the second-order properties in the statements (i) and (ii) are proved in a manner similar to the proof of Lemma \ref{thm:margins_for_nonlinear_mappings}.
This completes the proof.

\section{Proof of Theorem \ref{thm:HdotV_margin}} \label{pf:HdotV_margin}

Because of \eqref{eq:def_HdotV}, \eqref{eq:sum_margins_HdotV} holds clearly.
Note that ${\El{\partial\MValF/\partial\State}{\IDEl}}$ and ${\El{\Mmean(\wildcard,\Input(\wildcard))}{\IDEl}}$ are contained in ${\FSETfunc{\NumBasis}{\max\{{\fNumNest{\partial\MValF}},\fNumNest{\Mmean}\}}}$.
By choosing  
$[ {\subAbasis{1}}, \dots, {\subAbasis{\DimX}}]^{\MyTRANSPO}=\partial\MValF/\partial\State$,
$[ {\subBbasis{1}}, \dots, {\subBbasis{\DimX}}]^{\MyTRANSPO}=\Mmean(\wildcard,\Input(\wildcard))$, 
${\subAbasis{\IDEl}}={\subBbasis{\IDEl}}=0$ for $\IDEl>\DimX$, and
$\NLFunc(\yBasis)=\yBasis$,
we have $\MEdotV(\State)
=\NLFunc(   
\sum_{\IDEl=1}^{\NumBasis} 
{\subAbasis{\IDEl}}(\State) {\subBbasis{\IDEl}}(\State)    
)$, indicating $\MEdotV \in  {\FSETfunc{\NumBasis}{\fNumNest{\MEdotV}}}$.

Next, for each $\IDEl$, we choose  
${\subAbasis{1}}={\El{\partial\MValF/\partial\State}{\IDEl}}$,
${\subBbasis{1}}=1$,
${\subAbasis{\IDbEl}}={\subBbasis{\IDbEl}}=0$ for $\IDbEl>1$, and
$\NLFunc(\yBasis)=|\yBasis|$.	
We have $|{\El{   \partial \MValF(\State)/\partial\State       }{\IDEl}}|
=\NLFunc(   
\sum_{\IDbEl=1}^{\NumBasis} 
{\subAbasis{\IDbEl}}(\State) {\subBbasis{\IDbEl}}(\State)    
)$, indicating $|{\El{\partial\MValF/\partial\State}{\IDEl}}| \in {\pFSETfunc{\NumBasis}{{\fNumNest{\partial\MValF}}+1}}$.
Since $\Msd(\wildcard,\Input(\wildcard))$ is nonnegative, 
Proposition \ref{thm:pFSETfunc_properties} \ref{thm:pFSETfunc_nondec} indicates
${\El{\Msd(\wildcard,\Input(\wildcard))}{\IDEl}} \in  {\pFSETfunc{\NumBasis}{\fNumNest{\Msd}+1}}$.
Thus, $|{\El{\partial\MValF/\partial\State}{\IDEl}}|$ and ${\El{\Msd(\wildcard,\Input(\wildcard))}{\IDEl}}$ are contained in ${\pFSETfunc{\NumBasis}{\max\{{\fNumNest{\partial\MValF}},\fNumNest{\Msd}\}+1}}$.
By choosing  
${\subAbasis{\IDEl}}=|{\El{\partial\MValF/\partial\State}{\IDEl}}|$,
${\subBbasis{\IDEl}}={\El{\Msd(\wildcard,\Input(\wildcard))}{\IDEl}}$ for $\IDEl\leq\DimX$,
${\subAbasis{\IDEl}}={\subBbasis{\IDEl}}=0$ for $\IDEl>\DimX$, and
$\NLFunc(\yBasis)=\yBasis$,	
we have $\SDdotV(\State)
=\NLFunc(   
\sum_{\IDEl=1}^{\NumBasis} 
{\subAbasis{\IDEl}}(\State) {\subBbasis{\IDEl}}(\State)    
)$, indicating $\SDdotV \in  {\pFSETfunc{\NumBasis}{\fNumNest{\SDdotV}}}$.
This completes the proof.

\section{Proof of Theorem \ref{thm:general_control_design}} \label{pf:general_control_design}

For every $(\HJBweight, \OffSetpenaltyFunc,  \estMDomStable, \allDomSimpX)$, there exists $\optVCparam$ because $\ObjectiveF(\wildcard;\HJBweight, \estMDomStable, \allDomSimpX)$ is continuous on the bounded closed $\SetVCparam$.
For brevity of notation, we denote ${\PenaltyperX{\VCparam}{\SimpPt{\IDdom}{\IDPt}}{\VUBofMargin}{\WUBofMargin}{\OffSetpenaltyFunc}}$, $\ObjectiveF(\VCparam;\HJBweight, \estMDomStable, \allDomSimpX)$, and $\PerformanceTerm(\VCparam; \estMDomStable, \allDomSimpX)$ by ${\SimplePenaltyperX{\VCparam}{\SimpPt{\IDdom}{\IDPt}}}$, $\ObjectiveF(\VCparam;\HJBweight )$, and $\PerformanceTerm(\VCparam )$, respectively.
Firstly, we prove the statement (i).	
Since there exists $\VCparam$ satisfies \eqref{eq:estSR_is_a_SR_offsetFuncs} by the assumption {\ASSVCparamExist}, 
for every $\NewDist$, 
for every ${\SimpPt{\IDdom}{\IDPt}} \in  \estMDomStable$,
we have
${\SimplePenaltyperX{\VCparam}{\SimpPt{\IDdom}{\IDPt}}}
\leq 2\OffSetpenaltyFunc
$.
Using such a $\VCparam$ and the optimality of $\stabVCparam$ yields
$
\ObjectiveF(\stabVCparam; 0)
\leq
\ObjectiveF(\VCparam; 0)
=
\sum_{ \State \in  \VerticesestMDomStable }
{\SimplePenaltyperX{\VCparam}{\State}}
\leq
\sum_{ \State \in  \VerticesestMDomStable }
2 \OffSetpenaltyFunc 
$,
which implies 
$\ObjectiveF(\stabVCparam; \HJBweight)
\leq
\HJBweight \PerformanceTerm(\stabVCparam) 
+ \sum_{\State \in  \VerticesestMDomStable } 2\OffSetpenaltyFunc $.
By combining this result with the relations 
$\ObjectiveF(\optVCparam; \HJBweight) \leq \ObjectiveF(\stabVCparam; \HJBweight)$, 
$\PerformanceTerm(\optVCparam) \geq 0$, and
$\penaltyFunc(\ArgpenaltyFunc ;\OffSetpenaltyFunc)\geq 0$,
for every $(\HJBweight,\estMDomStable,\allDomSimpX,\OffSetpenaltyFunc)$, we obtain
\begin{align}
\max_{ \State \in  \VerticesestMDomStable }
{\SimplePenaltyperX{\optVCparam}{\State}}
\leq
\ObjectiveF(\optVCparam; \HJBweight) 
&\leq 
\ObjectiveF(\stabVCparam;\HJBweight)
\leq 
\HJBweight \PerformanceTerm(\stabVCparam) 
+ 
\sum_{ \State \in  \VerticesestMDomStable }
2
\OffSetpenaltyFunc
.\label{eq:pf_max_penalty_point}
\end{align}	
Because of $\estMDomStable \subseteq \DomX$ and the continuity of $\VUBofMargin$ and $\WUBofMargin$ on the bounded closed $\DomX$, there exists a positive lower bound 
$\LBUBofMargin := \min\{
\inf_{\State \in \estMDomStable} \VUBofMargin(\State)
,
\inf_{\State \in \estMDomStable} \WUBofMargin(\State)	 
\}/2>0$ independent of $\NewDist$.
For each $\allDomSimpX$, we set $\UBHJBweight>0$ and $\UBOffSetpenaltyFunc>0$ 
such that 
$
\UBHJBweight 
(
\sup_{\OffSetpenaltyFunc \in [0,\UBOffSetpenaltyFunc]}
\PerformanceTerm(\stabVCparam) 
)
+ 
\UBOffSetpenaltyFunc
\sum_{ \State \in  \VerticesestMDomStable } 2
\leq 
\LBUBofMargin
$ holds.
Then, for any $\NewDist>0$, $\HJBweight \in [0,\UBHJBweight]$, and $\OffSetpenaltyFunc \in [0,\UBOffSetpenaltyFunc]$, 
using \eqref{eq:pf_max_penalty_point} gives
\begin{align}
&
\max_{ \State \in  \VerticesestMDomStable }
\max\{
\penaltyFunc(\VUBofMargin(\State)  -\MValF(\State ; \optVCparam) ;\OffSetpenaltyFunc)
,
\penaltyFunc(\WUBofMargin(\State)  +\HdotV(\State ; \optVCparam) ; \OffSetpenaltyFunc)
\}
\nonumber\\&
\leq
\max_{ \State \in  \VerticesestMDomStable }
{\SimplePenaltyperX{\optVCparam}{\State}}
\leq
\LBUBofMargin
\leq
\penaltyFunc(\LBUBofMargin;\OffSetpenaltyFunc) 
.
\end{align}	
This indicates that for any $\IDdom \in \estSetStableIDdom $ and $\IDPt$, we obtain 
$\MValF(\SimpPt{\IDdom}{\IDPt}; \optVCparam) \geq \VUBofMargin(\SimpPt{\IDdom}{\IDPt}) - \LBUBofMargin \geq \LBUBofMargin$
and
$\HdotV(\SimpPt{\IDdom}{\IDPt}; \optVCparam)  \leq -\WUBofMargin(\SimpPt{\IDdom}{\IDPt}) + \LBUBofMargin \leq -\LBUBofMargin$.
Meanwhile,
because of the assumption \ASSboundedMargin,	
if $ \UBDist < \min \{ \ODEFUBDist, \sqrt{  {\LBUBofMargin}/ \max\{{\LUMconst{\MValF}},  {\LUMconst{\HdotV}}  \}   } \}$ holds, 
any $\NewDist \in (0, \UBDist]$ satisfies
$\zzzLM{\MValF;\optVCparam} \leq  {\LUMconst{\MValF}} \NewDist^{2} < {\LBUBofMargin}$ and 
$\zzzUM{\HdotV;\optVCparam} \leq  {\LUMconst{\HdotV}} \NewDist^{2} < {\LBUBofMargin}$.
Therefore, Algorithm {\MyLabelAlgSR} gives $\MDomStable$ satisfying
$\estMDomStable \subseteq \MDomStable$.
This implies the statement (i).

Next, we prove the statement (ii).
Supposing that \eqref{eq:Design_suboptimality} does not hold,
using the optimality $\ObjectiveF(\optVCparam; \HJBweight)\leq  \ObjectiveF(\stabVCparam; \HJBweight)$ gives the following inequality:
$
\ObjectiveF(\optVCparam;0)
=
\ObjectiveF(\optVCparam; \HJBweight) - \HJBweight\PerformanceTerm(\optVCparam)
<  
\ObjectiveF(\stabVCparam; \HJBweight) - \HJBweight\PerformanceTerm(\stabVCparam) 	
=
\ObjectiveF(\stabVCparam; 0)
$.
This contradicts the optimality $\ObjectiveF(\stabVCparam; 0)
\leq \ObjectiveF(\optVCparam; 0)$.
This completes the proof.

\end{document}